\documentclass[a4paper, 11pt]{amsart}   
\usepackage{comment}   
\usepackage{color}

\newtheorem{theorem}{Theorem}[section]

\newtheorem{lemma}[theorem]{Lemma}
\newtheorem{proposition}[theorem]{Proposition}

\newtheorem{definition}[theorem]{Definition}
\newtheorem{definition and theorem}[theorem]{Definition and
Theorem}
\newtheorem{remark}[theorem]{Remark}
\newtheorem{remarks}[theorem]{Remarks}
\newtheorem{*remark}[theorem]{$^* $Remark}

\newtheorem{*exercise}[theorem]{$^* $Exercise}
\newtheorem{**exercise}[theorem]{$^{** } $Exercise}
\newtheorem{example}[theorem]{Example}
\newtheorem{examples}[theorem]{Examples}

\newtheorem{corollary}[theorem]{Corollary}

\excludecomment{versionA}   
   
\def\rouge{}
\def\bleu{}
\def\noir{}   
   
\usepackage[a4paper]{geometry}   
\title[Convergence Stationary RBF schemes]{Convergence of Stationary Radial Basis Function-schemes for 
evolution equations}      
   
\author[Brad Baxter and Raymond Brummelhuis]{Brad Baxter\address{Brad Baxter, Department of Economics, Mathematics and Statistics, Birkbeck, University of London, Malet Street, London WC1E 7HX, UK, \lowercase{e-mail: b.baxter@bbk.ac.uk}}  and Raymond Brummelhuis\address{Raymond Brummelhuis, Laboratoire de Math\'ematiques de Reims, Universit\'e de Reims-Champagne-Ardenne, FR 3399 CNRS, BP 1039, 51687 cedex 2, Reims, France, \lowercase{e-mail: raymondus.brummelhuis@univ-reims.fr}}     
}   
      
\begin{document}   
   
\begin{abstract}We establish precise convergence rates for semi-discrete schemes 
based on Radial Basis Function interpolation, as well as approximate approximation results for such schemes. Our schemes use stationary interpolation on regular grids, with basis functions from a general class of functions generalizing one introduced earlier by M. Buhmann. Our results apply to parabolic equations such as the heat equation or Kolmogorov-Fokker-Planck equations associated to L\'evy processes, but also to certain hyperbolic equations.   
\end{abstract}   
   
\subjclass[2010]{Primary: 65M12, 65M15, 65M20, 65M70; Secondary: 35S05, 35S10, 41A05, 41A25, 60G51}   

\maketitle   
   
   
\section{\bf Introduction}   
   
We establish convergence rates for semi-discrete numerical schemes based on stationary Radial Basis Function (RBF) interpolation. We do this 
for the classical heat equation 
but also for more general constant coefficient (translation invariant) 
pseudo-differential evolution equations. 
Our analysis applies to (suitably regularized versions of) the fractional heat equation and other Kolmogorov-Fokker-Planck equations associated to L\'evy processes, 
but also to certain non-parabolic equations such as the free Schr\"odinger equation or  hyperbolic equations such as the half-wave and transport equations. 
   

The scheme we consider is an RBF-version of the method of lines, implemented on regular square grids, and we examine the convergence of the scheme when the grid size tends to 0, using basis functions which scale with the grid. We refer to Buhmann \cite{Bu_book}, Wendland \cite{W} and Fasshauer \cite{F} for general introductions to Radial Basis Function interpolation, and to the beginning of sections 5 and 6 below for a precise description of the scheme; roughly speaking, RBF-interpolation seeks to interpolate a function by a linear combination of translates of a given basis function, centered in the interpolation points. The coefficients can be computed by solving a linear system 
whose coeffient matrix, under appropriate conditions on the basis function, is always non-singular. To turn this into a scheme for solving a linear evolution equation on $\mathbb{R }^n $, one takes the coefficients to be time-dependent and asks that the equation be satisfied exactly in the interpolation points, leading to a system of ordinary differential equations. An advantage of using a regular grid for the interpolation points is that the interpolation problem can be reduced to the construction of a single Lagrange function (also called cardinal function), which greatly simplifies the theoretical analysis. The advantage of using basis functions which scale with the grid, known as stationary interpolation, is that the condition number of the coefficient matrix of the linear system which, in a practical implementation of the scheme, has to be solved numerically, becomes independent of the grid-size.  
   
The numerical performance of such a scheme 
was examined in 
\cite{RB_RC}, \cite{KVF}, \cite{PLMP} \cite{RC}, \cite{RC_SH}, motivated by applications in mathematical finance.   
Our main results will relate the order of convergence of the scheme to the degree of the operator and to the approximation order of the underlying RBF interpolation, as function of the grid-size. The latter will only be algebraic, since we use stationary interpolation, but can be arbitrarily large, depending on the basis function and 
the degree of smoothness of the initial condition. We will furthermore show that under certain circumstances we can have {\it approximate approximation}, in the sense that, in case of non-convergence of the scheme to the true solution, one can nevertheless get arbitrarily close to the real solution by an appropriate choice of basis function. More generally, if the scheme does converge, one can, for initial values which are sufficiently smooth and for grid-sizes which are small but not too small, observe an apparent order of convergence which is bigger than the actual one.  This was observed numerically in \cite{RB_RC}. 
The notion of approximate approximation was introduced by Maz'ya \cite{M_94}, and analyzed in detail for approximation using Gaussian kernels in Maz'ya and Schmidt \cite{M_S_96}, \cite{M_S}. Amongst  other things, we generalize their work to more general basis functions.   
\medskip   
   
\begin{versionA}   
\noindent \textcolor{blue}{Since we use interpolation on Cartesian grids $\Rightarrow $ can use Lagrange functions which simplify; no need to discuss/phrase things in terms of conditional positive definiteness and  adding polynomials to the translates of the RBF, although practical implementation would need to be done in this way: see \cite{RB_RC}, Appendix A, for formulation}   
\end{versionA}   
\medskip   
   
We will in fact not 
restrict ourselves to particular examples of Radial Basis Functions, such as the generalized multiquadrics or the polyharmonic basis functions, but carry out our analysis for a general class of basis functions 
which we 
introduce in section 2. Since this class is a generalization of one introduced by Martin Buhmann in \cite{Bu1} and by which it was inspired, we have called it the Buhmann class. We will analyze the properties of our scheme in Fourier space and, for that reason, first re-examine in section 3 the convergence of RBF-interpolation on regular grids from the Fourier point of view by deriving precise estimates for the Wiener norm of the difference between a function and its RBF-interpolant. Our convergence theorems have a non-zero intersection with classical results of Buhman and Powell (see \cite{Bu_book} and its references), 
generalizing these in some respects. 
Despite the use of the Wiener norm we can allow certain classes of polynomially increasing functions. The Fourier transform of such a function will have a non-integrable algebraic singularity at 0, and the allowed order of the singularity (and therefore the allowed 
growth of the function) will depend on the basis function which is used for the interpolation. In section 4 we show that the convergence rates which we found in section 3 are best possible, and discuss approximate approximation.   
   
The next two sections examine the convergence of the RBF-variant of the method of lines, first, in section 5, for the in many respects typical case of classical heat equation before indicating, in section 6, how these results extend 
to more 
general pseudo-differential evolution equations. We show that the scheme converges at a rate of $h^{\kappa - q } $, where $q $ is the order of the operator ($q = 2 $ for the heat equation) and $\kappa $ the order of convergence of the underlying RBF-interpolation scheme, which is also the order of the singularity in 0 of the Fourier transform of the basis function which is used. We show that this rate is in general optimal, 
but that we can also have 
approximate approximation, in the sense that for appropriate basis functions which are sufficiently "flat" and with sufficiently smooth initial data there can be an apparent higher order of convergence when $h $ is "small but not too small", which is determined by the degree of smoothness of the initial data. This 
provides an explanation for the empirical convergence rates which were 
observed in \cite{RB_RC}.   

\medskip   
      
One obvious limitation of the present paper is that we have restricted ourselves to interpolation on regular grids of scaled integer points, whereas one of the strengths of the RBF method is that one can use arbitrarily scattered interpolation points, opening up the way to adaptive methods. This flexibility may become important when dealing with variable coefficient linear differential operators or with non-linear ones.  We note however, in our defense, that the often-used Finite Difference methods are usually restricted to regular grids also, and that even on regular grids RBF methods 
can have definite advantages over FD methods when treating non-local operators, as they do not 
discretize the operator, and can therefore be better suited when the latter has a singular kernel, such as for the Kolmogorov backward equation of certain L\'evy processes: see \cite{RB_RC}. Another limitation is that we only have treated translation invariant evolution equations. 
These do however already include large classes of operators which are of interest for applications, such as 
the fractional heat equation or other Kolmogorov-Fokker-Planck equations associated to L\'evy processes. 
It would obviously be 
of interest to generalize our results to variable coefficient PDEs, but this will require other methods. We also have (primarily) examined convergence in Wiener norm: 
convergence of the scheme in the $L^2 $ norm or more general Sobolev norms will be treated elsewhere. We finally 
want to note that although we have used stationary RBF interpolation, our analysis 
can be extended to the non-stationary case, when one uses the same basis function for all grid-sizes. In this case one 
can have exponential rates of convergence, as first discoverd by Madych and Nelson \cite{M_N} for interpolation 
of functions in the so-called native space of the basis function. This also will be examined elsewhere. 
\medskip   

\noindent {\bf Notations
}: $C $ denotes the usual "variable constant", whose exact numerical value is allowed to change from one occurrence to the other.   
We use the following convention for the Fourier transform $\widehat{f } = \mathcal{F }(f) $ of an integrable function $f $ on $\mathbb{R }^n $:          
$$   
\widehat{f }(\xi ) = \int _{\mathbb{R }^n } f(x) e^{-i (x, \xi ) } dx ,   
$$   
 $(x, \xi ) $ being the Euclidean inner product on $\mathbb{R }^n . $ We will routinely use the extension of the Fourier transform $\mathcal{F } $ to the space of tempered distributions $\mathcal{S }' (\mathbb{R }^n ) $, where $\mathcal{S }(\mathbb{R }^n ) $ is the usual Schwarz-space of rapidly decreasing functions.   
      
For $s \in \mathbb{R } $, let $\mathring{L }^1 _s (\mathbb{R }^n ) $ be the space of measurable functions on $\mathbb{R }^n $ for which     
\begin{equation} \label{eq:L^1_s}   
|| f ||^{\circ } _{1, s } := \int _{\mathbb{R }^n } |\xi |^s |f(\xi ) | \, d\xi < \infty ,      
\end{equation}   
and $L^1 _{\kappa } (\mathbb{R }^n ) := L^1 (\mathbb{R }^n ) \cap \mathring{L }^1 _{\kappa }  (\mathbb{R }^n ) $, its non-homogeneous version, with norm $\int |f | (1 + \xi | )^s d\xi . $   
We will also need the "mixed" spaces   
\begin{equation} \label{eq:L^1_rs}   
\mathring{L }^1 _{r, s } := \{ f \in \mathring{L }^1 _s (\mathbb{R }^n ) : || f ||^{\circ } _{r, s } := \int _{\mathbb{R }^n } | f(\xi ) | \, ( |\xi |^r \wedge |\xi |^s ) d\xi < \infty \} ,   
\end{equation}   
for $r \leq s $, where $a \wedge b := \min (a, b ) . $ Note that this scale of spaces is increasing in $r $ for $r \leq s $, and that $\mathring{L }^1 _r (\mathbb{R }^n ) \subset \mathring{L }^1 _{r, s } (\mathbb{R }^n ) $, while $\mathring{L }^1 _{s,s } (\mathbb{R }^n ) = \mathring{L }^1 _s (\mathbb{R }^n ) . $   
\medskip   
   
\begin{versionA}   
\noindent \textcolor{blue}{Parenthetical remark for self: concretely, since $r \leq s $ by assumption,    
$$   
|| f ||^{\circ } _{r, s } := \int _{|\xi | \leq 1 } |\xi |^s  | f(\xi ) | d\xi + \int _{|\xi | \geq 1 } |\xi |^r | f(\xi ) | d\xi ,   
$$   
which makes it obvious that the norm increases in $r $ for $r \in [0 , s ] $ and also that   
$$   
|| f ||^{\circ } _{r, s } \leq || f ||^{\circ } _r .   
$$   
This norm should be contrasted with   
$$   
\int _{\mathbb{R }^n } | f(\xi ) | \, \max ( |\xi |^r , |\xi |^s ) d\xi = \int _{|\xi | \leq 1 } |\xi |^r  | f(\xi ) | d\xi + \int _{|\xi | \geq 1 } |\xi |^s | f(\xi ) | d\xi ,   
$$   
which also plays a r\^ole (notably in our approximate approximation estimates), although we don't really need to introduce this norm explicitly in our estimates (and we won't in fact in this latest version)   
}   
\medskip   
   
\end{versionA}   

Finally, we define weighted sup-norm spaces $L^{\infty } _s (\mathbb{R }^n ) $ of measurable functions such that   
\begin{equation} \label{eq:L_inf_s}   
|| f ||_{\infty , s } := \sup _{x \in \mathbb{R }^n } (1 + |x| )^s |f(x) | ,   
\end{equation}   
where the sup is the essential supremum, as usual.   
If $s < 0 $, an element $f $ of $L^{\infty } _s (\mathbb{R }^n ) $ is of polynomial growth of order at most $|s | $: $|f(x) | \leq C (1 + |x|)^{|s | } $ on $\mathbb{R }^n $, with $C = || f ||_{\infty , s } . $   
   
Derivatives of functions $f = f(x) $ on $\mathbb{R }^n $ will be denoted by $\partial _x ^{\alpha } f (x) $ or by $f^{(\alpha ) } (x) $, $\alpha \in \mathbb{N }^n $ a multi-index. If $K \in \mathbb{N } $ and $\lambda \in (0, 1 ] $, then 
$C_b ^{K, \lambda } (\mathbb{R }^n ) $ will denote the H\"older space of $K $-times differentiable functions on $\mathbb{R }^n $ with bounded derivatives of all orders, such that the derivatives of order $K $ satisfy a uniform H\"older condition on $\mathbb{R }^n $ with exponent $\lambda $, 
provided with the norm   
$$   
\sum _{|\alpha | \leq K } || f^{(\alpha ) } ||_{\infty } + \sum _{|\alpha | = K } || f^{(\alpha ) } ||_{0 , \lambda } ,   
$$   
where $|| g ||_{0 ; \lambda } := \sup _{\xi \neq \eta } | g (\xi ) - g (\eta ) |  / |\xi - \eta |^{\lambda } . $   
      
Finally, $\lfloor x \rfloor $ and $\lceil x \rceil $ denote the usual floor and ceiling functions, defined as the   
greatest, respectively smallest integer which is less than,
respectively greater than a real number $x $; note that $\lceil x \rceil =
\lfloor x \rfloor + 1 $ if $x \notin \mathbb{N } $, while $\lceil
x \rceil = \lfloor x \rfloor = x $ otherwise.   
   
\section{\bf A class of basis functions for interpolation on a regular grid}

\subsection{The Buhmann class }   
We introduce a flexible class of basis functions which is well-suited for stationary interpolation on regular grids. Since this class of functions is a 
generalisation of the one introduced earlier by Buhmann \cite{Bu1} (called {\it admissible} in there) we will call it the {\it Buhmann class}. From the onset, we will allow non-radial basis functions, radiality not being essential for most of the theory (as is of course well known). 
     
\begin{definition} \label{def:wBuhmann class} \rm{For $\kappa \geq 0 $ and $N > n $ we define the Buhmann class $\mathfrak{B }_{\kappa , N } (\mathbb{R }^n ) $ as the set of   functions $\varphi \in C(\mathbb{R }^n ) $ such that   
\medskip   
   
\noindent (i) $\varphi $ is of polynomial growth of order strictly less than $\kappa $, in the sense that $\varphi \in L^{\infty } _{- \kappa + \varepsilon } (\mathbb{R }^n ) $ for some $\varepsilon > 0 . $    
\medskip   
   
\noindent (ii) ({\it Regularity and strict positivity.}) The restriction to $\mathbb{R }^n \setminus 0 $ of the Fourier transform $\widehat{\varphi } := \mathcal{F } (\varphi ) $ (in the sense of tempered distributions) can be identified with a function in $C^{n + \lfloor \kappa \rfloor + 1 } (\mathbb{R }^n \setminus 0 ) $, which we will continue to denote by $\widehat{\varphi } $, which is pointwise strictly positive: $\widehat{\varphi } (\eta ) > 0 $ for all $\eta \in  \mathbb{R }^n \setminus 0 . $    
\medskip   
   
\noindent (iii) ({\it Elliptic singularity at 0.}) There exist positive constants $c,  C $ such that for all $|\alpha | \leq n + \lfloor \kappa \rfloor + 1 $,           
\begin{equation} \label{eq:wB2}   
|\partial _{\eta } ^{\alpha } \widehat{\varphi } | \leq C |\eta |^{-\kappa - |\alpha | } , \ \ |\eta | \leq 1 ,   
\end{equation}   
while also      
\begin{equation} \label{eq:wB3}   
\widehat{\varphi } (\eta ) \geq c |\eta |^{- \kappa } , \ \ |\eta | \leq 1 .   
\end{equation}   
\medskip   
   
\noindent (iv) ({\it Decay at infinity.}) There exist positive constants $C_{\alpha } $, $|\alpha | \leq n + \lfloor \kappa \rfloor + 1 $, such that   
\begin{equation} \label{eq:wB1}   
|\partial _{\eta } ^{\alpha } \widehat{\varphi }(\eta ) | \leq C_{\alpha } |\eta | ^{-N } , \ \ |\eta | \geq 1 .   
\end{equation}   
\medskip  
   
}   
\end{definition}   

We use the term "elliptic" for condition (iii) because of the resemblance of (\ref{eq:wB2}) and (\ref{eq:wB3}) with the ellipticity condition on symbols in pseudo-differential theory (where the singularity would be at infinity). The significance of $n + \lfloor \kappa \rfloor + 1  $ is that this is the smallest integer which is strictly greater than $n + \kappa . $ (Note that if $\kappa \notin \mathbb{N } $, then $n + \lfloor \kappa \rfloor + 1 = n + \lceil \kappa \rceil . $) Conditions (ii) and (iii) for derivatives up till this order will imply polynomial decay of order $n + \kappa $ of the associated Lagrange interpolation function which we will define below. Requiring higher order differentiability would not improve this rate of decay: $n + \kappa $ is best possible, under condition (iii).   
   
In most 
of the results of this paper, strict positivity of $\widehat{\varphi } $ on $\mathbb{R }^n \setminus 0 $ could have been replaced by the weaker condition that the "periodisation" of $\widehat{\varphi } $, $\sum _k \widehat{\varphi } (\eta + 2 \pi k ) $, be pointwise strictly positive on all of $\mathbb{R }^n $, as in \cite{Bu1}; note that by (\ref{eq:wB1}) with $\alpha = 0 $, this series converges absolutely on $\mathbb{R }^n \setminus \mathbb{Z }^n $ , given that $N > n $, while it can be set equal to $\infty $ on $\mathbb{Z }^n $, in view of (\ref{eq:wB3}). Since for most of the radial basis functions used in practice, $\widehat{\varphi } (\eta ) $ itself is already strictly positive, we have opted to impose the stronger condition, which also 
simplifies the proofs.   
\medskip   
   
\begin{remarks} \label{remark:relation_Bu1} \rm{(i) Buhmann \cite{Bu1} studied stationary RBF interpolation on regular grids for a slightly more restricted class of radial basis functions. The main difference between his original class and the one of our definition \ref{def:wBuhmann class} (besides, as already mentioned, Buhmann requiring strict positivity of the periodisation of $\widehat{\varphi } $ instead of of $\widehat{\varphi } $ itself) lies in condition (iii), where Buhmann asks that for small $|\eta | $, $\widehat{\varphi }(\eta ) $ be asymptotically equivalent to a positive multiple of $|\eta |^{- \kappa } $ modulo an relative error which has to be sufficiently small: $\widehat{\varphi } (\eta) = A |\eta |^{- \kappa } (1 + h(\eta ) ) $ with $A > 0 $ and $|\partial ^{\alpha } _{\eta } h(\eta ) | = O(|\eta |^{\varepsilon - |\alpha | } ) $ as $\eta \to 0 $ for  $|\alpha | \leq n + \lfloor \kappa \rfloor + 1 $, with an $\varepsilon > \lceil \kappa \rceil - \kappa . $ 
Under these conditions Buhmann proved the existence of a unique Lagrange function for interpolation on $\mathbb{Z }^n $, constructed as an infinite linear combination of translates of $\varphi $, which moreover decays as $|x |^{- \kappa - n } $ at infinity. This fundamental result remains true for $\varphi $'s in $\mathfrak{B }_{\kappa , N } (\mathbb{R }^n ) $:  see theorem \ref{thm_L_1} below and and its proof in Appendix \ref{Appendix:Lagrange_function}. The condition that $\varepsilon > \lceil \kappa \rceil - \kappa $ 
is in our treatment made unnecessary by lemma \ref{lemma:Stein}.   
   
\medskip   
   
\noindent (ii) All conditions in definition \ref{def:wBuhmann class} except the first are on the Fourier transform of $\varphi . $ One can show (cf. Appendix \ref{Appendix:Lagrange_function}) that if the Fourier transform of a polynomially increasing function $\varphi $ satsifies (ii), (iii) and (iv), then there exists a function $\widetilde{\varphi } (x) $ which grows at most as $\max (|x |^{\kappa - n } \log |x| , 1 ) $ at infinity (and, slightly better, as $\max (|x|^{\kappa - n } , 1 ) $ if $\kappa \notin \mathbb{N } $) and a polynomial $P(x) $ such that   
$$   
\varphi (x) = \widetilde{\varphi }(x) + P(x) .   
$$   
The function $\widetilde{\varphi } $ is unique modulo polynomials of degree $\lfloor \kappa \rfloor - n . $ If we moreover require $\varphi $ to have polynomial growth of order strictly less than $\kappa $, as in definition \ref{def:wBuhmann class}, then $P(x) $ will be a polynomial of degree of at most $\lceil \kappa \rceil - 1 $ (which is $\lfloor \kappa \rfloor $ if $\kappa \notin \mathbb{N } $, and $\kappa - 1 $ if $\kappa \in \mathbb{N } $). Note that the Fourier transform of a polynomial is a linear combination of derivatives of the delta-distriubution in 0, and therefore equals 0 on $\mathbb{R }^n \setminus 0 . $   
\medskip   
   
\noindent (iii) The condition that $N > n $ will suffice for convergence of the RBF interpolants on regular grids $h \mathbb{Z }^n $ as $h \to 0 $, but will have to be strengthened to $n > N + k $ for convergence of the RBF schemes for solving parabolic PDEs and PIDEs which are of order $k $ (in the space variables).   
\medskip    
   

}   
\end{remarks}

The usual examples of radial basis functions,  such as the generalised multi-quadrics, cubic and higher order splines, thin
plate splines, inverse multi-quadrics and Gaussians, are Buhmann class. 
\medskip   
   
   
One can show that if $\varphi \in L^{\infty } _{- p } (\mathbb{R }^n ) $, $p \in \mathbb{N } $, 
satisfies conditions (ii) - (iv) of definition \ref{def:wBuhmann class}, then $\varphi $ conditionally positive definite of order $\mu $, where $\mu $ is the smallest integer such that $2 \mu > \max ( \lfloor \kappa \rfloor - n , p , 0 ) $: for this it would in fact be sufficient that $\widehat{\varphi } |_{\mathbb{R }^n \setminus 0 } $ is locally integrable, satisfies (\ref{eq:wB2}) with $\alpha = 0 $ and is integrable on $\{ |\eta | \geq 1 \} . $ One can therefore, by standard RBF theory, interpolate an arbitrary function on a finite set $X $ of points by a linear combination of translates of $\varphi $ plus a  polynomial of degree $\mu - 1 $, provided the set $X $ is unisolvent for this class of polynomials: see for example \cite{Bu_book}. This  involves solving a linear system of equations.  The next theorem establishes the existence and main properties of a Lagrange function in terms of which the solution of the interpolation problem on $\mathbb{Z }^n $ can be simply expressed. 
   
\begin{theorem} \label{thm_L_1} Suppose that $\varphi \in \mathfrak{B }_{\kappa , N } (\mathbb{R }^n ) $ with $\kappa \geq 0 $ and $N > n . $ Then there exist coefficients $c_k $, $k \in \mathbb{Z }^n $, such that the series   
\begin{equation} \label{def:L_1}
L_1 (x) := L_1 (\varphi )(x) := \sum _{k \in \mathbb{Z }^n } c_k \varphi (x - k )
\end{equation}
converges absolutely and uniformly on compacta and defines a Lagrange function for interpolation on $\mathbb{Z }^n : $   
$$
L_1 (j ) = \delta _{0j }, \ \ j \in \mathbb{Z }^n .
$$
The function $L_1 $ satisfies the bound
\begin{equation} \label{bound_L_1}
|L_1 (x) | \leq C (1 + |x| )^{- \kappa - n } , \ \ x \in \mathbb{R }^n ,   
\end{equation}
and its Fourier transform is given by   
\begin{equation} \label{Fourier_L_1}   
\widehat{L }_1 (\eta ) = \frac{\widehat{\varphi } (\eta ) }{\sum _k \widehat{\varphi } (\eta + 2 \pi k ) } .   
\end{equation}   
Moreover, at the points of $2 \pi \mathbb{Z }^n $, $\widehat{L }_1 $   
satisfies the Fix-Strang conditions:   
\begin{equation} \label{eq:SFC}   
\widehat{L }_1 (2 \pi k + \eta ) = \delta _{0k } + O(|\eta |^{\kappa } ), \ \ \eta \to 0 .   
\end{equation}   
\end{theorem}   
\medskip   
   
\noindent See Appendix A for the proof. 
We make the trivial but for the sequel important observation that   
\begin{equation} \label{eq:sum_transl_L-hat}   
\sum _{k \in \mathbb{Z }^n } \widehat{L }_1 ( \eta + 2 \pi k ) = 1 .   
\end{equation}   
\medskip   
   
\begin{remarks} \rm{(i) We will write $L_1 (\varphi ) $ if we want to stress the dependence on the basis function $\varphi $,   
otherwise we will simply write $L_1 . $ The subindex 1 in $L_1 $ is a notational reminder that $L_1 $ is a Lagrange function for interpolation on the standard grid $\mathbb{Z }^n $ with width 1. For {\it stationary RBF interpolation} on the scaled grids $h \mathbb{Z }^n $ one uses the scaled basis functions $\varphi (x/h ) $, whose associated Lagrange functions then simply are $L_h (x) := L_1 (x/h ) . $ If $f \in L^{\infty } _{- p } (\mathbb{R }^n ) $ for some $p < \kappa $, then $s_h [f ] (x) := \sum _j f(hj ) L_1 (h^{-1 } x - j ) $ is an infinite linear combination of translates of $\varphi $ which will interpolate $f $ on $h \mathbb{Z }^n $, where the series converges absolutely and uniformly on compacta, in view of the growth restriction on $f . $   
\medskip   
   
\noindent (ii) One important point of the theorem is that the basis function $\varphi $ need not
decay at infinity, but is allowed to grow polynomially. A high order of growth will in fact lead to a high order convergence of the stationary RBF interpolants $s_h [f ] $ to $f $ as $h \to 0 $, since this will translate into a strong singularity in 0 of the Fourier transform of $\widehat{\varphi } $ in the form of a large $\kappa $, which implies that $\widehat{L }_1 $ will satisfy the Fix - Strang conditions to a high order. The latter then implies a convergence rate of $O(h^{\kappa } ) $ in sup-norm, as shown by Buhmann \cite{Bu1} (under suitable conditions on $f $); see also \cite{Bu_book}, Chapter 4. We will prove such convergence theorems for the Wiener norm instead of the sup-norm, using an entirely different approach: see theorems \ref{thm:convergence_RBF}, \ref{thm:convergence_RBF_bis} and \ref{thm:conv_RBF_3}. Note that, contrary to $\varphi $, the Lagrange function $L_1 $ will decay at infinity, as shown by (\ref{bound_L_1}), and this the more rapidly the higher $\kappa $ is. In particular, $L_1 $ is integrable if $\kappa > 0 $ and its Fourier transform then exists in the classical sense, as an
absolutely convergent integral. 
It is possible for $L_1 (x) $ to have faster decay: Buhmann \cite{Bu1} shows that if $\widehat{\varphi } (\eta ) \sim |\eta |^{- \kappa } $ as $\eta \to 0 $ with 
$\kappa \in 2 \mathbb{N } $, then   
$$   
L_1 (x) \leq C (1 + |x| )^{- \kappa - n - \varepsilon} ,   
$$   
while there are examples 
of $\varphi $ for which $L_1 (x) $ decays exponentially: see \cite{Bu_book} for details and references.      
\medskip

\noindent (iii) The proof of theorem \ref{thm_L_1} shows that the coefficients $c_{-k } $ are precisely the Fourier
coefficients of $(\sum _k \widehat{\varphi } (\eta + 2 \pi k )^{-1
} . $ They satisfy bounds analogous to the ones satisfied by $L_1 $:
$|c_k| = O(|k |^{- n - \kappa } ) . $   
This guarantees that the defining series for $L_1 (x) $ converges
absolutely and uniformly on compacta including when $\kappa = 0 $, in view of 
condition (i) of definition \ref{def:wBuhmann class}.   
}   
\end{remarks}   
   
Since the denominator of (\ref{Fourier_L_1}) is $2 \pi $-periodic and
positively bounded away from 0, $\widehat{L }_1 (\eta ) $ will have the
same decay as $\widehat{\varphi } (\eta ) $ as $|\eta | \to
\infty . $ We state this as a lemma, for later reference:   
   
\begin{lemma} \label{lemma: decay_L1_hat} There exists a constant $C = C_n > 0 $ such that for all $\ell \in \mathbb{Z }^n \setminus 0 $,   
\begin{equation} \label{eq:SFC1}
\max _{\eta \in [- \pi , \pi ]^n } |\eta |^{ - \kappa } |\widehat{L _1 }(\eta + 2 \pi \ell ) | \leq C | \ell |^{ - N } .   
\end{equation}   
\end{lemma}   
   
\noindent {\it Proof.} The function $\sum _k \widehat{\varphi } (\eta + 2\pi k ) $ is periodic and, by the positivity and ellipticity of $\widehat{\varphi  } $ at 0, bounded from below by $c |\eta |^{- \kappa } $ for some $c > 0 . $ Hence     
$$   
|\widehat{L _1 } (\eta + 2 \pi \ell ) | = \frac{\widehat{|\varphi } (\eta + 2 \pi \ell ) | }{\sum _k \widehat{\varphi } (\eta + 2 \pi k ) } 
\leq C |\eta |^{\kappa } |\eta + 2 \pi \ell |^{- N } ,      
$$   
which implies (\ref{eq:SFC1}). \hfill $\Box $   
\medskip   
   
Another useful lemma clarifies the smoothness properties of $\widehat{L }_1 $:   
   
\begin{lemma} \label{lemma:der_L_hat} There exist constants $C_{\alpha } $ such that for each multi-index $\alpha $ with $|\alpha | \leq n + \lfloor \kappa \rfloor + 1 $ and all 
$k \in \mathbb{Z }^n $,   
\begin{equation} \label{eq:der_L_hat}   
\left | \partial ^{\alpha } _{\eta } \left( \widehat{L }_1 (\eta + 2 \pi k ) - \delta _{0k } \right) \right | \leq \frac{C_{\alpha } }{(1 + |k | )^N } |\eta |^{\kappa - |\alpha | } ,   
\end{equation}   
for $\eta \neq 0 $ in a neighborhood of 0.     
In particular, if $\kappa > 0 $ then $\widehat{L }_1 $ belongs to the H\"older space $C_b ^{\lceil \kappa \rceil - 1 , \lambda }  (\mathbb{R }^n ) $, with $\lambda = \kappa - (\lceil \kappa \rceil - 1 ) . $   
\end{lemma}   
   
\noindent Note that $\lceil \kappa \rceil - 1 = \lfloor \kappa \rfloor $ if $\kappa $ is non-integer, but that it is equal to $\kappa - 1 $ if $\kappa $ is a positive integer, so that $\lambda = 1 $ then.   
\medskip   
   
\noindent {\it Proof of lemma \ref{lemma:der_L_hat}}. This is elementary: if we let $\widehat{\varphi }_{\rm per } (\eta ) := \sum _k \widehat{\varphi } (\eta + 2 \pi k ) $, then applying Leibnitz's rule to the product $\widehat{L }_1 \widehat{\varphi }_{\rm per } = \widehat{\varphi } $ yields that   
$$   
(\partial ^{\alpha } _{\eta } \widehat{L }_1 ) \widehat{\varphi }_{\rm per } = \partial ^{\alpha } _{\eta } \widehat{\varphi } - \sum _{\beta < \alpha } \left( \begin{array}{cc} \alpha \\ \beta \end{array} \right) \partial ^{\beta } _{\eta } \widehat{L }_1 \, \partial ^{\alpha - \beta } _{\eta } \widehat{\varphi } _{\rm per } .   
$$   
The estimate (\ref{eq:der_L_hat}) for $k \neq 0 $ now follows by induction on $\alpha $, using that $\partial _{\eta } ^{\alpha } \widehat{\varphi }_{\rm per } (\eta + 2 \pi k ) = \partial _{\eta } ^{\alpha } \widehat{\varphi }_{\rm per } (\eta ) = O(|\eta |^{- \kappa - |\alpha | } ) $, together with  (\ref{eq:wB1}) of definition \ref{def:wBuhmann class} and lemma \ref{lemma: decay_L1_hat} (to start the induction). If $k = 0 $, we use the same argument, starting from   
$$   
\widehat{\varphi }_{\rm per } \left( \widehat{L }_1 - 1 \right) = \widehat{\varphi } - \widehat{\varphi }_{\rm per } ,   
$$   
on observing that the right hand side is $C^{\lfloor \kappa \rfloor + n + 1 } $ near 0, since equal to $\sum _{k \neq 0 } \widehat{\varphi }(\eta + 2 \pi k ) . $ Finally, the fact that $\widehat{L }_1 (\eta + 2 \pi k ) - \delta _{0k} $ is $O(|\eta |^{\kappa } $ implies that all derivatives of order up to $\lfloor \kappa \rfloor $, if $\kappa \notin \mathbb{N } $, or $\kappa - 1 $, if $\kappa \in \mathbb{N } \setminus 0 $, exist and are 0. Their continuity in 0 follows from (\ref{eq:der_L_hat}). \hfill $\Box $   
   
\begin{remark} \rm{We briefly pause to examine the differentiability of $\widehat{L }_1 $ if $\kappa \in \mathbb{N } . $ Letting $g(\eta ) := \sum _{k \neq 0 } \widehat{\varphi } (\eta + 2 \pi k ) $ and $\psi (\eta ) := |\eta |^{\kappa } \widehat{\varphi } (\eta ) $, we have that   
$$   
\widehat{L }_1 (\eta ) = \frac{\psi (\eta ) }{\psi (\eta ) + |\eta |^{\kappa } g(\eta ) } .   
$$   
This shows that $\widehat{L }_1 $ cannot be $C^{\kappa } $ in 0 if $\kappa \in \mathbb{N } $ is not even, even if $\psi $ is (note that then $\psi (0) \neq 0 $ given that $\varphi $ is Buhmann class). If $\kappa \in 2 \mathbb{N } $, then $\widehat{L }_1 $ will be as smooth as $\psi (\eta ) $ is in 0, and therefore as $\widehat{\varphi } $ is away from 0.   
}   
\end{remark}   
   
We finally note that to construct numerical PDE schemes using RBF interpolation one will obviously need sufficient differentiability of $L_1 . $ The proof of theorem
\ref{thm_L_1} given in appendix A also yields existence and decay of
derivatives of $L_1 $, provided $N $ is chosen sufficiently large:   
   
\begin{theorem} \label{lemma_L_1} Suppose that $k \in \mathbb{N } $ and let
$\varphi \in \mathfrak{B }_{\kappa , N } (\mathbb{R }^n ) $ with
$N > n + k . $ Then $L_1 \in C^k (\mathbb{R }^n ) . $ Moreover,
$|\partial _x ^{\alpha } L_1 (x) | = O(|x |^{- \kappa - n } ) $ as
$|x| \to \infty $, for all $|\alpha | \leq k . $
\end{theorem}   
\noindent See also theorem \ref{lemma:a(D)L_1} below.   
   
   

\section{\bf Convergence of RBF-interpolants}   
   
As a preparation for our analysis of numerical RBF-schemes for the heat equation and other evolution equations, we first revisit the convergene of stationary RBF interpolants, providing an alternative perspective on classical convergence theorems of Buhmann and Powell. As stated in the introduction, we will limit ourselves 
to stationary interpolation on regular grids $h \mathbb{Z }^n $, meaning that we let   
the basis function scale with the grid-size: $\varphi _h (x) := \varphi
(x/h ) . $ The associated Lagrange function scales similarly, and
the RBF interpolant $s_h [f] $ of a given function $f : \mathbb{R
}^n \to \mathbb{R } $ can be conveniently written as
\begin{equation} \label{eq:s_h[f]}   
s_h [f] (x) = \sum _j f(hj ) L_1 \left( \frac{x }{h } - j \right) .   
\end{equation}
where $L_1 $ is the Lagrange function of theorem \ref{thm_L_1}. Here, and below, sums over $j $, $k $, $\ell $, etc. are understood to
be over $\mathbb{Z }^n . $ Note that the use of the Lagrange   
function eliminates the need for inverting the coefficient matrix $(\varphi _h (hj - hk ) )_{j, k } = (\varphi (j - k ) )_{j, k }
$ in the standard formulation of RBF interpolation\footnote{In the present, idealized, set-up of interpolation on $h \mathbb{Z }^n $ that coefficient matrix is infinite; in practice, one would have to truncate the matrix: $|j |, |k | \leq N $ (where, $| j | = | j |_{\infty } = \max _{\nu } |j_{\nu } | $) with $N \sim h^{-1 } $, taking larger and larger sections of the matrix as $h \to 0 . $ One would also have to truncate the series for $L_1 $, leading to quasi-interpolation.}. The decay at infinity of $L_1 $ easily implies that the series (\ref{eq:s_h[f]}) converges absolutely if 
$f $ is of polynomial growth of order strictly less than $\kappa $, in the sense that $f \in L^{\infty } _{- p } (\mathbb{R }^n ) $ for some $p < \kappa . $     
   
Throughout this section, we fix a basis function $\varphi = \varphi _1 \in \mathfrak{B }_{\kappa , N } (\mathbb{R }^n ) $ with 
$N > n $ and $\kappa > 0 $: 
see subsection 3.2 below for the case of $\kappa = 0 $: although there can be no convergence then, this case is nevertheless of interest since we can have approximate convergence in the sense of Maz'ja and Smith: see secion 4 below. 
We will systematically work in Fourier-space, and examine convergence of $s_h [f ] $ to $f $ in Wiener norm, 
\begin{equation} \nonumber   
|| f ||_A 
:= || \widehat{f } ||_1 ,          
\end{equation}   
except for the end of this section where we will briefly examine convergence in weigthed sup-norms. 
Convergence in Wiener norm of course trivially implies convergence in Chebyshev or uniform norm, since  $ || f ||_{\infty } \leq || f ||_A . $   
   
\subsection{Convergence in Wiener norm} We begin by computing the Fourier transform of $s_h [f ] $ for Schwarz-class functions $f . $ For sufficiently rapidly decaying functions $g $, let us define the function $\Sigma _h (g) $   
\begin{equation} \label{eq:def_Sigma_h}   
\Sigma _h (g ) (\xi ) := \left( \sum _k g(\xi + 2 \pi h^{-1 } k ) \right) \widehat{L_1 } (h \xi ) .  
\end{equation}   
The map $\Sigma _h : g \to \Sigma _h (g) $ will play an important r\^ole in what follows. We note that $\Sigma _h $ is a contraction with respect to the $L^1 $-norm: indeed, by the positivity of $\widehat{L }_1 $ and monotone convergence,   
\begin{eqnarray*}   
|| \Sigma _h (g ) ||_1 &\leq & \sum _k \int _{\mathbb{R }^n } \, | g (\xi + 2 \pi h^{-1 } k ) | \, \widehat{L }_1 (h \xi ) \, d\xi \\   
&=& \int _{\mathbb{R }^n } |g (\xi ) | \left( \sum _k \widehat{L }_1 (h \xi + 2 \pi k ) \right) d\xi \\   
&=& || g ||_1 ,   
\end{eqnarray*}   
in view of (\ref{eq:sum_transl_L-hat}); 
$\Sigma _h $ therefore extends to a contraction on $L^1 (\mathbb{R } ) . $ 
We also note that if $g \in L^1 (\mathbb{R }^n ) $, then the defining series for $\Sigma _h (g) $ converges absolutely a.e., since   
$$   
\int _{ ]0 ,\pi ] ^n } \sum _k | g(\xi + 2 \pi k ) | d\xi = \int _{\mathbb{R }^n } | g(\xi ) | d\xi < \infty .   
$$

\begin{lemma} \label{lemma:FT_s_h} If $f \in \mathcal{S }(\mathbb{R }^n ) $ then $s_h [f ] \in L^1 (\mathbb{R }^n ) $ and 
then 
\begin{equation} \label{eq:FT_s_h}   
\widehat{s_h [f ] } 
= \Sigma _h (\widehat{f} ) .   
\end{equation}   
\end{lemma}   
   
\noindent {\it Proof}. Since $\kappa > 0 $, $L_1 $ is integrable by theorem \ref{thm_L_1} and therefore $|| s_h [f ] ||_1 \leq \left( h^n \sum _j |f(hj ) | \right) || L_1 ||_1 . $ Applying Fubini's theorem to the function $(j, x ) \to f(hj ) L_1 (h^{-1 } x - j ) e^{- i (x, \xi ) } $ on $\mathbb{Z }^n \times \mathbb{R }^n $ one finds
\begin{eqnarray} \nonumber   
\widehat{s_h [f] } (\xi ) &=& \left( \sum _j f(jh ) e^{ - i h (j ,
\xi ) } \right) h^{n } \widehat{L }_1 (h \xi ) \nonumber  \\
&=& \left( \sum _k \widehat{f }(\xi + 2 \pi h^{-1 } k ) \right)
\widehat{L_1 } (h \xi ) , \label{eq:s_h_hat} \nonumber \\   
&=& \Sigma _h (\widehat{f } )(\xi ) ,      
\end{eqnarray}
where for the second line we used the Poisson summation formula: $\sum _j g(j)  = \sum _k \widehat{g } (2 \pi k ) $, with $g(x) := f(hx ) e^{- i h (x, \xi ) } . $ \hfill $\Box $   
\medskip   
      
We can then already state a first convergence theorem: 
      
\begin{theorem} \label{thm:convergence_RBF} Let $\kappa > 0 . $ Then there exists a constant $C = C_{\varphi } > 0 $ such that for all tempered functions $f $ for which $\widehat{f } \in L^1 _{\kappa } (\mathbb{R }^n ) $ and for all positive $h \leq 1 $,
\begin{equation} \label{eq:conv_estimate_1}      
|| \, f - s_h [f] \, ||_A 
\leq C h^{\kappa } || f ||^{\circ } _{\kappa } = C h^{\kappa } \int _{\mathbb{R }^n } |\widehat{f }(\xi ) | \, |\xi |^{\kappa } \,  d\xi .   
\end{equation}   
\end{theorem}   
   
\noindent The condition that $\widehat{f } \in L^1 _{\kappa } (\mathbb{R }^n ) $ implies a certain smoothness: $f $ must have continuous and bounded derivatives of order $\lfloor \kappa \rfloor . $   
\medskip

\noindent {\it Proof}. We first note that if $\widehat{f } \in L^1 (\mathbb{R }^n ) $, then $\widehat{s_h [f ] } = \Sigma _h (\widehat{f } ) $: the hypothesis on $\widehat{f } $ implies that $f $ is a bounded continuous function. 
It follows that $s_h [f ] $ is well-defined, by (\ref{bound_L_1}), and that there exists a constant $C > 0 $ such that for all $h \leq 1 $,     
\begin{equation}   
|| \, s_h [f ] \, ||_{\infty } \leq C ||f ||_{\infty } .    
\end{equation}   
Indeed, $ | s_h [f ](x) | \leq || f ||_{\infty } \sum _j | L_1 (h^{-1 } x - j ) | $; the right had side is $h $-periodic, and its sup on $\{ |x | \leq h / 2 \} $ can be estimated by a constant times $\sum _j | L_1 (j) | $, which converges since $\kappa > 0 . $      
   
The Fourier transform of $s_h[f ] $ therefore exists as a temperered distribution. We show using a density argument that $\widehat{s_h [f ] } = \Sigma _h (\widehat{f } ) $:  since $\widehat{f } \in L^1 (\mathbb{R }^n ) $, then there exists a sequence $f_{\nu } \in \mathcal{S } (\mathbb{R }^n ) $ such that $|| \widehat{f }_{\nu } - \widehat{f } ||_1 \to 0 . $ Consequently $\Sigma _h (\widehat{f }_{\nu } ) \to \Sigma _h (f) $ in $L^1 $ and therefore also as tempered distributions. On the other hand, $|| s_h [f_{\nu } ] - s_h [f ] ||_{\infty } \leq C || f - f_{\nu } ||_{\infty } \leq C || \widehat{f } - \widehat{f }_{\nu } ||_1 \to 0 $, so $s_h [f_{\nu } ] \to s_h [f ] $ as tempered distributions also.    Hence $\Sigma _h (\widehat{f }_{\nu } ) = \widehat{s _h [f_{\nu } ] } \to \widehat{s_h [f ] } $, and consequently $\widehat{s_h [f ] } = \Sigma _h (\widehat{f } ) . $ As a consequence, $s_h [f ] $ has finite Wiener norm if $f $ has. 
   
We now observe that since $0 \leq \widehat{L }_1 \leq 1 $, and using the montone convergence,   
\begin{eqnarray}   
|| f - \widehat{s }_h [f ] ||_1 &\leq & \int _{\mathbb{R }^ n } \, |\widehat{f }(\xi ) | (1 - \widehat{L }_1 (h\xi ) ) d \xi + \sum _{k \neq 0 } \int _{\mathbb{R }^ n } |\widehat{f }(\xi + 2 \pi h^{-1 } k ) | \, \widehat{L }_1 (h \xi ) \, d\xi \nonumber \\   
&=& \int _{\mathbb{R }^ n } \, |\widehat{f }(\xi ) | \left( (1 - \widehat{L }_1 (h\xi ) ) + \sum _{k \neq 0 } \widehat{L }_1 (h \xi + 2 \pi k ) \right) \, d\xi \nonumber \\   
&=& 2 \int _{\mathbb{R }^n } (1 - \widehat{L }_1 (h \xi ) ) \, |\widehat{f }(\xi ) | \, d\xi , \label{eq:proof_convergence_RBF}
\end{eqnarray}   
where we used (\ref{eq:sum_transl_L-hat}) once more.      
The Fix-Strang condition (\ref{eq:SFC}) in 0 
then implies (\ref{eq:conv_estimate_1}) with $C = 2 \sup _{\eta \neq 0 } (1 - \widehat{L }_1 (\eta ) ) / |\eta |^{\kappa } $ (a number which, in principle at least, is explicitly computabe for  a given $\varphi $).   
   
\hfill $\Box $   
\medskip   
   
\begin{remark} \rm{\noindent The theorem generalizes to the case when $\widehat{f } = \nu $ is a finite Borel measure such that  $|\xi |^{\kappa } \in L^1 (\mathbb{R }^n , d |\nu | ) $: in that case,     
\begin{equation} \label{eq:conv_estimate_1a}      
|| \, f - s_h [f] \, ||_{\infty } 
\leq C h^{\kappa } \int _{\mathbb{R }^n } \, |\xi |^{\kappa } \, d | \nu | (\xi ) .   
\end{equation}   
To show this, one first defines $\Sigma _h (\nu ) $ by duality: if $\psi \in \mathcal{S } (\mathbb{R }^n ) $, then   
$$   
\langle \Sigma _h (\nu ) , \psi \rangle := \langle \nu , \Sigma ' (\psi ) \rangle ,   
$$  
where $\Sigma _h ' (\psi ) := \sum _k \psi (\xi + 2 \pi h^{-1 } k ) \, \widehat{L }_1 (h \xi + 2 \pi k ) \in C_b (\mathbb{R }^n ) $, and one checks that $\widehat{s_h [f ] } = \Sigma _h (\nu ) $ as tempered distributions. Since $|| \Sigma _h ' (\psi ) ||_{\infty } \leq C ||\psi ||_{\infty } $, on account of the decay of $\widehat{L }_1 $, $\Sigma _h (\nu ) $ is a finite Borel measue. Using (\ref{eq:sum_transl_L-hat}) again, 
one estimates   
$$   
\left | \langle \Sigma _h (\nu ) - \nu , \psi \rangle \right | = \left | \langle \nu , \Sigma _h ' (\psi ) - \psi \rangle \right | \leq 2 || \psi ||_{\infty } \int _{\mathbb{R }^n } (1 - \widehat{L }_1 (h\xi ) ) d | \nu \ (\xi ) ,   
$$   
where we can take $\psi \in C_b (\mathbb{R }^n ) . $ It follows that the variation norm of $\Sigma _h (\nu ) - \nu $ is bounded by $C h^{\kappa } $, which implies (\ref{eq:conv_estimate_1a}).   
}   
\end{remark}      
   
We next observe that the right hand side of (\ref{eq:conv_estimate_1}) still makes sense for certain $\widehat{f } $ having a non-integrable singularity at 0. Allowing such singularities means allowing $f $'s which grow at a certain polynomial rate, and we can 
prove the following extension of theorem \ref{thm:convergence_RBF}:   
   
   
\begin{theorem} \label{thm:convergence_RBF_bis} Let $f $ be a tempered function on $\mathbb{R }^n $ such that $|f(x) | \leq C (1 + |x | )^p $ for some $p < \kappa $, and such that         
\begin{equation} \label{eq:convergence_RBF_bis}   
\widehat{f } |_{\mathbb{R } ^n \setminus 0 } \in 
\mathring{L }^1 _{\kappa } (\mathbb{R }^n ) .   
\end{equation}   
Then $\widehat{s _h [f ] } - \widehat{f } $ is in $L^1 (\mathbb{R }^n ) $, and   
\begin{equation} \label{eq:conv_RBF_bis}   
|| s_h [f ] - f ||_A \leq C h^{\kappa } || \widehat{f } ||^{\circ } _{1, \kappa } . 
\end{equation}   
\end{theorem}   
   
\noindent Equation (\ref{eq:convergence_RBF_bis}) means that the restriction to $\mathbb{R }^n \setminus 0 $ of the tempered distribution $\widehat{f } $  can be identified with a locally integrable function which is integrable with respect to the weight $|\xi |^{\kappa } $ and therefore is in $\mathring{L }^1 _{\kappa } (\mathbb{R }^n ) $, if we interpret it as an almost everywhere defined function on $\mathbb{R }^n $, whose norm we then simply denote by $|| \widehat{f } ||^{\circ } _{1 , \kappa } $ instead of the more correct $|| \widehat{f } |_{\mathbb{R }^n \setminus 0 } ||^{\circ } _{1, \kappa } . $   
\medskip   
   
\noindent {\it Proof}. We first check that $s_h [f ] $ is a tempered distribution: this is a consequence of the estimate   
\begin{equation} \label{eq:s_h_assertion_1}   
|| \, s_h [f ] \, ||_{\infty , -p } \leq C || \, f \, ||_{\infty , -p } , \ \ p \geq 0 ,   
\end{equation}   
which can be shown as follows: first note that $f \to s_h [f ] $ commutes with translations by elements of $h \mathbb{Z }^n $: if $k \in \mathbb{Z }^n $, then   
$$   
s_h [f ] (x - kh ) = s_h [f (\cdot - hk ) ] (x) .   
$$   
Let $| \cdot  | = | \cdot |_{\infty } $ be the $\ell ^{\infty } $-norm on $\mathbb{R }^n . $ If $|x | \leq h/2 $ and $f \in L^{\infty } _{-p } (\mathbb{R }^n ) $ with $p < \kappa $, then      
$$     
| s_h [f ](x) | \leq || f ||_{\infty , -p } \left( 1 + \sum _{| j | \geq 1 } \frac{ (1 + h |j | )^p }{(1 + | h^{-1 } x - j | )^{\kappa + n } } \right) \leq C || f ||_{\infty , -p } ,     
$$   
since $| h^{-1 } x - j | \geq |j |/2 $ if $|j | \geq 1 . $ Next, if $|x - hk | \leq h/2 $ with $k \in \mathbb{Z }^n $, then by translation invariance,   
$$   
| s_h [f ] (x) | \leq C || f (\cdot + hk ) ||_{\infty , -p } \leq C (1 + h |k | )^p || f ||_{\infty , -p } ,   
$$
which implies (\ref{eq:s_h_assertion_1}). \noir The next lemma identifies the Fourier transform if $s_h [f ] . $   
   
\begin{lemma} \label{lemma:convergence_RBF_bis} Suppose that $|f(x) | \leq C (1 + |x | )^p $ for some $p < \kappa $ and that $\widehat{f } |_{\mathbb{R } \setminus 0 } \in L^1 (\mathbb{R }^n \setminus 0 , \min ( |\xi |^{\kappa } , 1 ) d\xi ) . $ Then  
the tempered distribution $\widehat{s_h [f ] } - \widehat{f } $ can be identified with the function 
\begin{equation} \label{eq:lemma_convergence_RBF_bis}   
\left( \widehat{L }_1 (h \xi ) - 1 \right) \widehat{f }(\xi ) + \sum _{k \neq 0 } \widehat{f } (\xi + 2 \pi h^{-1 } k ) \widehat{L }_1 (h \xi ) , \ \ \xi \neq 0 ,   
\end{equation}   
which is in $L^1 (\mathbb{R }^n ) . $   
\end{lemma}   

The proof of the lemma involves extending $\widehat{f } $ to a continuous linear functional on the H\"older spaces $C_b ^{\lceil \kappa \rceil - 1 , \lambda }  (\mathbb{R }^n ) $ with $\lambda = \kappa - (\lceil \kappa \rceil - 1 ) $ (so that $\lambda = \kappa - \lfloor \kappa \rfloor $ if $\kappa \notin \mathbb{N } $, and $\lambda = 1 $ otherwise), and using this to define $\Sigma _h (\widehat{f } ) $ as a tempered distribution.  In order not to interrupt the flow of the argument with distribution-theoretical technicalities, we postpone the proof     
to Appendix \ref{Appendix_proof_thm_conv_RBF_bis}. Note that the individual terms of (\ref{eq:lemma_convergence_RBF_bis}) are integrable on account of the Fix-Strang conditions satisfied by $\widehat{L }_1 $, and that the $L^1 $-norm of (\ref{eq:lemma_convergence_RBF_bis}) can be bounded by the $L^1 $-norm of $2 ( \widehat{L }_1 (h \xi ) - 1 ) |\widehat{f } (\xi ) | $, using once more that the sum of translates of $\widehat{L }_1 $ by elements of $(2 \pi ) \mathbb{Z }^n $ is identically equal to one.   
   
The lemma implies the estimate (\ref{eq:proof_convergence_RBF}), and the theorem follows as before. 
   
\hfill $\Box $   
   
\begin{example} \rm{If $p \geq 0 $ and if $f \in C^{\lfloor p \rfloor + n + 1 } (\mathbb{R }^n ) $ satisfies   
\begin{equation} \label{eq:S^p}   
| \partial _x ^{\alpha } f (x) | \leq C_{\alpha } (1 + |x| )^{p - |\alpha | } , \ \ |\alpha | \leq \lfloor p \rfloor + n + 1 ,   
\end{equation}   
with $p \geq 0 $ then one can show that $\widehat{f } |_{\mathbb{R }^n \setminus 0 } \in C(\mathbb{R }^n \setminus 0 ) $, and that $|\widehat{f } (\xi ) | \leq C |\xi | ^{- p - n } $ near 0 while $\widehat{f } (\xi ) = O(|\xi |^{ - \lfloor p \rfloor - n - 1 } ) $ at infinity: if (\ref{eq:S^p}) holds for all $\alpha $, this follows for example from Stein \cite{Stein},  proposition 1 of Chapter VI. Examination of the proof shows that we only need the number of derivatives indicated. It follows that $|\xi |^{\kappa } \widehat{f } (\xi ) $ is integrable if $p < \kappa $ and 
theorem \ref{thm:convergence_RBF_bis} therefore applies to such functions.   
\medskip   
   
We briefly compare our theorem \ref{thm:convergence_RBF_bis} with convergence results of Buhmann and Powell, cf. \cite{Bu_book} and its references. Theorem 4.6 of \cite{Bu_book} states that if $\kappa \notin \mathbb{N } $ and if $f \in C^{\lceil \kappa \rceil } (\mathbb{R }^n ) $ such that   
\begin{equation} \label{eq:cond_Buhmann}   
\max _{|\alpha | = \lceil \kappa \rceil, \lceil \kappa \rceil - 1 } || \partial _x ^{\alpha } f ||_{\infty } < \infty ,   
\end{equation}   
then $|| s_h [f] - f ||_{\infty } \leq C h^{\kappa } . $ If $\kappa $ is an odd integer, Buhmann {\it loc cit.}, theorem 4.7, needs an additional degree of differentiability, with the derivatives of order $\kappa + 1 $ again bounded. 
Note that these condititons imply that $f $ is of polynomial growth of order at most $\lceil \kappa \rceil - 1 . $ On the other hand, theorem \ref{thm:convergence_RBF_bis} covers cases when $f $ can have stronger growth, and derivatives of order $\lceil \kappa \rceil - 1 $ do not not need to be bounded, for example  $f(x) = (1 + |x | )^{p/2 } $ with $\lfloor \kappa \rfloor < p < \lceil \kappa \rceil $ if $\kappa \notin \mathbb{N } $ and $\kappa - 1 < p < \kappa $ if $\kappa \in \mathbb{N }^* $. We also note that condition (\ref{eq:convergence_RBF_bis}) implies that $f $ is $C^{\lfloor \kappa \rfloor } $ (a consequence of the integrability of $\widehat{f } (\xi ) |\xi |^{\kappa } $ on $|\xi | \geq 1 $ and the smoothness of the inverse Fourier transform of any compactly supported distribution) without the derivatives of order $\lfloor \kappa \rfloor $ necessarily being bounded (because of the singularity at 0) and we believe 
this order of differentiability to be the maximal one necessary.   
   
}   
\end{example}   
   
\begin{remark} \rm{If $\kappa \notin \mathbb{N } $ then theorem \ref{thm:convergence_RBF_bis} remains true if $\widehat{f } | _{\mathbb{R }^n \setminus 0 } $ can be identified with a Borel measure $\nu $ on $\mathbb{R }^n \setminus 0 $ for which $|\xi |^{\kappa } \in L^1 (\mathbb{R }^n , d|\nu | ) . $ The estimate (\ref{eq:conv_RBF_bis}) then generalises to an estimate for the variation norm of $\Sigma (\widehat{f } ) - \widehat{f } $ (as measure on $\mathbb{R }^n $) which then implies a uniform estimate   
\begin{equation}   
|| s_h [f ] - f ||_{\infty } \leq C h^{\kappa } \int _{\mathbb{R }^n } |\xi |^{\kappa } \, d |\nu |(\xi ) .   
\end{equation}   
}   
\end{remark}   
\medskip   
      
We next observe that if one is satisfied with a slower rate of convergence, the growth condition on $\widehat{f } $ at infinity can be weakened. Recall the spaces $\mathring{L }^1 _{r, s } (\mathbb{R }^n ) $ defined in (\ref{eq:L^1_rs}). 
         
\begin{corollary} \label{cor:convergence_RBF} (of the proof). Suppose that $0 \leq r \leq \kappa $ and that   
$f \in L^{\infty } _{ - \kappa + \varepsilon } (\mathbb{R }^n ) $for some $\varepsilon > 0 $ 
such that   
$$   
\widehat{f } |_{\mathbb{R }^n \setminus 0 }  \in \mathring{L }^1 _{r, \kappa } (\mathbb{R }^n ) . 
$$   
Then  for all $h \leq 1 $,   
$$   
|| \widehat{s_h [f] } - \widehat{f } ||_1 \leq C \,  h^r \, || \widehat{f } ||^{\circ } _{r, \kappa } .   
$$   
\end{corollary}   
   
\noindent Here, the constant $C $ can be taken the same as in theorems \ref{thm:convergence_RBF_bis} and \ref{thm:convergence_RBF} (cf. the proof of the latter).   
\medskip   
   
\noindent {\it Proof}. 
The hypotheses on $\widehat{f } $ certainly imply that $\widehat{f } |_{\mathbb{R }^n } \in L^1 \left( \mathbb{R }^n , \min (|\xi |^{\kappa } , 1 ) d\xi \right) $, so we can apply lemma \ref{lemma:convergence_RBF_bis}. In particular, the estimate (\ref{eq:proof_convergence_RBF}) still holds. Now split this integral into three parts 
over the ranges $h|\xi | \leq h $, $h \leq h |\xi | \leq 1 $ qnd $h |\xi | \geq 1 $, and use the Fix - Strang condition in 0, together with the trivial bound $|\eta |^{\kappa } \leq |\eta |^r $ if $|\eta | \leq 1 $ to estimate:   
\begin{eqnarray*}   
&&\int _{\mathbb{R }^n } \left( 1 - \widehat{L }_1 (h \xi ) \right) | \widehat {f }(\xi ) | \, d\xi \\   
&\leq & C \left( \int _{|h \xi | \leq h } \, (h |\xi | )^{\kappa } |\widehat{f } (\xi ) | \, d\xi  + \int _{h < |h \xi | \leq 1 } \, (h |\xi | )^{\kappa } |\widehat{f } (\xi ) | \, d\xi + \int _{|h \xi | > 1 } \, | \widehat{f } (\xi ) | \, d\xi \right) \\   
&\leq & C \left( h^{\kappa } \int _{|\xi | \leq 1 } \, |\xi |^{\kappa } |\widehat{f } (\xi ) | \, d\xi + h^r \int _{1 \leq |\xi | \leq h^{-1 } } \, |\xi |^r |\widehat{f } (\xi ) | \, d\xi + h^r \int _{|\xi | \geq h^{-1 } } \, | \widehat{f } (\xi ) | \, |\xi |^r \, d\xi \right) \\   
&\leq & C h^r \int _{\mathbb{R }^n } \, | \widehat{f }(\xi ) | \, |\xi |^r \wedge |\xi |^{\kappa } \, d\xi .   
\end{eqnarray*}   
   
\hfill $\Box $   
   
\medskip   
   
The corollary shows that there is an interplay between the order of convergence of the RBF interpolants and the smoothness of the function which is interpolated, as quantified by the decay of $\widehat{f }(\xi ) $ at infinity. The singularity of $\widehat{f }(\xi ) $ at 0 can be of order $|\xi |^{- \kappa - n + \varepsilon } $, as before, which is compatible with $f(x) $ having polynomial growth  of order less than $\kappa . $   
\medskip       
   
\subsection{The case of $\kappa = 0 $} 
The results of this section, and in particular the basic estimate  (\ref{eq:proof_convergence_RBF}), remain true if $\kappa = 0 $, which includes for example Gaussian basis functions.     The reader may of course wonder why one would want to consider this case, since this estimate then does not imply that $s_h [f ] $ will converge to $f $ and, as we will see in section 4, it won't. The reason is that we can still have approximate approximation, in the sense that $\limsup _{h \to 0 } || s_h [f ] - f ||_A $ can be made  arbitrarily small for functions $f $ whose Fourier transform decay sufficiently rapidly, by an appropriate choice of basis function: see theorem \ref{thm:approx_approx_RBF} and the discussion following it. There is a minor technical problem, in that theorem \ref{thm_L_1} no longer guarantees that $L_1 (x) $ decays sufficiently rapidly to be integrable, though it may for particular examples such as the Gaussian, or under stronger conditions on $\varphi $, as per Buhmann's result for integer pair $\kappa $mentioned before. The interpolating function $s_h [f ] $ is therefore no longer be guaranteed to be in $L^1 (\mathbb{R }^n ) $, even if $f $ is rapidly decreasing. Its Fourier transform will still exist as a tempered distribution, and lemma \ref{lemma:FT_s_h} will still be true, 
as an easy approximation argument shows:   
$\sum _{|j | \leq N } f(hj ) L_1 (h^{-1 } x - j ) \to s_h [f ] $   
in $\mathcal{S } ' (\mathbb{R }^n ) $ so $\left( \sum _{|j | \leq N } f(hj ) e^{ - i h (j ,\xi ) } \right) h^{n } \widehat{L }_1 (h \xi ) \to \widehat{s_h [f ] } $   
as tempered distributions. But $\left( \sum _j f(jh ) e^{ - i h (j , \xi ) } \right) h^n \widehat{L }_1 (h \xi ) $ converges to $\Sigma _h (\widehat{f } ) . $   
   
Next, the series for $s_h [f ] $ will not necessarily converge if $\widehat{f } \in L^1 (\mathbb{R }^n ) $ and we need to impose an additional condition that $f \in L^{\infty } _{\varepsilon } (\mathbb{R }^n ) $ for some $\varepsilon > 0 $, as in theorem \ref{thm:convergence_RBF_bis}. Alternatively, we can observe that if $\widehat{f } $ is integrable, then the defining series for $s_h [f ] $ is summable 
in the sense that  if $\chi \in \mathcal{S } (\mathbb{R }^n ) $ with $\chi (0) = 1 $ and if we let     
\begin{equation} \nonumber \label{eq:s_h_kappa = 0}   
s_h ^{\varepsilon } [f ] (x) := \sum _{j } \chi (\varepsilon jh ) f(hj ) L_1 (h^{-1 } x - j ), 
\end{equation}   
then, as $\varepsilon \to 0 $, $s_h ^{\varepsilon } [f ] $ converges uniformly on $\mathbb{R }^n $ to a continuous function whose Fourier transform is $\Sigma _h (\widehat{f } ) $, independently of the choice of $\chi . $ 
To show this, observe that the Fourier transform of the left hand side is equal to $(2 \pi )^{-n } \Sigma _h (\widetilde{\chi } _{\varepsilon } * \widehat{f } ) $, where $\widetilde{\chi } _{\varepsilon } (\xi ) = \varepsilon ^{-n } \widehat{\chi } (- \varepsilon \xi ) . $ Since (by a classical result on convolution with approximate identities),   
$$   
(2 \pi )^{-n } \widetilde{\chi }_{\varepsilon } * \widehat{f } \to \widehat{f }   
$$   
in $L^1 $ (observing that $(2 \pi )^{- n } \int _{\mathbb{R }^{-n } } \widehat{\chi } (- \xi ) d\xi = \chi (0) = 1 $), and since $\Sigma _h $ is a contraction, it follows that 
$s_h ^{\varepsilon } [f ] $ 
converges in Wiener norm, and therefore in sup-norm, to the inverse Fourier transform of the, integrable, function $\Sigma _h (\widehat{f } ) . $ If we now define $s_h[f ] $ as the limit of the $s_h ^{\varepsilon } [f ] $, the estimate (\ref{eq:proof_convergence_RBF}) follows as before. We in fact only need the Fourier transform of $\chi $ to be integrable, but this excludes taking for $\chi $ the characteristic function of a cube centered at 0, which would entail ordinary convergence of the series for $s_h [f ] (x) . $   
\medskip

\subsection{Convergence in weighted sup-norms} Although our focus in this paper is on convergence in the Wiener norm, one can allow $f $'s for which 
$\widehat{f } $ does not necessarily coincide with a function or a measure on $\mathbb{R } \setminus 0 $, but is a more general distribution, if we replace the Wiener norm by one of the weighted sup-norms $|| \cdot  ||_{\infty , -p } . $ We give an example which can be deduced from theorem \ref{thm:convergence_RBF} by an approximation argument.
   
\begin{theorem} \label{thm:conv_RBF_3} Let $p \in \mathbb{N } $, $p< \kappa $ and let $f $ be a tempered function on $\mathbb{R }^n $ whose Fourier transform can be written as   
\begin{equation} \label{eq:conv_RBF_3_i}   
\widehat{f } = \sum _{|\alpha | \leq p } \partial _{\xi } ^{\alpha } \nu_{\alpha } ,   
\end{equation}   
with $\nu _{\alpha } $ complex Borel measures on $\mathbb{R }^n $ satisfying   
\begin{equation} \label{eq:conv_RBF_3_ii}   
\int _{\mathbb{R }^n } \, (1 + |\eta | )^{\kappa } d|\nu _{\alpha } | (\eta ) < \infty .   
\end{equation}   
Then   
\begin{equation} \label{eq:conv_RBF_3_iii}   
|| s_h [f ] - f ||_{\infty , -p } 
\leq C_f h^{\kappa } ,   
\end{equation}   
where we can take      
\begin{equation} \label{eq:conv_RBF_3_iv}   
C_f = C \cdot \sum _{|\alpha | \leq p } || (1 + |\eta | )^{\kappa } ||_{L^1 (|\nu _{\alpha | } ) } ,   
\end{equation}   
for some positive constant $C $ independent of $f . $   
\end{theorem}   
   
\noindent {\it Proof}. The hypothesis on $\widehat{f } $ imply that $f $ is continuous and of polynomial growth of order at most $p $: $|| f ||_{\infty , -p } < \infty . $ Let $\chi _R (x) := \chi (x/R ) $, where $\chi = \chi _1 \in C^{\infty } _c (\mathbb{R }^n ) $, $0 \leq \chi \leq 1 $, $\chi (x) = 1 $ on $B(0, 1 ) . $ We will apply theorem \ref{thm:convergence_RBF} to $\chi _R f $ and for that purpose first bound $|| \widehat{\chi _R f } ||_{1, \kappa } . $   
       
\begin{lemma} For $f $ as in theorem \ref{thm:conv_RBF_3} and $R \geq 1 $,   
\begin{equation}   
|| \widehat{\chi _R f } ||_{1, \kappa } \leq C_f R^p ,   
\end{equation}   
with $C_f $ as in (\ref{eq:conv_RBF_3_iv}).   
\end{lemma}   
   
\noindent {\it Proof}. Since $\widehat{f \chi _R } = (2 \pi )^{-n } \widehat{f } * \widehat{\chi _R } = (2 \pi )^{-n } \sum _{\alpha } (-1 )^{|\alpha | } \nu _{\alpha } * \partial ^{\alpha } _{\xi } \widehat{\chi _R }  $, we find that (writing $\widehat{\chi }^{(\alpha ) } $ for $\partial _{\xi } ^{\alpha } \widehat{\chi } $)   
\begin{eqnarray*}   
(2 \pi )^n || \widehat{\chi _R f } ||_{1, \kappa } &\leq & \sum _{|\alpha | \leq p } R^{n + |\alpha | } \int _{\mathbb{R }^n } \int _{\mathbb{R }^n }  |\widehat{\chi }^{(\alpha ) } (R(\xi - \eta ) ) | (1 + |\xi | )^{\kappa } d |\nu _{\alpha } | (\eta ) d\xi \\   
&=& \sum _{|\alpha | \leq p } R^{|\alpha | } \int _{\mathbb{R }^n } \int _{\mathbb{R }^n }  |\widehat{\chi }^{(\alpha ) } (\zeta ) | (1 + | \eta  + R^{-1 } \zeta | )^{\kappa } d |\nu _{\alpha } | (\eta ) d\zeta .   
\end{eqnarray*}   
The lemma follows by observing that  $(1 + | \eta + R^{-1 } \zeta  | ) 
\leq (1 + |\eta | ) (1 + R^{-1 } |\zeta | ) \leq (1 + |\eta | ) (1 + |\zeta | ) $ and using the rapid decay of $\widehat{\chi } . $   
   
\hfill $\Box $   
\medskip   
      
\noindent {\it Proof of theorem \ref{thm:conv_RBF_3} (continued)}. By theorem \ref{thm:convergence_RBF},   
$$   
|| s_h [ f \chi _R ] - f \chi _R ||_{\infty } \leq C_f R^p h^{\kappa } .   
$$   
We next compare $s_h [f ] $ with $s_h [\chi _R f ] $: since $\chi _R (x) = 1 $ for $|x | \leq R $,        
\begin{eqnarray*}   
|s_h [f] (x) - s_h (\chi _R f ) (x) | &=& \sum _{h |j | \geq R } | (f(hj ) - \chi _R (hj ) f(hj ) ) L_1 (h^{-1 } x - j ) | \\   
&\leq & 2 \sum _{h |j | \geq R } |f(hj ) | \, |L_1 (h^{-1 } x - j ) | .     
\end{eqnarray*}     
Now if $|x | \leq R/2 $, then $|hj | \geq R $ implies that $|x - hj | \geq |hj | / 2 $ so that $|h^{-1 } x - j | \geq |j | / 2 . $ Hence, by the decay at infinity of $L_1 $,   
\begin{eqnarray*}   
\sup _{|x | \leq R/2 } \sum _{h |j | \geq R } \, |f(hj ) | \, L_1 (h^{-1 } x - j ) &\leq & || f ||_{\infty , -p }  \, h^p \sum _{h |j | \geq R } |j |^{p - \kappa - n } \\   
&\leq & C || f ||_{\infty , -p } \, h^{\kappa } R^{p - \kappa } ,      
\end{eqnarray*}   
since we can for example bound the sum by a constant times $\int _{|y | \geq R / h } \, |y| ^{p - \kappa - n } dy $ (recall that $p < \kappa $).   
   
Writing $s_h [f] - f = s_h [f]  - s_h [\chi _R f ] + s_h [\chi _R f ] - \chi _R f + \chi _R f - f $, these estimates imply that 
$$   
R^{-p } \sup _{|x | \leq R/2 } |s_h [f ] (x) - f(x) | \leq C_f h^{\kappa } ,   
$$   
for $R \geq 1 $, which implies the theorem. \hfill $\Box . $   
\medskip   
   
Examples of functions $f $ which satisify the hypothesis of theorem \ref{thm:conv_RBF_3} are the inverse Fourier transforms of compactly supported distributions of order $p < \kappa $ since, by a structure theorem going back to Laurent Schwartz, such a compactly supported distribution can be written in the form (\ref{eq:conv_RBF_3_i})). 
   
\section{\bf Approximate approximation}   
   
   
It is easy to show that the approximation error in the Wiener norm cannot in general go to 0 faster than $h^{\kappa } $: if $\widehat{f } \in L^1 (\mathbb{R }^n ) $ has compact support, then the supports of $\widehat{f } (\cdot + 2 \pi k / h ) $ will be disjoint if $h $ is sufficiently small. It follows that   
\begin{eqnarray*}   
|| \widehat{s_h [f] } - \widehat{f } ||_1 &=& \int _{\mathbb{R }^n } |\widehat{f } (\xi ) (\widehat{L }_1 (h \xi ) - 1 ) | d\xi + \sum _{k \neq 0 } \int _{\mathbb{R }^n } |\widehat{f }(\xi + 2 \pi k / h ) | \, \widehat{L }_1 (h \xi ) d\xi \\   
&=& \int _{\mathbb{R }^n } |\widehat{f } (\xi ) | (1 - \widehat{L }_1 (h \xi ) ) d\xi + \sum _{k \neq 0 } \int _{\mathbb{R }^n } |\widehat{f }(\xi ) | \, \widehat{L }_1 (h \xi + 2 \pi k ) d\xi \\   
&=& 2 \int _{\mathbb{R }^n } |\widehat{f }(\xi ) | \left( 1 - \widehat{L }_1 (h \xi ) \right) d\xi ,   
\end{eqnarray*}   
since $\sum _k \widehat{L }_1 (\eta + 2 \pi k ) = 1 . $ If we define      
\begin{equation} \label{eq:def_underline_l}   
\underline{l }_{\kappa } := \underline{l }_{\kappa } (\varphi ) := 2 \liminf _{\eta \to 0 } \frac{1 - \widehat{L }_1 (\eta ) }{|\eta |^{\kappa } } .   
\end{equation}   
then Fatou's lemma implies that   
\begin{equation} \label{eq:lb_approx_error}   
\liminf_{h \to 0 } h^{- \kappa } || s_h [f] - f ||_A \geq \underline{l }_{\kappa } \int _{\mathbb{R }^n } \, |\xi |^{\kappa } |\widehat{f }(\xi ) | \, d\xi .   
\end{equation}   
We will see below that $\underline{l }_{\kappa } > 0 . $ The inequality (\ref{eq:lb_approx_error}) remains valid if $\widehat{f } $ is not compactly supported but decays sufficiently fast at infinity:    
see theorem \ref{thm:lower-bound_approx_error} below. Here we first examine the corresponding upper bound. \medskip   
   
As we just noted, one cannot in general do better that $O(h^{\kappa }) $ for the approximation error. However, for suitable basis functions $\varphi $ and for $\widehat{f }(\xi ) $ which decay sufficiently fast at infinity we may observe a higher {\it apparent} rate of convergence for $h $'s which are small but not too small. If $\varphi \in \mathfrak{B }_{\kappa , N } (\mathbb{R }^n ) $, we let\footnote{The index $\kappa $ is a reminder of the degree of the singularity of $\widehat{\varphi } $ at 0, and therefore of the natural convergence rate of the RBF interpolants.}      
\begin{equation}   
\overline{l }_{\kappa } := \overline{l }_{\kappa } (\varphi )) := 2 \limsup _{|\eta | \to 0 } \frac{1 - \widehat{L }_1 (\eta ) }{|\eta |^{\kappa } } .   
\end{equation}   
A slight modification of the proof of theorems \ref{thm:convergence_RBF} and \ref{thm:convergence_RBF_bis} then gives the following more precise estimate for the approximation error.   
   
   
\begin{theorem} \label{thm:approx_approx_RBF} Let $\varphi \in \mathfrak{B }_{\kappa , N } (\mathbb{R }^n ) $ with $\kappa \geq 0 $ and 
suppose that $f \in L^{\infty } _{- p } (\mathbb{R }^n ) $ for some $p < \kappa $ such that for some $s > \kappa $,   
$$   
\widehat{f } |_{\mathbb{R }^n \setminus 0 } \in \mathring{L }^1  _{\kappa } (\mathbb{R }^n ) \cap \mathring{L }^1 _s (\mathbb{R }^n ) =  L^1 (\mathbb{R }^n , \max \left( |\xi |^{\kappa } , |\xi |^s ) \right) .   
$$   
Then there exists for each $\varepsilon > 0 $, a constant $C_{\varepsilon } = C_{\varepsilon , \varphi } $ such that      
\begin{equation} \label{eq:approx_approx_s_h[f]}   
|| s_h [f] - f ||_A \leq (1 + \varepsilon ) \overline{l }_{\kappa } (\varphi ) h^{\kappa } ||\widehat{f } ||^{\circ } _{1, \kappa } + C_{\varepsilon} h^s || \widehat{f } ||^{\circ } _{1, s } .   
\end{equation}   
\end{theorem}   
\medskip   
   
\noindent 
Concretely, the condition on $\widehat{f } $ means that   
$$   
\int _{|\xi | \leq 1 } |\widehat{f } | \, |\xi |^{\kappa } d\xi + \int _{|\xi | \geq 1 } |\widehat{f } | \, |\xi |^s d\xi < \infty .   
$$   
The theorem implies that if $\overline{l }_{\kappa } (\varphi ) $ is very small, and $\widehat{f } \in L^1 _s (\mathbb{R } ^ n ) $ with $s > \kappa $ 
then the rate of convergence for small, but not too small $h $'s will at first appear to be $h^s \ll h^{\kappa } $, up to the point that 
the first term dominates and the error saturates at a level comparable to $\overline{l }_{\kappa } (\varphi ) h^{\kappa } . $ This is the phenomenon of {\it approximate approximation} which was discovered 
by Maz'ya \cite{M_94} in the context of quasi-interpolation: see also Maz'ya and Schmidt \cite{M_S}. The quasi-interpolants these authors consider are, in our notation, $\sum _j f(jh) \varphi _h (x - jh ) $ with $\varphi (x) $ of the form $\phi (x/c ) $, where $\phi $ is rapidly decreasing (for example, a Gaussian) and $c > 0 $ is called a shape parameter. Since we are in the case of $\kappa = 0 $, their quasi-interpolants will not converge to $f(x) $, but it is shown in \cite{M_S} that if $\phi $ is smooth, satisfies certain moment conditions and decays sufficiently rapidly at infinity,  and if $f $ has bounded derivatives up till order $L $, then by choosing $c $ sufficiently large one can achieve an apparent order of convergence of $h^L $ up to a small saturation error which goes to 0 as $c $ tends to infinity. This should be compared with theorem \ref{thm:approx_approx_RBF} when $\kappa = 0 $, in which case there will also be no actual convergence and where the required smoothness of $f $ is formulated in terms of its Fourier transform. Of course, this theorem concerns the exact interpolants instead of the quasi-interpolants. We will use shape parameters below to construct basis functions $\varphi \in \mathfrak{B }_{\kappa , N } (\mathbb{R }^n ) $ with small $\overline{l }_{\kappa } (\varphi ) . $ 
 
\begin{versionA}   
\medskip   
   
\noindent \textcolor{blue}{N.B. If $\widehat{f } \in L^1 _s (\mathbb{R }^n ) $, then $f \in C^{\lfloor s \rfloor } (\mathbb{R }^n ) $ and the conclusion of theorem \ref{thm:approx_approx_RBF} for such $f $'s is stronger that that implied by Maz'ya and Schmidt's result in \cite{M_S}, theorem 2.17. Also note that theorem \ref{thm:approx_approx_RBF} we get true convergence of order $\kappa $ and apparent higher order convergence if $\widehat{f } $ decays sufficiently rapidly.} 
\medskip   
   
\end{versionA}   
We will encounter similar approximate approximation phenomena when studying convergence rates of RBF schemes in sections 5 and 6 below.      
\medskip   
   
\noindent {\it Proof of theorem \ref{thm:approx_approx_RBF}}. It suffices, by (\ref{eq:proof_convergence_RBF}), to bound $|| (1 - \widehat{L }_1 (h \xi ) ) \widehat{f }(\xi ) ||_1 . $ Let $\overline{l } := \overline{l }_{\kappa } (\varphi ) . $ Then if $\varepsilon > 0 $, there exists a $\rho (\varepsilon ) > 0 $ such that if $h |\xi | < \rho (\varepsilon ) $, then   
$0 \leq 1 - \widehat{L }_1 (h\xi ) \leq \frac{1 }{2 } (1 + \varepsilon ) \overline{l } \cdot h^{\kappa } |\xi |^{\kappa } $, and   
\begin{eqnarray*}     
&&|| (1 - \widehat{L }_1 (h \xi ) ) \widehat{f }(\xi ) ||_1 \\   
&\leq & \int _{|h \xi | \leq\rho (\varepsilon )  } \, (1 - \widehat{L }_1 (h \xi ) ) |\widehat{f }(\xi ) | \, d\xi + 2 \int _{|h \xi | \geq \rho (\varepsilon ) } \, |\widehat{f }(\xi ) | \, d\xi \\   
&\leq & \frac{1 }{2 } (1 + \varepsilon ) \overline{l } \, h^{\kappa } \int _{|\xi | \leq \rho (\varepsilon ) / h } |\xi |^{\kappa } |\widehat{f } (\xi ) | \, d\xi + 2 \rho (\varepsilon )^{- s } h^s \int _{|\xi | > \rho (\varepsilon ) / h } \, |\widehat{f } (\xi ) | \, |\xi |^s  \, d\xi ,      
\end{eqnarray*}   
which implies the theorem. 
\hfill $\Box $   
   
\begin{corollary} \label{cor:approx_approx_RBF} If $\widehat{f } $ satisfies the conditions of theorem \ref{thm:approx_approx_RBF}, 
then   
\begin{equation}   
\limsup _{h \to 0 } h^{- \kappa } || s_h [f] - f ||_A \leq \overline{l }_{\kappa } (\varphi ) \int _{\mathbb{R }^n } \, |\xi |^{\kappa } |\widehat{f }(\xi ) | \, d\xi .   
\end{equation}   
\end{corollary}   
   
The next theorem complements this upper bound by the lower bound (\ref{eq:lb_approx_error}) when $\widehat{f } $ is not necessarily compactly supported.   
   
\begin{theorem} \label{thm:lower-bound_approx_error} Let $f $ satisfy the hypothesis of theorem \ref{thm:approx_approx_RBF}:   
$f \in L^{\infty } _{- p } (\mathbb{R }^n ) $ for some $p < \kappa $ and $\widehat{f } |_{\mathbb{R }^n \setminus 0 } 
\in L^1 (\mathbb{R }^n , \max \left( |\xi |^{\kappa } , |\xi |^s ) \right) $   
for some 
$s > \kappa . $   
Then     
\begin{eqnarray} \nonumber   
\underline{l }_{\kappa } (\varphi ) \int _{\mathbb{R }^n } |\xi |^{\kappa } |\widehat{f } (\xi ) | d\xi &\leq & \liminf _{h \to 0 } h^{- \kappa } || s_h [f] - f ||_A \\   
&\leq & \limsup _{h \to 0 } h^{- \kappa } || s_h [f] - f ||_A \leq \overline{l }_{\kappa } (\varphi ) \int _{\mathbb{R }^n } |\xi |^{\kappa } |\widehat{f } (\xi ) | d\xi . \nonumber   
\end{eqnarray}   
\end{theorem}   
   
\noindent {\it Proof}. We only need to establish the lower bound. If $|\xi |_{\infty } = \max _j |\xi _j | $ is the $\ell ^{\infty } $-norm on $\mathbb{R }^n $, let $Q_h = \{ \xi \in \mathbb{R }^n : |\xi |_{\infty } \leq \pi / h \} = [ - \pi / h , \pi / h ]^n $, the cube centered at 0 with sides $2\pi / h $, and let $Q_h (\ell ) = h^{-1 } \ell + Q_h . $ 
Then 
\begin{eqnarray}   
|| \widehat{s_h [f ] } - \widehat{f } ||_1 &=& \sum _{\ell } \int _{Q_h (\ell ) } \left | \sum _k \left( \widehat{L }_1 (h \xi ) - \delta _{0, k } \right) \widehat{f } (\xi + 2 \pi k / h ) \right | d\xi \nonumber \\   
&=& \sum _{\ell } \int _{Q_h } \left | \sum _k \left( \widehat{L }_1 (h \xi + 2 \pi \ell ) - \delta _{0, k } \right) \widehat{f } (\xi + 2 \pi (k + \ell ) / h ) \right | d\xi \nonumber   
\end{eqnarray}      
so that   
\begin{eqnarray}   
|| \widehat{s_h [f ] } - \widehat{f } ||_1 &\geq & \ \ \sum _{\ell } \int _{Q_h } \left | \left( \widehat{L }_1 (h \xi + 2 \pi \ell ) - \delta _{0 , -\ell } \right) \widehat{f }(\xi ) \right | d\xi \nonumber \\   
&& - \sum _{\ell } \int _{Q_h } \sum _{k \neq - \ell } \left | \left( \widehat{L }_1 (h \xi + 2 \pi \ell ) - \delta _{0 , -  k } \right) \right | \, |\widehat{f } (\xi + 2 \pi (k + \ell ) / h ) | d\xi . \label{eq:proof_lb_ae}   
\end{eqnarray}   
The double sum in the second line 
can be bounded by   
\begin{eqnarray*}   
&&\sum _{\ell } \sum _{k \neq - \ell , 0 } \, \int _{Q_h } \left | \widehat{L }_1 (h \xi + 2 \pi \ell )\widehat{f } (\xi + 2 \pi (k + \ell ) / h ) \right | d\xi \\   
&& + \sum _{\ell \neq 0 } \int _{Q_h } \left | \left( \widehat{L }_1 (h \xi + 2 \pi \ell ) - 1 \right) \widehat{f } (\xi + 2 \pi (k + \ell ) / h ) \right | d\xi \\   
&\leq & \left( \sum _{\ell } \frac{C }{(1 + |\ell | )^N } + 2 \right) \int _{|\xi |_{\infty } \geq \pi / h } \, |\widehat{f }(\xi ) | \, d\xi \\   
&\leq & C h^s \int _{\mathbb{R }^n } \, \widehat{f }(\xi ) | \, |\xi |^s \, d\xi , 
\end{eqnarray*}   
where we used that       
$$   
\sup _{Q_h } | \widehat{L }_1 (h \xi + 2 \pi \ell ) | = \sup _{\eta \in Q_1 } \widehat{L }_1 (\eta + 2 \pi \ell ) \leq \frac{C }{(1 + |\ell | )^N } .   
$$   
The first line of (\ref{eq:proof_lb_ae}), on account of $\widehat{L }_1 $ taking values in $[0, 1 ] $, equals   
$$   
\int _{Q_h } \left( (1 - \widehat{L }_1 (h \xi ) + \sum _{\ell \neq 0 } \widehat{L }_1 (h \xi + 2 \pi \ell ) \right) |\widehat{f }(\xi ) | \, d\xi = 2 \int _{\mathbb{R }^n } \, (1 - \widehat{L }_1 (h\xi ) ) |\widehat{f } | \mathbf{1 }_{Q_h } \, d\xi ,   
$$   
where $\mathbf{1 }_{Q_h } $ is the indicator function of $Q_h . $ Since $\mathbf{1 }_{Q_h } \to 1 $ as $h \to 0 $, the lower bound now follows once more by Fatou's lemma and the definition of $\underline{l }_{\kappa } (\varphi ) . $ \hfill $\Box $   
\medskip   
   
The next proposition gives a simple explicit formula for $\underline{l }_{\kappa } $ and $\overline{l }_{\kappa } : $      
   
\begin{proposition} \label{prop:formula_l_kappa}   
For $\varphi \in \mathfrak{B }_{\kappa , N } (\mathbb{R }^n ) $ with $\kappa \geq 0 $ and $N > n $, let 
\begin{equation} \label{eq:def_A(varphi)}   
\underline{A } = \underline{A }(\varphi ) := \liminf _{\eta \to 0 } |\eta |^{\kappa } \widehat{\varphi }(\eta ), \ \ 
\overline{A } :=   
\overline{A }(\varphi ) := \limsup _{\eta \to 0 } |\eta |^{\kappa } \widehat{\varphi } (\eta ) .   
\end{equation}   
%
Then if $\kappa > 0 $,   
\begin{equation} \label{eq:formula_l_kappa}   
\overline{l }_{\kappa } (\varphi ) = \frac{2 }{\underline{A } } \sum _{k \neq 0 } \widehat{\varphi } (2 \pi k ) , \ \ \underline{l }_{\kappa } (\varphi ) = \frac{2 }{\overline{A } } \sum _{k \neq 0 } \widehat{\varphi } (2 \pi k ) ,   
\end{equation}   
while if $\kappa = 0 $,   
\begin{equation} \label{eq:formula_l_kappa=0}   
\overline{l }_0 (\varphi ) = \frac{ 2 \sum _{k \neq 0 } \widehat{\varphi } (2 \pi k ) }{\underline{A } + \sum _{k \neq 0 } \widehat{\varphi }(2 \pi k ) } , \ \ \underline{l }_0 (\varphi ) = \frac{ 2 \sum _{k \neq 0 } \widehat{\varphi } (2 \pi k ) }{\overline{A } + \sum _{k \neq 0 } \widehat{\varphi }(2 \pi k ) } .   
\end{equation}   
\end{proposition}   
   
\noindent Note that $\underline{A }(\varphi ) > 0 $ by  definition \ref{def:wBuhmann class} (iii) and that the series in these formulas converges absolutely since $N > n . $   
\medskip   
   
\noindent {\it Proof.} If $L_1 = L_1 (\varphi ) $ is the Lagrange function associated to $\varphi $, then   
$$   
0 \leq 1 - \widehat{L }_1 (\eta ) = \frac{ \sum _{k \neq 0 } \widehat{\varphi } (\eta + 2 \pi k ) }{ \sum _k \widehat{\varphi } (\eta + 2 \pi k ) } .   
$$   
If we let $R(\eta ) := \sum _{k \neq 0 } \widehat{\varphi } (\eta + 2 \pi k ) $, then $R $ is continuous (even $C^{\lfloor \kappa \rfloor + n + 1 } $) in a neighborhood of 0. Since     
$$   
\frac{1 - \widehat{L }_1 (\eta ) }{|\eta |^{\kappa } } = \frac{R(\eta ) }{|\eta |^{\kappa } \widehat{\varphi }(\eta ) + |\eta |^{\kappa } R(\eta ) } ,   
$$   
(\ref{eq:formula_l_kappa}) and (\ref{eq:formula_l_kappa=0}) follow upon letting $\eta \to 0 . $ \hfill $\Box $   
   
\begin{corollary} \label{corollary:approx_RBF} If $\lim _{\eta \to 0 } |\eta |^{- \kappa } \widehat{\varphi } (\eta ) $ exists, then $\underline{l }_{\kappa } (\varphi ) = \overline{l }_{\kappa } (\varphi ) = l_{\kappa } (\varphi ) $, and   
\begin{equation}   
\lim _{h \to 0 } h^{- \kappa } || s_h [f ] - f ||_A = l_{\kappa } (\varphi ) \, || f ||^{\circ }_{1 , \kappa } ,   
\end{equation}   
for $f $ as in theorem \ref{thm:approx_approx_RBF} with $s > \kappa . $   
\end{corollary}   
   
\noindent We can often construct basis functions with small $\overline{l } _{\kappa } (\varphi ) $ by introducing a so-called {\it shape-parameter} $c $ and taking $\varphi $ of the form $\varphi (x) = \phi (x/c ) := \phi _c (x) $ with $c $ large, for suitable $\phi \in \mathfrak{B }_{\kappa , N } (\mathbb{R }^n ) $:   
   
\begin{proposition} Suppose that $\phi \in \mathfrak{B }_{\kappa , N } (\mathbb{R }^n ) $ with $\kappa \geq 0 $ and $N > \max (\kappa , n ) . $ Then $\lim _{c \to \infty } \overline{l }_{\kappa } (\phi _c ) = 0 . $   
\end{proposition}     
   
\noindent {\it Proof}. Since $\widehat{\phi }_c (\eta ) = c^n \widehat{\phi } (c \eta ) $, it follows that $\underline{A }(\phi _c ) = c^{n - \kappa } \underline{A }(\phi ) $, and therefore, by (\ref{eq:formula_l_kappa}), if $\kappa > 0 $,    
\begin{equation} \label{eq:l_phi_c}   
\overline{l }_{\kappa } (\phi _c ) = 2 \frac{c^{\kappa } }{\underline{A }(\phi ) } \sum _{k \neq 0 } \widehat{\phi } (2 \pi c k ) \leq C \, c^{\kappa - N } \sum _{k \neq 0 } |k |^{-N } ,   
\end{equation}   
which tends to 0 as $c \to \infty $ under the stated conditions on $\kappa . $ The case of $\kappa = 0 $ follows by observing that 
$$   
\overline{l }_0 (\phi _c ) \leq \frac{2 }{\underline{A }(\phi _c ) } \sum _{k \neq 0 } \widehat{\phi _c }(2 \pi k ) ,   
$$   
and proceeding as before. \hfill $\Box $   
\medskip   
   
\noindent Examples of basis functions $\phi $ which satisfy the conditions of the corollary are the Gaussians (for which $\kappa = 0 $) and the generalized multiquadrics, whose Fourier transforms decay exponentially at infinity, but none of the homogeneous basis functions, since for these $\kappa = N $: see examples \ref{examples_multiquadric} below for more discussion. In fact, for a Gaussian or a multiquadric, $\widehat{\phi }(\xi ) $ decays expontially at infinity, and $\overline{l }_{\kappa } (c) $ will decay exponentially in $c . $   
   
\medskip   
   
   
\subsection{Dependence of the constants on the shape parameter} \label{subsection:c-dependence} If we want to apply (\ref{eq:approx_approx_s_h[f]}) with $\varphi = \phi _c $ as above (and some fixed $\varepsilon $), it becomes interesting to ask how the constant $C_{\varepsilon } = C_{\varepsilon , \varphi } $ depends on the shape-parameter $c . $ On expects it to go to infinity with $c $, and the question then is at what rate exactly. We will see that under reasonable assumptions on $\widehat{\phi } (\eta ) $, this constant behaves like $c^s . $   
   
We first analyse its dependence on $\varphi $ for a general $\varphi \in \mathfrak{B }_{\kappa , N } (\mathbb{R }^n ) . $ From the proof, $C_{\varepsilon , \varphi } $  is proportional to $\rho (\varepsilon , \varphi ) ^{-s } $, where $\rho (\varepsilon , \varphi ) $ is any positive number such that $2 | \eta |^{- \kappa } (1 - \widehat{L }_1 (\eta ) ) < (1 + \varepsilon ) 
\overline{l }_{\kappa } $ for $|\eta | \leq \rho (\varepsilon ) . $ If we put 
$$   
R(\eta ) := R_{\varphi } (\eta ) = \sum _{k \neq 0 } \widehat{\varphi } (\eta + 2 \pi k ) ,   
$$   
then   
$$   
0 \leq \frac{1 - \widehat{L }_1 (\eta ) }{|\eta |^{\kappa } } = \frac{R(\eta ) }{|\eta |^{\kappa } ( \widehat{\varphi } (\eta ) + R(\eta ) ) } \leq \frac{R(\eta ) }{ |\eta |^{\kappa } \widehat{\varphi } (\eta ) }   
$$   
If we introduce the family of semi-norms    
$$   
p_r (\psi ) := \sum _{k \neq 0 } \sup _{|\eta | \leq r } | \psi (\eta + 2 \pi k ) | , \ \ r > 0 ,   
$$   
then 
$| R(\eta ) - R(0 ) | \leq p_r (|\nabla \widehat{\varphi } | ) \,  |\eta | $ for $|\eta | \leq r $,    
where $|\nabla \varphi | $ the euclidean norm of the gradient, so that $0 \leq 2 R(\eta ) < 2 (1 + \varepsilon ) R(0) $ if $|\eta | \leq \rho _1 $, where   
\begin{equation} \label{eq:estimate_rho_1}   
\rho _1 
:= \rho _1 (\varepsilon ; \varphi ) := \min \left( r , \varepsilon R_{\varphi } (0) / p_r (|\nabla \widehat{\varphi } | \right) .   
\end{equation}   
On the other hand, there exists a $\rho _2 := \rho _2 (\varepsilon ; \varphi ) $ such that $ | \eta |^{\kappa } \widehat{\varphi } (\eta ) \geq (1 - \varepsilon ) \underline{A }(\varphi ) $ if  $|\eta | \leq \rho_2 $, where $\underline{A } = \underline{A } (\varphi ) $ was defined in proposition \ref{prop:formula_l_kappa} above.   
Combining these two estimates and observing that $\overline{l }_{\kappa } = 2 R(0) / \underline{A } $ we see that   
$$   
0 \leq \frac{1 - \widehat{L }_1 (\eta ) }{ |\eta |^{\kappa } } \leq \left( \frac{1 + \varepsilon }{1 - \varepsilon } \right) \overline{l }_{\kappa }   
$$   
if $|\eta | \leq \rho $ where   
\begin{equation} \label{eq:estimate_rho}   
\rho:= \rho (\varepsilon ; \varphi ) := \min (\rho_1 , \rho_2 ) ,   
\end{equation}   
and we get an approximate approximation estimate in the form   
\begin{equation} \label{eq:approx_approx_expl}   
|| s_h [f] - f ||_A \leq \frac{1 + \varepsilon }{1 - \varepsilon } \overline{l }_{\kappa } (\varphi ) h^{\kappa } ||\widehat{f } ||^{\circ } _{1, \kappa } + \frac{h^s }{\rho ( \varepsilon , \varphi )^s } || \widehat{f } ||^{\circ } _{1, s } .   
\end{equation}   
We now 
derive lower bounds for $\rho $ when $\varphi = \phi _c $ 
for a fixed $\varepsilon $, which we will sometimes drop from the notations. 
First of all, there exists a $r(\phi ) = r(\phi ; \varepsilon ) $ such that   
$\widehat{\phi } (\eta ) \geq (1 - \varepsilon ) \underline{A } (\phi ) $ if $|\eta | \leq r(\phi ) . $ This implies that $|\eta |^{\kappa } \widehat{\phi _c } (\eta) = c^n |\eta |^{\kappa } \widehat{\phi } (c \eta ) \geq c^{n - \kappa } \frac{1 }{2 } \underline{A } (\phi ) = \frac{1 }{2 } \underline{A } (\phi _c ) $ if $c |\eta | \leq r(\phi ) $, so we can take $\rho_2 (\varepsilon ; \phi _c ) = r(\phi ) / c . $ 
   
The behavior of $R_{\phi _c } (0) / p_{\pi, N }( \nabla {\widehat{\phi _c } } ) $ for large $c $ depends on the asymptotic behavior at infinity of $\widehat{\phi } $ and of $\nabla \widehat{\phi } . $ We consider the case of a polynomially decaying $\widehat{\phi } $ and that of an exponentially decaying one.   
   
\begin{example} (Polynomially decaying $\widehat{\phi } $) \rm{Suppose that   
\begin{equation} \label{eq:asymp_phi}   
\inf _{|\eta | \geq 2 \pi } |\eta |^N \widehat{\phi }(\eta ) =: a > 0 .      
\end{equation}   
Then   
\begin{equation}   
R_{\phi _c } (0) \geq a c^{n - N } \sum _{k \neq 0 } (2 \pi |\eta | ) ^{-N } .   
\end{equation}   
On the one hand, we have for 
any $\varphi \in \mathfrak{B }_{\kappa , N } (\mathbb{R }^n ) $ that   
$$   
p_{\pi } (\nabla \widehat{\varphi } ) \leq \left( \sum _{k \neq 0 } (2 \pi (|k | - \tfrac{1 }{2 } ) ^{-N } \right) \, \sup _{|\eta | \geq \pi } |\eta |^N |\nabla \widehat{\varphi }(\eta ) | ,   
$$   
which implies that   
$$   
p_{\pi } ( |\nabla \widehat{\phi }_c | ) \leq C \cdot c^{n + 1 - N } \sup _{|\eta | \geq c \pi } |\eta |^N |\nabla \widehat{\phi } (\eta ) |
$$   
which can certainly be bounded from above by a constant times $c^{n + 1 - N } $ if $c \geq 1 . $ 
Taking $r = \pi $ in (\ref{eq:estimate_rho_1}) we therefore find 
that $\rho _1 (\varepsilon , \phi _c ) \geq C c^{-1 } $ for some constant $C . $  It follows that $\rho (\varepsilon , \phi _c ) \geq C c^{-1 } $ for some 
constant $C $ and 
(\ref{eq:approx_approx_expl}) implies an estimate   
$$   
|| s_h [f] - f ||_A \leq C_1 c^{\kappa - N } h^{\kappa } ||\widehat{f } ||^{\circ } _{1, \kappa } + C_2 ^s c^s h^s || \widehat{f } ||^{\circ } _{1, s } ,   
$$   
with constants $C_1 $ and $C_2 $ which depend on $\phi $ (and on $\varepsilon $), 
but not on the shape parameter $c . $ Neglecting the numerical values of these constants, the second term will dominate the first as long as $h \gg c^{-1 - N/(s - \kappa ) } $, 
and the estimate will be relevant for $h $'s in the range $c^{-1 - N/(s - \kappa ) } \ll h \ll  c^{-1 } $ (since we want $ (ch)^s \ll 1 $). 
}   
\end{example}   
   
\begin{example} (Exponentially decaying $\widehat{\phi } $) {\rm Suppose that $\phi \in \mathfrak{B }_{\kappa , N } (\mathbb{R }^n ) $ is such that there exists constants $p \geq 0 $ and $C > 0 $ for which 
%
%
\begin{equation} \label{eq:phi_exp_dec}      
C^{-1 }  |\eta |^{-p } e^{- |\eta | } \leq \widehat{\phi } (\eta ) \leq C |\eta |^{-p } e^{- |\eta | } , \ \ |\eta | \geq 1 ,        
\end{equation}    
and   
$$   
|\nabla \widehat{\phi } | \leq C |\eta |^{-p } e^{- |\eta | } , \ \ |\eta | \geq 1 .   
$$   
An example of such a $\phi $ is given by the multi-quadric on $\mathbb{R }^n $: see Example \ref{examples_multiquadric} below.   
   
Since $\widehat{\phi _c } (\eta ) = c^n \widehat{\phi } (c \eta ) $, $p_r (|\nabla \widehat{\phi _c } | ) $ will then be bounded by a constant times 
$$   
c^{n + 1 } \sum _{k \neq 0 } \sup _{|\eta | \leq r } \frac{e^{ - c | \eta + 2 \pi k | } }{ | c \eta + 2 \pi c k |^p }   
$$   
Using that if $|\eta | \leq r $, $e^{ - c | \eta + 2 \pi k | } \leq e^{c r } e^{- 2 \pi c |k | } $ and, assuming wlog that $r \leq \pi $, that $|\eta + 2 \pi k | \geq 2 \pi (|k | - \frac{1 }{2 } ) \geq \pi |k | $ for $|k | \geq 1 $, we find that this expression is bounded by   
$$   
c^{n + 1 } e^{c r } \sum _{k \neq 0 } \frac{e^{- 2 \pi c |k | } }{ (\pi c |k | )^p } ,   
$$   
and using the first inequality of (\ref{eq:phi_exp_dec}), it follows that $p_r (|\nabla \widehat{\phi _c } | ) $ is bounded by a constant times $c R_{\phi _c } (0) . $ Remembering (\ref{eq:estimate_rho_1}), it follows that $\rho _1 (\phi _c ) $ is bounded from below by a constant times $\min ( \pi , r , c^{-1 } e^{- c r } ) . $ This still has $r $ as a free parameter, and taking $r = c^{-1 } $ and remembering the estimate for $\rho _2 (\phi _c ) $ from the previous example, we find that we can take $\rho (\phi _c ) \geq C \cdot c^{-1 } . $ Since (\ref{eq:phi_exp_dec}) implies that $\underline{l }_{\kappa } \leq C \cdot c^{\kappa - p } e^{- 2 \pi c } $, we now have the approximate approximation estimate   
$$   
|| s_h [f] - f ||_A \leq C_1 c^{\kappa - p } e^{- 2 \pi c } h^{\kappa } ||\widehat{f } ||^{\circ } _{1, \kappa } + C_2 ^s c^s h^s || \widehat{f } ||^{\circ } _{1, s } .   
$$   
which is relevant for the range $c^{- 1 - p / (s - \kappa ) } e^{- 2 \pi c / (s - \kappa ) } \ll h \ll c^{-1 } $, where $c \gg 1 . $   
\medskip   
   
\begin{example} (The Gaussian) \rm{A final interesting example is that of the Gaussian, with $\widehat{\phi } (\eta ) = e^{- |\eta |^2 } . $ In this case, $\nabla \widehat{\phi _c } (\eta ) = c^{n + 2 } \eta e^{- c^2 |\eta |^2 } $, and $p_r (|\nabla \widehat{\phi _c } | ) $ for $r \leq \pi $ is bounded by a constant times   
$$   
c^{n + 2 } \sum _{k \neq 0 } |k | \sup _{|\eta | \leq r } e^{ - c^2 | \eta + 2 \pi k |^2 } \leq c^{n + 2 } e^{(\epsilon ^{-1 } - 1 ) c^2 r^2 } \cdot \sum _{k \neq 0 } e^{ - c^2 (1 - \epsilon ) |2 \pi k |^2 } ,   
$$   
for any $0 < \epsilon < 1 $, where we used the inequality $(a - b )^2 \geq (1 - \epsilon ) a^2 - (\epsilon ^{-1 } - 1 ) b^2 $ with $a = 2 \pi |k | $ and $b = |\eta | . $ Asymptotically, for large $c $, this behaves like $c^{n + 2 } e^{ - c^2 (1 - \epsilon ) 4\pi ^2 + (\epsilon ^{-1 } - 1 ) c^2 r^2 } $, while $R_{\phi _c } (0) c^{n + 2 } \simeq e^{- 4 \pi ^2 c^2 } $, and we can choose $\rho _1 (\phi _c ) $ (with a fixed $\varepsilon $) to be bounded from below by a constant times   
$$   
\min \left( \pi , r , c^{-2 } e^{- 4 \pi ^2 \epsilon c^2 } e^{(\epsilon ^{-1 } - 1 ) c^2 r^2 } \right) .   
$$   
If we now take the free parameters $r $ and $\epsilon $ equal to $c^{- 2 } $ and remember that $\rho _1 (\phi _c ) \simeq c^{-1 } $, we conclude that we can choose $\rho (\phi _c  ) \simeq c^{-2 } $, and we have the following approximate approximation estimate for Gaussian RBF interpolation (recalling that $\kappa = 0 $ in this case):      
$$   
|| s_h [f] - f ||_A \leq C_1 e^{- 4 \pi ^2 c^2 } ||\widehat{f } ||_1 + C_2 ^s c^{2s } h^s || \widehat{f } ||^{\circ } _{1, s } ,   
$$  
where the constant $C_1 $ can be taken arbitrarily close to 1.   
}   
\end{example}

}   
\end{example}   
   

\section{\bf Convergence of stationary RBF schemes for 
the heat equation}   
   
\subsection{An RBF scheme for the heat equation} We introduce an RBF scheme for the 
Cauchy problem for the classical heat equation, 
\begin{equation} \label{eq:heat}
\left \{ \begin{array}{ll}
\partial _t u (x, t ) = \Delta u (x, t ) , \ x \in \mathbb{R }^n , t > 0 \\
u(x, 0 ) = f(x) ,
\end{array}
\right.
\end{equation}
$\Delta = \sum _{j = 1 } ^n  \partial _{x_j } ^2 $ being the Laplace operator, and examine its convergence. The scheme is a variant of the classical method of lines, and looks for approximate solutions $u_h $ of the form   
\begin{equation} \label{eq:def_u_h}   
u _h (x, t ) = \sum _{k \in \mathbb{Z }^n } c_k (t; h ) L_1 (h^{-1 } x - k ) ,   
\end{equation}
where the $c_k (\cdot ; h ) : [0, \infty ) \to \mathbb{R } $ are differentiable functions. Here, $L_1 $ is the Lagrange interpolation function of theorem \ref{thm_L_1}, associated to a given basis function $\varphi \in \mathfrak{B }_{\kappa , N } (\mathbb{R }^n ) $ which we fix. {\it We assume throughout this section that $N > n + 2 $ and $\kappa > 0 . $} We will see below that for the scheme to converge  we will need that $\kappa > 2 $ while for $\kappa = 2 $ we can have approximate convergence. The coefficients $c_k (t; h ) $ of $u_h $ are determined by requiring that $u_h $ solve (\ref{eq:heat}) exactly in the points of $h \mathbb{Z }^n $:
\begin{equation} \label{eq:heat_1}   
\partial _t u_h (hj , t ) = \Delta u_h (jh , t ) , \ \ \forall j \in \mathbb{Z }^n ,   
\end{equation}   
while $u_h (x, 0 ) $ is taken to be equal to $s_h [f ] (x) $, the RBF-interpolant of $f . $ 
This leads to the following initial value problem for the coefficients $c_j (t; h ) $:
\begin{equation} \label{eq:ODE-system}
\left \{ \begin{array}{ll}
\displaystyle{\frac{d c_j  }{dt } (t; h ) = h^{-2 } \sum _k \Delta L_1 (j - k ) c_k (t; h ) } \\
c_j (0; h ) = f(jh ) .
\end{array}
\right.
\end{equation}   
Since this is an infinite system of ODEs we first discuss existence and uniqueness of solutions 
in suitable Banach spaces.   
\medskip   
   
For $s \in \mathbb{R } $, let      
\begin{equation} \nonumber   
\ell ^{\infty } _{s } := \ell ^{\infty } _{s } (\mathbb{Z }^n ) := \{ (c_j )_{j \in \mathbb{Z }^n } : || c ||_{\infty , s } := \sup _j (1 + |j | )^{ s } |c_j | < \infty \} .   
\end{equation}   
One easily verifies 
that the convolution operator   
$$   
A := A_L : (c_j )_j \to \left( \sum _k \Delta L_1 (j - k ) \, c_k \right)_j ,   
$$   
is a bounded operator on $\ell ^{\infty } _{-p } $ if $0 \leq p < \kappa . $   
%
Indeed,          
\begin{eqnarray*}   
(1 + |j | )^{-p } \left | \sum _k \Delta L_1 (j - k ) c_k \right | &\leq & \sum _k \left( (1 + |j | )^{- p } | (1 + |k | )^p |\Delta L_1 (j - k ) | \right) || c ||_{\infty , - p } \\   
&\leq & \left( \sum _k (1 + |j - k | )^p |\Delta L_1 (j - k ) | \right) || c ||_{\infty , - p } ,   
\end{eqnarray*}   
using that $(1 + |k| ) \leq (1 + |j - k | )( 1 + |j | ) . $ The sum of the series on the right is independent of $j $ and finite if $p < \kappa $, by theorem \ref{lemma_L_1}.   
%
\noindent It follows that if we let $c(t) := (c_j (t) )_j $, and if the initial value $c(0) \in \ell ^{\infty } _{- p } $ for some $p \in [0 , \kappa ) $, the system (\ref{eq:ODE-system}) 
has a unique $\ell ^{\infty } _{- p } $-valued solution which is given by $c(t) = e^{h^{-2 } t A_L } (c(0) ) . $   
%
Next, if for $c \in \ell ^{\infty } _{-p } $ we let (with some abuse of notation)   
$$   
s_h [c](x) := \sum _{j \in \mathbb{Z }^n } c_j L_1 (h^{-1 } x - j ) ,  
$$   
then $s_h : c \to s_h [c] $ is a bounded linear operator from $\ell ^{\infty } _{-p } \to L^{\infty } _{-p } (\mathbb{R } ) $ if $0 \leq p < \kappa . $   
Indeed, using the decay of $L_1 $,   
\begin{eqnarray*}   
\frac{ || s_h [c ] ||_{\infty, -p } }{||c ||_{\infty , -p } } &\leq& \sup _{x \in \mathbb{R }^n } (1 + |x| )^{-p } \sum _j \frac{(1 + |j | )^p }{(1 + |h^{-1 } x - j | )^{\kappa + n } } \\   
&=& \sup _{y \in \mathbb{R }^n } (1 + |hy |)^{-p } \sum _j \frac{(1 + |j | )^p }{(1 + |y - j | )^{\kappa + n } } \\   
&\leq & \sup _{y \in \mathbb{R }^n } \left( \frac{1 + |y | ) }{1 + h|y | } \right)^p \sum _j \frac{(1 + |y - j | )^p }{(1 + |y - j | )^{\kappa + n } } .   
\end{eqnarray*}   
The sum on the right converges and defines a 1-periodic continuous function on $\mathbb{R }^n $ which is therefore uniformly bounded, while the factor in front can be estimated by $\max (1 , h^{-n } ) . $   
\medskip   
   
If $f \in L^{\infty } _{-p } (\mathbb{R }^n ) $, we can in particular take $c(0) = f |_{h \mathbb{Z }^n } $, and  
\begin{equation} \label{eq:def_u_h_bis}   
u_h [f ] (x , t ) := s_h \left[ e^{h^{-2 } t A_L } \left( f |_{h \mathbb{Z }^n } \right) \right] (x) ,   
\end{equation}
is the unique function (\ref{eq:def_u_h}) whose coefficients satisfy (\ref{eq:ODE-system}). We summarize this discussion in the following lemma:     
   
\begin{lemma} \label{lemma:ODE-system} If $0 \leq p < \kappa $ and if $f \in L^{\infty } _{-p } (\mathbb{R }^n ) $ then there is a unique function $u_h = u_h [f ] \in C^1 \left( [0 , \infty ) ; L^{\infty } _{-p } (\mathbb{R }^n ) \right) $ of the form (\ref{eq:def_u_h}) which satisfies (\ref{eq:heat_1}) 
and      
$f \to u_h [f ] (\cdot , t ) $ is a bounded linear map on $L^{\infty } _{-p } (\mathbb{R } ) $, for each $t \geq 0 $ and $h > 0 . $  
\end{lemma}   
   
\noindent In particular, for each fixed $t $, $u_h (x, t ) $ has tempered growth in $x $, and thus possesses a well-defined Fourier transform, which we will study next.   
\medskip   
   
\begin{versionA}   
\noindent \textcolor{blue}{{\bf Comment}: The case $\kappa = 0 $, which may arise if we consider evolution equations with a  $\psi $DO of order 0 on the right hand side, would need a separate discussion if we'd want to study approximate convergence for these. If $\kappa = 0 $, we could only prove for the moment that if $q > $ (corresponding to $p = - q < 0 $), then   
$$   
A_L : \ell _{q , \infty } \to \ell _{0 , \infty } ,   
$$   
and the series for $s_h \left[ e^{h^{-2 } t A_L } \left( f |_{h \mathbb{Z }^n } \right) \right] (x) $ does not a-priori converge (though it might be summable?).   
}   
\end{versionA}   
     
\subsection{Convergence of the scheme in Wiener norm}   
We 
start by computing the Fourier transform of $u_h = u_h [f ] . $ Let us introduce the auxiliary function $G (\eta ) $ on $\mathbb{R }^n $ by
\begin{eqnarray} \label{eq:def_G}   
G(\eta ) := G_{\varphi } (\eta ) &:= & \sum _k |\eta + 2 \pi k |^2 \widehat{L }_1 (\eta + 2 \pi k ) \\  
&=& \frac{\sum _k |\eta + 2 \pi k |^2 \widehat{\varphi } (\eta + 2 \pi k ) }{\sum _k \widehat{\varphi }(\eta + 2 \pi k ) } \nonumber   
\end{eqnarray}   
where the series converges absolutely, since $N > n + 2 . $    
   
\begin{lemma} \label{lemma:u_h_hat} Let $\varphi \in \mathfrak{B }_{\kappa , N } (\mathbb{R } ) $ with $N > n + 2 $ and $\kappa > 0 . $ If 
$\widehat{f } \in L^1 (\mathbb{R }^n ) $, then the Fourier transform of $u_h (x, t ) $ with respect to $x $ is given by   
\begin{equation} \label{u_h_hat}
\widehat{u }_h (\xi , t ) = e^{- t h^{-2 } G(h \xi ) } \widehat{s_h [f] }(\xi ) . 
\end{equation}
\end{lemma}   
   
\noindent {\it Proof.} Since $( \Delta L_1 (j ) )_{j \in \mathbb{Z }^n } $ is in $\ell ^1 := \ell ^1 (\mathbb{Z }^n ) $, it follows that $A_L $ is a bounded operator on $\ell ^1 $, and hence the system (\ref{eq:ODE-system}) has an $\ell ^1 $-valued solution   $c(t) = (c_j (t) )_j $ if the initial value $c(0) \in \ell _1 . $ 
In particular, if $c(0) = f |_{h \mathbb{Z }^n } $ with $f \in \mathcal{S } (\mathbb{R }^n ) $, then the function $(x, k ) \to c_k (t ; h ) L_1 ( h^{-1 } x - k ) $ is absolutely integrable on $\mathbb{Z }^n \times \mathbb{R }^n $ for each fixed $t \geq 0 , h > 0 $, and an application of Fubini's theorem shows that 
the Fourier transform of $u_h ( \cdot , t ) $ is given by $\widehat{u }_h (\xi , t ) = h^n \widehat{L }_1 (h \xi ) \, \gamma _h (\xi , t ) $, where   
$$
\gamma _h (\xi , t ) := \sum _{j \in \mathbb{Z } } c_j (t; h ) e^{- i h (j, \xi ) } ,
$$
the series being absolutely convergent. By (\ref{eq:ODE-system}),   
\begin{eqnarray*}
\partial _t \gamma _h (\xi , t ) &=& h^{-2 } \left( \sum _j \Delta L_1 (j ) e^{- i h (j , \xi ) } \right) \gamma _h (\xi , t ) \\
&=& h^{-2 } \left( \sum _k \widehat{\Delta L }_1 (h \xi + 2 \pi k ) \right) \gamma _h (\xi , t ) \\
&=& - h^{-2 } G(h \xi ) \gamma _h (\xi , t ) ,
\end{eqnarray*}
where the second line follows from the Poisson summation formula, whose application is justified by the decay at infinity of $\Delta L_1 $ and its Fourier transform.     
Hence 
$\partial _t \widehat{u }_h (\xi , t ) = - h^{-2 } G(h \xi ) \widehat{u }_h (\xi , t ) $ which, together with the initial condition $u_h (x, 0 ) = s_h [f ](x) $ implies (\ref{u_h_hat}).   
   
If $\widehat{f } \in L^1 (\mathbb{R }^n ) $, (\ref{u_h_hat}) follows by an 
approximation argument: 
if $\widehat{f }_{\nu } \to \widehat{f } $ in $L^1 $ with $f_{\nu } $ rapidly deceasing, then $f_{\nu } \to f $ in $L^{\infty } $, so by lemma \ref{lemma:ODE-system}, $u_h [f_{\nu } ] (\cdot , t ) \to u_h [f ] (\cdot , t ) $ in $L^{\infty } $ also, since $\kappa > 0 . $ 
Hence their Fourier transforms converge in $\mathcal{S }' := \mathcal{S } ' (\mathbb{R }^n ) . $ On the other hand, $\widehat{s_h [f_{\nu } } ] = \Sigma _h (\widehat{f_{\nu } } ) \to \Sigma _h  (\widehat{f } ) = \widehat{s_h [f ] } $ in $L^1 $ and since $G $ is non-negative, $e^{- h^{-2 } t G(h \xi ) } \widehat{s_h [f_{\nu } ] } (\xi ) 
\to e^{- h^{-2 } t G(h \xi ) } \widehat{s_h [f ] } (\xi ) $ 
in $L^1 $, and therefore in $\mathcal{S } ' . $ 
\hfill $\Box $   
\medskip   
   
The following proposition lists some useful properties of $G . $   
   
\begin{proposition} \label{prop_G} Suppose that $\varphi \in \mathfrak{B }_{\kappa , N } $ with $N > n + 2 $ and let $G := G_{\varphi } $ be defined by (\ref{eq:def_G}). Then    
\medskip   
   
\noindent (i) $G $ is a positive $2 \pi $-periodic function, and $G(\eta ) = 0 $ iff $\eta \in 2\pi \mathbb{Z }^n . $   
\medskip   
   
\noindent (ii) 
There exists a constant $C > 0 $ such that $0 \leq G(\eta ) - | \eta |^2 \leq C | \eta |^{\kappa } $ for $| \eta | \leq \pi . $     
\medskip   
   
\noindent (iii) $G  $ belongs to the H\"older space $C_b ^{\lceil \kappa \rceil - 1 , \lambda }  (\mathbb{R }^n ) $ with $\lambda = \kappa - (\lceil \kappa \rceil - 1 ) . $
\end{proposition}   
   
\noindent {\it Proof}. (i) The periodicity is obvious and the positivity of $G $ is an immediate consequence of the positivity of $\widehat{\varphi } . $ 
Next, 
$G(\eta ) = 0 $ iff $|\eta + 2\pi k |^2 \widehat{L }_1 (\eta + 2 \pi k ) = 0 $ for all $k . $ Since $\widehat{L }_1 $ is non-zero outside of $(2 \pi ) \mathbb{Z }^n \setminus 0 $, this implies that $\eta \in 2 \pi \cdot \mathbb{Z }^n . $ Conversely, any such $\eta $ is a zero, given that $\widehat{L }_1 (2 \pi k ) = \delta _{0k } . $   
\medskip   
   
Assertion (ii) follows from   
\begin{eqnarray*}      
G(\eta ) - | \eta |^2 &=& \frac{ \sum _{k \neq 0 } ( | \eta + 2 \pi k |^2 - | \eta |^2 ) \widehat{\varphi } ( \eta + 2 \pi k ) }{ \sum _k \widehat{\varphi } ( \eta + 2 \pi k ) } \\   
&=& \frac{ \sum _{k \neq 0 } ( 4 \pi ^2 | k |^2 + 4 \pi (\eta , k ) ) \widehat{\varphi } ( \eta + 2 \pi k ) }{ \sum _k \widehat{\varphi } ( \eta + 2 \pi k ) }   
\end{eqnarray*}   
and the behaviour of $\widehat{\varphi } (\eta ) $ for small $\eta . $   
\medskip

Finally, (iii) follows from lemma \ref{lemma:der_L_hat}. \hfill $\Box $   
\medskip   
   
Property (iii) allows us to extend lemma \ref{lemma:u_h_hat} to functions of polynomial growth: compare lemma \ref{lemma:convergence_RBF_bis}.   
   
\begin{lemma} \label{lemma:u_h_hat_bis} Suppose that $f \in L^{\infty } _{-p } (\mathbb{R }^n ) $ for some $p < \kappa $ such that $\widehat{f } |_{\mathbb{R }^n \setminus 0 } \in L^1 \left( \mathbb{R }^n , (|\xi |^{\kappa } \wedge 1 ) d\xi \right) . $ Then the identity (\ref{u_h_hat}) holds in the sense of distributions for each $t \geq 0 . $   
Moreover, $\widehat{u_h [f ] } (\cdot , t ) - \widehat{u } (\cdot , t ) $ can be identified with the function
\begin{equation} \label{eq:proof_thm:conv_RBF_scheme}   
e^{- t h^{-2 } G(h \xi ) } \left( \, \widehat{s_h [f ] } - \widehat{f } \, \right) (\xi )+ \left(e^{- t (h^{-2 } G(h \xi ) - |\xi |^2 ) } - 1 \right) e^{ - t |\xi |^2 }\widehat{f }(\xi ) ,   
\end{equation}       
which moreover is integrable on $\mathbb{R }^n . $   
\end{lemma}   
\smallskip   
   
\noindent Here, and below, $\widehat{f } $ without argument will indicate the distribution, and $\widehat{f }(\xi ) $ the function with which it can be identified on $\mathbb{R }^n \setminus 0 $; we recall that by lemma \ref{lemma:convergence_RBF_bis}, $\widehat{s_h [f ] } - \widehat{f } $ can be identified with an $L^1 $-function. That the second term of (\ref{eq:proof_thm:conv_RBF_scheme}) is in $L^1 $ follows from lemma \ref{prop_G}(ii).   
\medskip  
   
\noindent {\it Proof}. We just clarify the statement of the lemma, and refer to Appendix \ref{Appendix_proof_thm_conv_RBF_bis} for the 
proof, which uses elements of the proof of lemma \ref{lemma:convergence_RBF_bis}. The proof of that lemma shows that $\widehat{s_h [f ] } = \Sigma _h (\widehat{f } ) $ extends to a continuous linear functional on $C_b ^{\lceil \kappa \rceil - 1 , \lambda }  (\mathbb{R }^n ) . $ Since (i) and (iii) of proposition \ref{prop_G} imply that the function $e^{- h^{-2 } t G(h \cdot ) } $ is in $C_b ^{\lceil \kappa \rceil - 1 , \lambda }  (\mathbb{R }^n ) $, its product with $\Sigma _h (\widehat{f } ) $ is well-defined as an element of the dual of $C_b ^{\lceil \kappa \rceil - 1 , \lambda }  (\mathbb{R }^n ) . $ \hfill $\Box $   
\medskip   
   
   
We can now show convergence in Wiener norm of the $u_h $ to the solution of the Cauchy problem (\ref{eq:heat}): recall again (\ref{eq:L^1_rs}). 
   
\begin{theorem} \label{thm:conv_RBF_scheme_1} Let $\varphi \in \mathfrak{B }_{\kappa , N } (\mathbb{R }^n ) $ with $N > n + 2 $ and $\kappa > 2 $ and suppose that $f \in L^{\infty } _{-p } (\mathbb{R }^n ) $ for some $p < \kappa $ such that
$$   
\widehat{f } |_{\mathbb{R }^n \setminus 0 } \in  \mathring{L }^1_{\kappa - 2 , \kappa } (\mathbb{R }^n )     
$$   
Let $u_h := u_h [f ] $ and let $u $ be the solution to the Cauchy problem (\ref{eq:heat}) with initial value $f . $ Then there exists a constant $C = C_{\varphi } $ independent of $h $ and $f $ such that for $0 < h \leq 1 $,   
\begin{equation}   
|| u_h (\cdot , t ) - u(\cdot , t ) ||_A \leq C \max (t , 1 ) \, || \widehat{f } ||^{\circ }_{\kappa - 2 , \kappa } \, h^{\kappa - 2 } .   
\end{equation}   
In particular, $u_h $ converges to $u $ in sup-norm at a rate of $h^{\kappa - 2 } . $   
\end{theorem}   
   
\noindent {\it Proof.} 
Note that the conditions on $f $ are weaker than those of theorem 
\ref{thm:convergence_RBF_bis}. Indeed, we will be applying corollary \ref{cor:convergence_RBF} with $r = \kappa - 2 . $ By lemma \ref{lemma:u_h_hat_bis}, we can estimate $|| \widehat{u }_h (\xi , t ) - \widehat{u }(\xi , t ) ||_1 $ by   
\begin{equation}  \label{eq.proof_thm_conv_RBF_scheme_1}   
\left | \left | \, e^{- h^{-2 } t G(h \xi ) } \left( \, \widehat{s_h [f] } - \widehat{f } \, \right) \, \right | \right |_1 + \left | \left | \, 
\left( 1 - e^{- t h^{-2 } ( G(h \xi ) - |h \xi |^2 ) } \right) e^{ - t |\xi |^2 } \widehat{f } \, \right | \right |_1 . 
\end{equation}   
Since $G $ is positive, the first term can be bounded by $ | | \, \widehat{s_h [f] } - \widehat{f } \, | |_1 \leq C || \widehat{f } ||_{\kappa - 2 , \kappa } h^{\kappa - 2 } $, 
by corollary \ref{cor:convergence_RBF}. To estimate the second term, 
we first use proposition \ref{prop_G}(ii) together with the inequality $(1 - e^{-x } ) \leq x $ for $x \geq 0 $ to estimate the integral over $ h | \xi | \leq \pi $ by   
\begin{eqnarray} \label{eq.proof_thm_conv_RBF_scheme_2}   
&& \ \ \ \ \ \ \int _{h|\xi | \leq \pi } \, 
\left( 1 - e^{- t h^{-2 } ( G(h \xi ) - |h \xi |^2 ) } \right) e^{ - t |\xi |^2 } \, |\widehat{f }(\xi ) | \, d\xi \\   
&\leq & C \, h^{\kappa - 2 } \int _{|\xi | \leq h^{-1 } \pi } \, t |\xi |^{\kappa } \, |\widehat{f } (\xi ) | e^{- t |\xi |^2 } \, d\xi \nonumber \\   
&\leq & C h^{\kappa - 2 } \left(  \int _{|\xi | \leq 1 } \, t |\xi |^{\kappa } \, |\widehat{f } (\xi ) | d\xi + \sup _z \left( |z |^2 e^{- \frac{1 }{2 } |z |^2 } \right) \cdot \int _{ 1 \leq |\xi | \leq \pi / h } \, |\xi |^{\kappa - 2 } |\widehat{f } (\xi ) | \, d\xi \right) \nonumber \\   
&\leq & C \max (t , 2 e^{-1 } ) \, h^{\kappa - 2 } || \widehat{f } ||_{\kappa - 2 , \kappa } , \nonumber   
\end{eqnarray}   
assuming wlog that $h \leq \pi . $ 
Since the integral over $|\xi | \geq h^{-1 } $ can be bounded by       
\begin{eqnarray*}   
\int _{|\xi | \geq h^{-1 } \pi } \, \left | e^{- h^{-2 } t G(h \xi ) } - e^{- t |\xi |^2 } \right | \, |\widehat{f }(\xi ) | \, d\xi   
&\leq & 2 \pi ^{- (\kappa - 2 ) } h^{\kappa - 2 } \int _{|\xi | \geq h^{-1 } \pi } \, |\xi |^{\kappa - 2 } \, |\widehat{f } (\xi ) | \, d\xi \\   
&\leq & 2 \pi ^{- (\kappa - 2 ) } h^{\kappa - 2 } || \widehat{f } ||_{\kappa - 2 , \kappa } ,   
\end{eqnarray*}   
the theorem follows, where in fact we have established the more precise bound        
$$   
\left( C_{1, \varphi } \max (t, 2e^{-1 } ) + C_{2, \varphi } + 2 \pi ^{- (\kappa - 2 ) } \right) \, || \widehat{f } ||^{\circ }_{\kappa - 2 , \kappa } \, h^{\kappa - 2 } ,   
$$   
with $C_{1, \varphi } := \sup _{0 < |\eta | \leq \pi } (G(\eta ) - |\eta |^2 ) / |\eta |^{\kappa } $ and $C_{2, \varphi } := \sup _{0 < |\eta | } (1 - \widehat{L }_1 (\eta ) ) / |\eta |^{\kappa } $, the constant of corollary \ref{cor:convergence_RBF}.   
\hfill $\Box $   
   
\begin{remarks} \rm{(i) If we strengthen the hypothesis on $\widehat{f } $ to $\widehat{f } |_{\mathbb{R }^n \setminus 0 } \in \mathring{L }^1 _{\kappa - 2 } (\mathbb{R }^n ) $, we obtain an error bound of $C h^{\kappa - 2 } || \widehat{f } || ^{\circ } _{\kappa - 2 } $ with a constant $C $ which is independent of $t . $   
\medskip   
   
\noindent (ii) The estimate for the integral over $|\xi | \geq \pi / h $ is quite rough, but note that since $h^{-2 } G(h\xi ) $ is $2 \pi / h $-periodic and equal to 0 in points of $2 \pi h^{-1 } \mathbb{Z }^n $, $e^{ - h^{-2 } G (h\xi ) } - e^{ - t |\xi |^2 } $ can get arbitrarily close to $1 $ on this set. A similar remark applies to the first term of (\ref{eq.proof_thm_conv_RBF_scheme_1}). There may be room for improvement, by further analyzing the contribution of a small neighborhood of this set of points.   
}   
\end{remarks}   
   
\medskip   
   
It is not difficult to verify that $h^{\kappa - 2 } $ is the exact order of approximation if $\kappa > 2 $, and that the scheme does not converge if $\kappa = 2 . $ Let      
\begin{equation} \label{eq:underline_g}   
\underline{g }_{\kappa } = \underline{g }_{\kappa } (\varphi ) := \liminf _{\eta \to 0 } \frac{ G_{\varphi } (\eta ) - |\eta |^2 }{|\eta |^{\kappa } } .   
\end{equation}   
We will see in proposition \ref{prop:formula_a_kappa} below that $\underline{g }_{\kappa } > 0 . $ 
   
\begin{theorem} \label{thm:lower_bound_RBF_scheme} 
Let $f \in L^{\infty } _{-p } (\mathbb{R }^n ) $ for some $p < \kappa $ such that $\widehat{f } |_{\mathbb{R }^n \setminus 0 } \in L^1 _{k, \kappa } (\mathbb{R }^n ) $ 
for some $k \in ( \kappa - 2 , \kappa ] . $ 
Then if $\kappa > 2 $,   
\begin{equation} \label{eq:lower_bound_RBF_scheme_kappa >2}    
\liminf _{h \to 0 } h^{- \kappa + 2 } || u_h (\cdot , t ) - u (\cdot , t ) ||_A \geq \underline{g }_{\kappa } \, t \int _{\mathbb{R }^n } |\xi |^{\kappa } |\widehat{f } (\xi ) | e^{- t |\xi |^2 } \, d\xi ,    
\end{equation}   
while if $\kappa = 2 $,   
\begin{equation} \label{eq:lower_bound_RBF_scheme_kappa=2}   
\liminf _{h \to 0 } || u_h - u ||_A \geq \int _{\mathbb{R }^n } \left( 1 - e^{ - \underline{g }_2 t |\xi |^2 } \right) e^{- t |\xi |^2 } |\widehat{f }(\xi ) | \, d\xi .   
\end{equation}   
\end{theorem}   
   
\noindent {\it Proof}. Since $|| \widehat{s_h [f] } - \widehat{f } ||_1 = O (h^k ) $ and $k > \kappa - 2 $, lemma \ref{lemma:u_h_hat_bis}, 
corollary \ref{cor:convergence_RBF} and Fatou's lemma imply that   
\begin{eqnarray*}   
&&\liminf _{h \to 0 } h^{- \kappa + 2 } || \widehat{u }_h (\cdot , t ) - \widehat{u }(\cdot , t ) ||_1 \\   
&\geq & \int _{\mathbb{R }^n } \, \liminf _{h \to 0 } h^{- \kappa + 2 } \left ( 
1 - e^{- h^{-2 } t R(h \xi ) } \right ) e^{- t |\xi |^2 } 
|\widehat{f } (\xi ) | \, d\xi   
\end{eqnarray*}   
where we have put $R ( \eta )= G(\eta ) - |\eta |^2 . $ By 
the mean value theorem applied to the exponential, there exist $\zeta = \zeta ( h , \xi , t ) \in [0, 1 ] $ such that 
$$   
\left( 1 - e^{- h^{-2 } t R(h \xi ) } \right) 
= h^{-2 } t R(h \xi ) e^{\zeta h^{ - 2 } R (h \xi ) } .   
$$   
If $\kappa > 2 $, then  
$R( h \xi ) \to 0 $ as $h \to 0 $ by proposition \ref{prop_G}(ii), and (\ref{eq:lower_bound_RBF_scheme_kappa >2}) follows from 
$$   
\liminf _{h \to 0 } h^{- \kappa } | R(h \xi ) | = \underline{g }_{\kappa } |\xi |^{\kappa } .   
$$   
If $\kappa = 2 $, then   
$$   
\liminf _{h \to 0 } \left( 1 - e^{ - h^{-2 } t R (h \xi ) } \right) = 1 - e^{ - \underline{g }_2 t  |\xi |^2 } ,   
$$   
which proves (\ref{eq:lower_bound_RBF_scheme_kappa=2}).   
\hfill $\Box $   
   
   
   
\subsection{Approximate approximation properties of the scheme} As we have just seen, our RBF scheme for the heat equation does not converge if $\kappa = 2 $, which is for example the case if $n = 1 $ and 
the basis function is the Hardy multiquadric. 
It turns out that we can then still achieve an arbitrarily small absolute error by an appropriate choice of the basis function, 
e.g. by introducing a shape parameter. 
This is again an approximate approximation phenomenon of the type encountered in section 4, and which if $\kappa > 2 $ will again take the form of an 
apparent rate of convergence better than than $O( h^{\kappa - 2 } ) $ up to some threshold $h_0 $, 
for initial conditions $f $ whose Fourier transform decay sufficiently rapidly at infinity. Suppose that $\varphi \in \mathfrak{B }_{2 , N } (\mathbb{R }^n ) $ with $N > n + 2 $ and $\kappa \geq 2 $, and 
let the initial condition $f $ 
be as in lemma \ref{lemma:u_h_hat_bis}. Let   
\begin{equation} \label{eq:def_alpha_kappa}   
\overline{g }_{\kappa } := \overline{g }_{\kappa } (\varphi ) := \limsup _{\eta \to 0 } \frac{G(\eta ) - |\eta |^2 }{|\eta |^{\kappa } } .   
\end{equation}   
   
\begin{theorem} \label{thm:approx_approx_kappa > 2} Let $\kappa \geq 2 $ and $s > \kappa - 2 . $ If 
$s > \kappa $ 
then there exist 
for any $\varepsilon > 0 $ a constant $C_{\varepsilon } $ which does not depend on $t \geq 0 $ 
such that if $\widehat{f } |_{\mathbb{R }^n \setminus 0 } \in \mathring{L }^1 _{\kappa } \cap \mathring{L }^1 _s $,   
\begin{eqnarray} \nonumber 
|| u_h ( \cdot , t ) - u (\cdot , t ) ||_A &\leq & (1 + \varepsilon ) \overline{g }_{\kappa } \, h^{\kappa - 2 } \int _{\mathbb{R }^n } \, t |\xi |^{\kappa } \, |\widehat{f } (\xi ) | e^{- t |\xi |^2 } \, d\xi \\   
&\ & 
(1+ \varepsilon ) \overline{l }_{\kappa } \, h^{\kappa } \, || \widehat{f } ||^{\circ } _{1 , \kappa } + C_{\varepsilon } h^s 
|| \widehat{f } ||^{\circ } _{1, \kappa } , \nonumber   
\end{eqnarray}    
while if $\kappa - 2 < s \leq \kappa $ and $\widehat{f } \in \mathring{L }^1 _{s , \kappa } $, then   
\begin{equation} \nonumber   
|| u_h ( \cdot , t ) - u (\cdot , t ) ||_A \leq (1 + \varepsilon ) \overline{g }_{\kappa } \, h^{\kappa - 2 } \int _{\mathbb{R }^n } \, t |\xi |^{\kappa } \, |\widehat{f } (\xi ) | e^{- t |\xi |^2 } \, d\xi +   
C_{\varepsilon } h^s || \widehat{f } ||^{\circ } _{1, s, \kappa }   
\end{equation}   
for sufficiently small $h . $ If $\kappa = 2 $, we can replace the first term on the right in these inequalities by   
$$   
\int _{\mathbb{R }^n } \, |\widehat{f }(\xi ) | \left( 1 - e^{- (1 + \varepsilon ) \underline{g }_2 t |\xi |^2 } \right) e^{- t |\xi |^2 } d\xi .   
$$   
\end{theorem}   
   
\noindent {\it Proof}. We adapt the proof of theorem \ref{thm:conv_RBF_scheme_1}.  If $\varepsilon > 0 $, there exists a 
$r(\varepsilon ) > 0 $ such that $ | G(\eta ) - | \eta |^2 | >  (1 + \varepsilon ) \overline{g }_{\kappa } | \eta |^{\kappa } $ if $ | \eta | < 
r(\varepsilon ) . $ Hence, assuming wlog that $r(\varepsilon ) \leq \pi $,      
\begin{eqnarray*}   
&&\int _{h|\xi | \leq r (\varepsilon ) } \, 
\left( 1 - e^{- t h^{-2 } ( G(h \xi ) - |h \xi |^2 ) } \right) e^{ - t |\xi |^2 } \, |\widehat{f }(\xi ) | \, d\xi \\   
&\leq & \int _{h|\xi | \leq r(\varepsilon ) } \left( 1 - e^{ - (1 + \varepsilon ) \overline{g }_{\kappa } h^{\kappa - 2 } | \xi |^{\kappa }  } \right) e^{ - t |\xi |^2 } \, |\widehat{f }(\xi ) | \, d\xi \\   
&\leq & (1 + \varepsilon ) \overline{g }_{\kappa } h^{\kappa - 2 } \int _{\mathbb{R }^n } \, t |\xi |^{\kappa } \, |\widehat{f } (\xi ) | e^{- t |\xi |^2 } \, d\xi .   
\end{eqnarray*}   
%
%
We estimate the integral over $|\xi | \geq r (\varepsilon) / h $ 
by $2r (\varepsilon ) ^{-s } h^{-s } \int _{|\xi | \geq r(\varepsilon ) / h } \, |\widehat{f } (\xi ) | \, |\xi |^s d\xi $, which we bound by $2 r(\varepsilon )^{-s } h^s || \widehat{f } ||^{\circ } _{1, s } $ if $s > \kappa $, and by $|| \widehat{f } ||^{\circ } _{1, s , \kappa } $ if $s \leq \kappa $, assuming in the latter case that $h \leq r (\varepsilon) . $ If we finally use theorem \ref{thm:approx_approx_RBF} to estimate first term of (\ref{eq.proof_thm_conv_RBF_scheme_1}) when $s > \kappa $, and corollary \ref{cor:convergence_RBF} when $s \leq \kappa $, 
the theorem follows. \hfill $\Box $   
   
\medskip   
   
   
\begin{corollary} Let $\kappa \geq 2 $ and suppose that 
$g_{\kappa } := \lim _{\eta \to 0 } |G(\eta ) - |\eta |^2 | / |\eta |^{\kappa } $ exists, so that   
$\underline{g }_{\kappa } = \overline{g }_{\kappa } =: g_{\kappa } . $ Suppose that $\widehat{f } |_{\mathbb{R }^n \setminus 0 } \in L^1 _{s , \kappa } (\mathbb{R }^n ) $ 
for some 
$s > \kappa . $ Then   
\begin{equation} \nonumber   
\lim _{h \to 0 } h^{-(\kappa - 2 ) } || u_h (\cdot , t ) - u(\cdot , t ) ||_A = \left \{   
\begin{array}{lll}   
g_{\kappa } t \int _{\mathbb{R }^n } \, |\widehat{f }(\xi ) | \, |\xi |^{\kappa } e^{- t |\xi |^2 } d\xi , &\kappa > 2, \\   
\ \\   
\int _{\mathbb{R }^n } \, |\widehat{f }(\xi ) | \left( 1 - e^{- g_2 t |\xi |^2 } \right) e^{- t |\xi |^2 } d\xi , &\kappa = 2 .   
\end{array}   
\right.   
\end{equation}   
\end{corollary}   
   
Compare with corollary \ref{corollary:approx_RBF}.  It follows from proposition \ref{prop:formula_a_kappa} below that $g_{\kappa } $ exists iff $\lim _{|\eta | \to 0 } |\eta |^{\kappa } \widehat{\varphi }(\eta ) $ exists. One can also give a direct proof of this corollary using Lebesgue's dominated convergence theorem: see the proof of theorem \ref{thm:cov_RBF_gen_symbols_bis} below.      
\medskip   
   
\begin{versionA}      
\noindent \textcolor{blue}{{\it Alternative, direct proof of the previous proposition using Lebesgue's dominated convergence theorem } (to be deleted from final version). By corollary \ref{cor:convergence_RBF},  $\lim _{h \to 0 } h^{\kappa - 2 } || \widehat {s }_h [f ] - \widehat{f } ||_1 = 0 $, so that the limit equals that of   
$$   
h^{\kappa - 2 } || \left( e^{- t h^{-2 } (G (h\xi ) } - e^{- t |\xi |^2 } \right) \widehat{f } ||_1 ,   
$$   
as $h \to 0 . $ As before, we split this integral into an integral over $h|\xi | \leq r $ and one over $h |\xi | > r $,  and observe that the contribution of the second integral to the limit vanishes, given the hypothesis on $\widehat{f  } . $ For the first integral we can use Lebesgue's dominated convergence theorem, if we note that 
$$   
h^{ - (\kappa - 2 ) } \left( e^{- t h^{-2 } (G (h\xi ) - |h \xi |^2 ) } - 1 \right) e^{- t |\xi |^2 } \mathbf{1 }_{ \{ h |\xi | \leq 1 \} } \to g_{\kappa } t |\xi |^{\kappa } e^{- t |\xi |^2 },   
$$   
if $\kappa > 2 $, while it can be uniformly bounded by $C |\xi |^{\kappa } e^{- t |\xi |^2 / 2 } . $ If $\kappa =2 $, then   
$$   
\left( e^{- t h^{-2 } (G (h\xi ) - |h \xi |^2 ) } - 1 \right) e^{- t |\xi |^2 } \mathbf{1 }_{ \{ h |\xi | \leq 1 \} } \to \left( e^{- g_2 t |\xi |^2 } - 1 \right) e^{- t |\xi |^2 } .      
$$   
\hfill $\Box $   
}   
\medskip   
   
\end{versionA}   
   
If $\overline{g }_{\kappa } $ and $\overline{l }_{\kappa } $ 
are small, then theorem \ref{thm:approx_approx_kappa > 2} with $\varepsilon $ for example equal to $\max (\overline{g }_{\kappa } , \overline{l }_{\kappa } ) $ shows one may observe a higher apparent rate of convergence than the actual rate for small but not too small $h $'s if $\widehat{f }(\xi ) $ decays sufficiently rapidly at infinity. As in section 4, we can construct basis functions $\varphi $ with small $g_{\kappa } (\varphi ) $ by taking these of the form $\varphi (x) = \phi (c^{-1 } x ) $ and letting $c \to \infty . $ We start by deriving explicit formulas for $\overline{g }_{\kappa } (\varphi ) $ and $\underline{g }_{\kappa } (\varphi ) . $ 
Recall the definition of $\underline{A } (\varphi ) $ and $\overline{A } (\varphi ) $ in proposition \ref{prop:formula_l_kappa}.   
\medskip   
   
   
   
\begin{proposition} \label{prop:formula_a_kappa}   We have   
\begin{equation} \label{eq: estimate_a_kappa}   
\overline{g }_{\kappa } (\varphi ) = \frac{1 }{\underline{A }(\varphi ) } \sum _{k \neq 0 } | 2 \pi k |^2 \widehat{\varphi } (2 \pi k ) , \ \ \underline{g }_{\kappa } (\varphi ) = \frac{1 }{\overline{A }(\varphi ) } \sum _{k \neq 0 } | 2 \pi k |^2 \widehat{\varphi } (2 \pi k ) 
\end{equation}   
\end{proposition}   
   
\noindent {\it Proof}.   
   
\begin{eqnarray*}   
G(\eta ) - |\eta|^2 &=& \frac{\sum _k |\eta + 2 \pi k |^2 \widehat{\varphi } (\eta + 2 \pi k ) }{\sum _k \widehat{\varphi } (\eta + 2 \pi k ) } - |\eta |^2 \\   
&=& \frac{ \sum _{k \neq 0 } \left( 4 \pi (\eta , k ) + 4 \pi ^2 |k |^2 \right) \widehat{\varphi } (\eta + 2 \pi k ) }{\sum \widehat{\varphi } (\eta + 2 \pi k ) } \\   
&=& \frac{g(\eta ) }{\widehat{\varphi } (\eta ) + h(\eta ) } ,   
\end{eqnarray*}   
where $g(\eta ) := \sum _{k \neq 0 } \left( 4 \pi (\eta , k ) + 4 \pi ^2 |k |^2 \right) \widehat{\varphi } (\eta + 2 \pi k ) $ and $h(\eta ) := \sum _{k \neq 0 } \widehat{\varphi } (\eta + 2 \pi k ) $ are continuous 
in a neighborhood of 0 . It follows that   
$$   
\limsup _{\eta \to 0 } \left | \frac{G(\eta ) - |\eta |^2 }{|\eta |^{\kappa } }  \right | = \limsup _{\eta \to 0 } \frac{ | g(\eta ) | }{|\eta |^{\kappa } \widehat{\varphi } + |\eta |^{\kappa } h(\eta ) } = \frac{g(0) }{\underline{A } } = \frac{\sum _k 4 \pi ^2 |k |^2 \widehat{\varphi } (2 \pi k ) }{\underline{A } } ,   
$$   
with $\underline{A } = \underline{A } (\varphi ) . $ The formula for $\underline{g }_{\kappa } (\varphi ) $ follows similarly. \hfill $\Box $   
\medskip   
   
Note that $\overline{l }_{\kappa } (\varphi ) \leq \overline{a }_{\kappa } (\varphi ) . $ 
Also note that $\underline{g }_{\kappa } > 0 $ since $\overline{A } < \infty $ and $\overline{g }_{\kappa } < \infty $ since $\underline{A } > 0 $, by the ellipticity condition on $\widehat{\varphi } $ at 0. 
\medskip   
   
If we take $\varphi (x) := \phi _c (x) = \phi (x/c ) $, with $\phi \in \mathfrak{B }_{\kappa , N } (\mathbb{R }^n ) $, then $\widehat{\varphi }(\eta ) = c^n \widehat{\phi }(c \eta ) $, and $\underline{A }(\phi _c ) = c^{n - \kappa } \underline{A }(\phi ) . $ It follows that   
\begin{eqnarray*}     
\overline{g }_{\kappa } (\phi _c ) &=& c^{\kappa } \underline{A }(\phi ) \sum _{k \neq 0 } |k |^2 \widehat{\phi } (2 \pi c k ) \\   
&\leq & C c^{\kappa - N } \sum _{k \neq 0 } |k |^{2 - N } ,   
\end{eqnarray*}   
where the series converges since $N > n + 2 . $ 
   
\begin{corollary} \label{corr:formulas_l_kappa-a_kappa} If $\phi \in \mathfrak{B }_{\kappa , N } (\mathbb{R }^n ) $ with $N > \max (n + 2 , \kappa ) $, then $\overline{g }_{\kappa } (\phi _c ) \to 0 $ as $c \to \infty . $   
\end{corollary}   
   
\begin{examples} \label{examples_multiquadric} \rm{(i) Hardy's multiquadric with shape parameter $c $ is defined by   
\begin{equation}   
\varphi (x) := - \sqrt{|x |^2 + c^2 } , \ \ x \in \mathbb{R }^n .   
\end{equation}   
where the minus sign serves to make $\widehat{\varphi } (\eta ) $ positive. Note that $   
\varphi (x) = c \phi (x/c ) $ with $\phi (x) := - \sqrt{|x |^2 + 1 } $, so that we are in the situation of corollary \ref{corr:formulas_l_kappa-a_kappa}, except for an irrelevant multiplicative factor of $c . $ The Fourier transform of $\varphi $ on $\mathbb{R }^n \setminus 0 $ is given by   
$$   
\widehat{\varphi }(\eta ) = \pi ^{-1 } \left(2 \pi c \right)^{(n + 1 ) / 2 } |\eta | ^{- (n + 1 ) / 2 } K_{(n + 1 ) / 2 } (c |\eta | ) ,   
$$   
where $K_{\nu } $ is the MacDonald function, or modified Bessel function of the 
second kind: see for example cf. Baxter \cite{Ba1}. The limiting form of $K_{\nu }  $ for small values of the argument implies that  as $\eta \to 0 $, $\widehat{\varphi }(\eta ) \simeq  A_n |\eta |^{- n - 1 } $ (with $A_n = 2^{n } \pi ^{(n - 1 )/2 } \Gamma \left( \frac{n + 1 }{2 } \right) $), so that $\kappa = n + 1 $, and our RBF-scheme for the heat equation will converge if $n \geq 2 $, at a rate of $h^{n - 1 } . $ The MacDonald function is known to decay exponentially at infinity, so that we can apply corollary \ref{corr:formulas_l_kappa-a_kappa} to conclude that $\overline{l }_{\kappa } (\phi _c ) $ and $\overline{g }_{\kappa } (\phi _c ) \to 0 $ as $c \to \infty . $ In fact, these will converge to 0 at an exponential rate since $\sum _k |k |^2 \widehat{\phi } (2 \pi c k ) $ does.   
   
If $n =1 $, then $\kappa = 2 $, 
and the scheme will not converge. However, corollary \ref{corr:formulas_l_kappa-a_kappa} together with 
theorem \ref{thm:approx_approx_kappa > 2} shows that we can make the error arbitrarily small by taking the shape parameter $c$ sufficiently large, with moreover an arbitrarily large apparent order of convergence for small but not-too-small $h $'s if the Fourier transform of the initial value decays sufficiently rapidly at infinity. At first sight, this may seem strange, because we are after all simply performing an additional scaling by $c $, and we are already using scaled basis functions $\varphi _h (x) = \varphi (h^{-1 } x ) $ for our interpolation. Note, 
though, that we are interpolating with $\phi _{ch } $ on $h \mathbb{Z }^n $, and not on $ch \mathbb{Z }^n . $   
\medskip   
   
\noindent (ii) If we take a homogeneous basis function, $\phi (x) = |x |^p $ with $p > 0 $, then $\widehat{\phi } (\eta ) $ is proportional to $|\eta |^{- p - n } $ on $\mathbb{R }^n \setminus 0 $, so that $\kappa = n + p = N $ and corollary \ref{corr:formulas_l_kappa-a_kappa} does not apply, as indeed it shouldn't: if $\phi $ is homogeneous, then  $\widehat{L }_1 (\phi _c ) $ and $G(\phi _c ) $ are independent of $c $, and therefore $\overline{g }_{\kappa } (\phi _c ) $ and $\overline{l }_{\kappa } (\phi _c ) $ also.     
}   
\end{examples}   
   
\begin{remark} \rm{One can perform a similar analysis of the $c $-dependence of the constant $C_{\varepsilon } $ of theorem \ref{thm:approx_approx_kappa > 2} 
as the one we did in subsection \ref{subsection:c-dependence}, with similar conclusions; we skip the details.   
}   
\end{remark}   
   
\section{\bf Convergence of stationary RBF schemes for pseudo-differential evolution equations}   
   
The results of the previous section remain valid for a large class of constant coefficient pseudo-differential evolution equations 
\begin{equation} \label{eq:PS_1}   
\partial _t u + a(D ) u = 0 , \ \ t > 0 , \\   
\end{equation}   
under suitable conditions on the symbol $a = a(\xi ) $, notably ${\rm Re } \, a(\xi ) \geq 0 . $ The operatore $a(D ) $ is defined by $\widehat{a(D ) f  } = a \widehat{f } $ 
initially with domain $\mathcal{S }(\mathbb{R }^n ) $, for example. We are in fact restricting ourselves to a rather special class of pseudo-differential operators, the Fourier multiplier operators 
or convolution operators: if $a $ is a tempered distribution and if $f $ is a Schwarz class function, then $a(D ) f $ is the convolution of $f $ with the inverse Fourier transform of $a . $ These can also be considered as constant coefficient pseudo-differential operators, general pseudo-differential operators having symbols which also depend on $x . $ The latter are outside of the scope of this paper, but the  multiplier operators we consider here already contain many interesting examples, such as the fractional Laplacians or the generators of 
L\'evy processes. Regarding the latter, the equation (\ref{eq:PS_1}) occurs in mathematical finance, for derivative pricing in exponential L\'evy models,  and has been treated numerically in \cite{RB_RC} using the RBF scheme we investigate here, with good results. We note that, from a theoretical point of view,  convergence of this type of stationary scheme for convolution operators is 
not obvious, since these, as integral operators, are in general non-local (except if $a(\xi ) $ is a polynomial), and we already know from 
section 3 that to obtain good convergence we will need to use 
basis functions with polynomial growth. 
To understand 
why such schemes nevertheless perform well numerically,  
as for example observed in \cite{RB_RC}, was a main motivation for this paper.   
   
   
As regards the conditions on the symbols, we will work with the class $S^q _0 :=  S^q _0 (\mathbb{R }^n ) $ 
of $C^{\infty } $-functions $a : \mathbb{R }^n \to \mathbb{C } $ for which there exists for each multi-index $\alpha $ a constant $C_{\alpha } $ such that   
\begin{equation} \label{cond_symbol_1}   
|\partial _{\xi } ^{\alpha } a(\xi ) | \leq C_{\alpha } (1 + |\xi | )^q , \ \ \xi \in \mathbb{R }^n .   
\end{equation}   
We will not need the faster $(1 + |\xi | )^{q - |\alpha | } $-decay for 
$\partial _{\xi }  ^{\alpha } a $ which is a standard requirement in much of pseudo-differential theory and which for example is satisfied by the symbols of partial differential operators. The requirement of having $C^{\infty } $-symbols is slightly restrictive, and a priori excludes symbols such as $|\xi |^q $ for non-integer but positive $q $, but the theory presented here will apply to regularized versions of such symbols, for example, replacing $|\xi |^q $ by $(1 - \chi ) |\xi |^q $ with $\chi \in C^{\infty }_c (\mathbb{R }^n $ equal to 1 on a neigborhood of 0. Modifying a symbol on a compact set will not affect the singularities of the (distributional) kernel of $a(D ) $ but will change its decay properties at infinity. Care then has to be taken with the growth properties of the initial values $f $ which we allow for (\ref{eq:PS_1}) (equivalently, the singularities at 0 of $\widehat{f } $), in the various convergence theorems, which makes the statements of these theorems more complicated. On the other hand, if $a \in S^q _0 $ has non-negative real part, then the solution of (\ref{eq:PS_1}) with initial value $f $ makes sense for any tempered distibution $f $, being the inverse Fourier transform of $e^{- t a } \widehat{f } . $  We therefore limit ourselves here to smooth symbols.      
   
We first examine the action of $a(D) $ on $L_1 $:   
\medskip   
   

   
   
\begin{theorem} 
\label{lemma:a(D)L_1} Let $\kappa \geq 0 . $ If $a \in S^q _0 (\mathbb{R }^n ) $ and if $\varphi \in\mathfrak{B }_{\kappa , N } (\mathbb{R }^n ) $ with $N > n + q $, then $a(D ) L_1 $ is a bounded continuous function and there exists a constant $C > 0 $ such that   
$| a(D ) L_1 (x) | \leq C (1 + |x | )^{- (n + \kappa ) } . $   
\end{theorem}   
   
The proof is similar to that of the bound (\ref{bound_L_1}) of theorem \ref{thm_L_1}: see appendix \ref{Appendix:Lagrange_function}. In fact, here, and below, it would have sufficed to require (\ref{cond_symbol_1}) only for $|\alpha | \leq \lceil \kappa \rceil + n + 1 .  $ We suppose from now on that $N > n + q . $     
\medskip   
   
The second condition we will need to put on the symbol is that it has a non-negative real part:   
\begin{equation}\label{cond_symbol_2}   
{\rm Re } \, a(\xi ) \, \geq \, 0 .  
\end{equation}    
Perhaps curiously, we do not need $a(\xi ) $ or ${\rm Re } \, a(\xi ) $ to be elliptic. In particular, our results below will for example also apply to the free Schr\"odinger operator, for which 
$a(\xi ) = i |\xi |^2 $, or the regularized "half-wave equation", with $a(\xi ) = |\xi | $ outside of a neighborhood of 0.  
The heat equation obviously also falls within the class of allowed evolution equations, as do the Kolmogorov-Fokker-Planck equations associated to certain L\'evy processes: see example \ref{ex:pseudo-diff_evolution_eq}(ii) below. 

The proofs in this section will be similar to the ones for the classical heat equation in section 5, and we will only signal the differences.   
\medskip   
   
We will extend the results of section 5 to the Cauchy problem for (\ref{eq:PS_1}) for symbols satisfying (\ref{cond_symbol_1}), (\ref{cond_symbol_2}). The proofs will be similar to the ones for the classical heat equation in section 5, and we will mainly signal major differences. We are again interested in solving (\ref{eq:PS_1}) with initial value $f $ using a semi-discrete scheme which is the RBF-variant of the classical method of lines, looking for approximate solutions of the form (\ref{eq:def_u_h}), where $L_1 $ is the Lagrange function on $\mathbb{Z }^n $ associated to a basis function 
$\varphi \in \mathfrak{B }_{\kappa , N } (\mathbb{R }^n ) $ with $\kappa > 0 $ and $N > n + q $: cf. 
theorem \ref{lemma:a(D)L_1}. The coefficients $c_k (t; h ) $ of $u_h $ are again determined by requiring that $u_h $ solve (\ref{eq:PS_1}) exactly in the interpolation points: $\partial _t u_h (hj , t ) = - a(D) u_h (jh , t ) $ for $j \in \mathbb{Z }^n . $ This now leads to the (infinite) system of ODEs  
\begin{equation} \label{eq:ODE-system_a}   
\frac{d c_j (t; h ) }{dt } = - \sum _k a(h^{-1 } D_x ) ( L_1 ) (j - k ) c_k (t; h )   
\end{equation}   
where $a(h^{-1 } D_x ) $ has symbol $a(h^{-1 } \xi ) . $  
We again have to solve this system with initial condition $c_k (0) = f(hj ) . $ 
One shows as in lemma \ref{lemma:ODE-system} that if $p < \kappa $ then there exists a unique solution in $C^{\infty } ([0, \infty ) , \ell ^{\infty } _{ - p } ) $ and that, as a consequence, $u_h [f ] (\cdot , t ) $ is in $L^{\infty } _{-p } (\mathbb{R }^n ) $ if $f \in L^{\infty } _{- p } (\mathbb{R }^n ) $ with norm bounded by a constant times that of $f . $   
   
\begin{remark} \rm{A noteworthy feature of the RBF-scheme is that we do not need to discretize the operator $a(D ) $, contrary to for example Finite Difference schemes, but only need to know its action on $L_1 $ (or on $\varphi $ when working with irregularly spaced interpolation points). This is an advantage when the operator is a singular integral operator: 
see \cite{RB_RC} for concrete examples and further discussion.   
}   
\end{remark}   
   
To analyze the RBF scheme we introduce the auxiliary function $G_a ^* $ on $\mathbb{R }^n \times \mathbb{R }_{> 0 } $ defined by   
\begin{eqnarray} \label {eq:def_G_a}   
G_a ^* (\xi ; h ) 
&:= & \sum _k a(\xi + 2 \pi h^{-1 } k ) \widehat{L }_1 (h \xi + 2 \pi k ) \\  
&=& \frac{\sum _k a(\xi + 2 \pi h^{-1 } k ) \widehat{\varphi } (h \xi + 2 \pi k ) }{\sum _{\nu } \widehat{\varphi }(h \xi + 2 \pi \nu ) } , \nonumber   
\end{eqnarray}   
where the series converges absolutely since $N > n + q . $ To make the connection with the previous section note that if $a(\xi ) = |\xi |^2 $ then $G_a ^* (\xi ) = h^{-2 } G (h\xi ) $, with $G $ given by (\ref{eq:def_G}).   

One shows, 
similarly to lemma \ref{lemma:u_h_hat}, that if the initial condition $f $ is a Schwarz-class function, or more generally if $\widehat{f } \in L^1 (\mathbb{R }^n ) $, and if $a \in S^q _0 $ satisfies 
(\ref{cond_symbol_2}), then the Fourier transform with respect to $x $ of $u_h (x , t ) $ is given by   
\begin{equation} \label{u_h_hat_bis}   
\widehat{u }_h (\xi , t ) = e^{- t  G_a ^* (\xi ; h ) } \widehat{s_h [f] }(\xi ) .
\end{equation}   
Since $G_a ^* (\xi ; h ) $ is in $C_b ^{\lceil \kappa \rceil - 1 , \kappa - (\lceil \kappa \rceil - 1 ) } (\mathbb{R }^n ) $, we can extend this formula to initial values  $f $ of polynomial growth strictly less than $\kappa $ whose Fourier transform coincides on $\mathbb{R }^n \setminus 0 $ with an element of $L^1 (\mathbb{R }^n , (|\xi |^{\kappa }\wedge 1 ) d\xi ) . $   
%
We also note that $G_a ^* (\xi ; h ) $ is $2\pi / h $-periodic in $\xi $ and has non-negative real part. Its zero-set contains $2\pi h^{-1 } \mathbb{Z } \setminus 0 $ but may be bigger. 
We have the following basic estimate which generalizes proposition \ref{prop_G}(ii).   
   
\begin{proposition} \label{prop_G_a} Suppose that $a \in S^q _0 (\mathbb{R }^n ) $ for some 
$q \in \mathbb{R } $, and let $\varphi \in \mathfrak{B }_{\kappa , N } (\mathbb{R }^n ) $ with $N > n + q . $ Then there exists a constant $C $ such that 
for all $h < 1 $, 
\begin{equation}   
|G_a ^* (\xi ; h ) - a(\xi ) | \leq C h^{\kappa - \max (q , 0 ) } |\xi |^{\kappa } , \ \ |\xi | \leq \pi / h .   
\end{equation}   
\end{proposition}   
   
\noindent {\it Proof}. We have   
$$   
G_a ^* (\xi ; h ) - a(\xi ) = a(\xi ) \left( \widehat{L }_1 (h\xi ) - 1 \right) + \sum _{k \neq 0 } a \left( \xi + \frac{2 \pi }{h } \right) \widehat{L }_1 (h \xi + 2 \pi k ) .   
$$   
The first term is bounded by a constant times $(1 + |\xi | )^q |h \xi |^{\kappa } \leq C h^{\kappa - \max (q , 0 ) } |\xi |^{\kappa } $ if $|h \xi | \leq \pi . $ As for the other terms, $|\xi + 2 \pi k / h | $ is comparable to $|k |/h $ if $|\xi | \leq \pi / h $, so $|a(\xi + 2 \pi k / h ) | \leq C h^{-q } |k |^q $, and by (\ref{eq:SFC1}),   
\begin{eqnarray*}   
\widehat{L }_1 (h \xi + 2 \pi k ) 
\leq C |h \xi |^{\kappa } |k |^{-N } , \ \ |h \xi | \leq \pi ,   
\end{eqnarray*}   
so that   
$$   
\sum _{k \neq 0 } | a \left( \xi + \frac{2 \pi }{h } \right) \widehat{L }_1 (h \xi + 2 \pi k ) | \leq C h^{\kappa - q } |\xi |^{\kappa } \sum _{k \neq 0 } |k |^{q - N } 
\leq C h^{\kappa - q } |\xi |^{\kappa } .   
$$   
\hfill $\Box $   
\medskip   
   
   
We will suppose from now on that $f \in L^{\infty } _{ -p } (\mathbb{R }^n ) $ 
for some $p < \kappa $ such that $\widehat{f } |_{\mathbb{R }^n \setminus 0 } \in L^1 (\mathbb{R }^n , (|\xi |^{\kappa } \wedge 1 ) d\xi ) . $ The unique solution of the initial value problem in the space of tempered distributions is given by $\widehat{u } (\xi , t ) = e^{- t a(\xi ) } \widehat{f } $ which is well-defined as a product of a tempered distribution and a $C^{\infty } $-function all of whose derivatives  have at most polynomial growth. 
One shows, using the arguments of Appendix B, that $u_h (x, t ) $ is a continuous function of polynomial growth of order at most $\kappa $  for each $t > 0 $, 
and that the Fourier transform of $u_h (x , t ) - u(x, t ) $ is given by     
$$   
e^{- t h^{-2 } G_a ^* (\xi ; h ) } \left( \widehat{s_h [f ] }(\xi ) - \widehat{f }(\xi ) \right) + \left(e^{- t ( G_a ^* (\xi ; h ) - a(\xi ) ) } - 1 \right) e^{ - t a(\xi ) } \widehat{f }(\xi ) .   
$$     
We then can state the following convergence theorem. 
 
\begin{theorem} \label{thm:cov_RBF_gen_symbols} 
Let $a \in S^q _0 (\mathbb{R }^n ) $, 
${\rm Re } \, a(\xi ) \geq 0 $, and suppose that 
$\kappa \geq q . $ Then there exists a constant $C $ such that if $f $ is a function of polynomial growth of order strictly less than $\kappa $ such that $\widehat{f }|_{\mathbb{R }^n \setminus 0 } \in \mathring{L }^1 _{\kappa } (\mathbb{R }^n ) $,  
then   
\begin{equation}   
|| u_h (\cdot , t ) - u (\cdot , t ) ||_A \leq C \, 
\max(t, 1 ) \cdot h^{\kappa - \max (q , 0 ) } || \widehat{f } ||^{\circ } _{\kappa } .   
\end{equation}   
\end{theorem}   
   
\noindent {\it Proof}. The proof is similar to the proof of theorem \ref{thm:conv_RBF_scheme_1}, with a small twist. By proposition \ref{prop_G_a} and the elementary inequality $|e^z - 1 | \leq |z | e^{\max ( {\rm Re } \, z , 0 ) } $ for $z \in \mathbb{C } $, 
and since\footnote{In the case of the heat equation ($a(\xi ) = |\xi |^2 $), $a - G_a ^* $ was negative for small enough $h |\xi | $, but this is not clear for ${\rm Re } (a - G_a ^* ) $ for general $a $, in particular when ${\rm Re } \, a $ is not elliptic. 
}      
$$   
{\rm Re } (a(\xi ) - G_a ^* (\xi ; h ) ) \leq {\rm Re } \, a(\xi ) (1 - \widehat{L }_1 (h \xi ) ) \leq C h^{\kappa } |\xi |^{\kappa }  {\rm Re } \, a(\xi ) \leq \frac{1 }{2 } {\rm Re } \, a(\xi ) ,   
$$   
for sufficiently small $h |\xi | $, there exists an $\rho > 0 $ such that if $ h |\xi | \leq \rho $, then 
\begin{eqnarray} \nonumber   
\left | e^{- t G_a ^* (\xi ; h ) } - e^{ - t  a(\xi ) ) } \right | &=& e^{ - t \, {\rm Re } \, a(\xi ) } \left | e^{- t (G(\xi ; h ) - a(\xi ) ) } - 1 \right | \\   
&\leq & C \, t \, h^{\kappa - \max (q, 0 ) } |\xi |^{\kappa } e^{ - \frac{1 }{2 } t \, {\rm Re } \, a(\xi ) } , \label{eq:cov_RBF_gen_symbols1}   
\end{eqnarray}   
which in absence of further hypotheses on the symbol $a(\xi ) $ we simply bound by $C h^{\kappa - p } |\xi |^{\kappa } . $ The rest of the proof proceeds as before: 
since ${\rm Re } \, G_a ^* \geq 0 $,    
$$   
|| u_h (\cdot , t ) - u (\cdot , t ) ||_A \leq \int \left | e^{- t G_a ^* } - e^{ - t  a } \right | \, |\widehat{f } | d\xi  + || \widehat{s_h [f ] } - \widehat{f } ||_1 . 
$$   
The last term can be estimated using theorem \ref{thm:convergence_RBF_bis} while the integral can be bounded by 
$$      
\int _{|\xi | \leq \rho / h } \left | e^{- t G_a ^* } - e^{ - t  a } \right | \, | \widehat{f } | d\xi + 2 \int _{h |\xi | \geq \rho } |\widehat{f } | d\xi   
$$   
Now 
use (\ref{eq:cov_RBF_gen_symbols1}) for the first integral,   
and 
bound the second one 
by $(h / \rho )^{\kappa } \int _{|\xi | \geq \rho / h }|\xi | ^{\kappa }  |\widehat{f } | d\xi \leq (h / \rho )^{\kappa } || \widehat{f } ||^{\circ } _{1 , \kappa } . $   

%
   
\hfill $\Box $   
\medskip   
   
We note that we do not get an improved convergence rate beyond $O(h^{\kappa } ) $ if $a $ is of negative order $q < 0 $, even if $\widehat{f } $ were rapidly decreasing and $ a $ were  smoothing. We therefore limit ourselves to operators of non-negative order: $q \geq 0 . $   
   
Also note that, on comparing theorem \ref{thm:cov_RBF_gen_symbols} for $a(\xi ) = |\xi |^2 $ with theorem \ref{thm:conv_RBF_scheme_1}, we require a stronger decay of $\widehat{f } $ at infinity, which translates into two additional degrees of smoothnes (two extra derivatives) of $f . $ If we assume that ${\rm Re } \, a(\xi ) $ is elliptic, then   
\begin{equation} \label{eq:hypo_a1}   
\sup _{\mathbb{R }^n } t |\xi |^q e^{- \frac{1 }{2 } t {\rm Re } \, a (\xi ) } < \infty ,   
\end{equation}   
is independent of $t $ and theorem \ref{thm:cov_RBF_gen_symbols} remains valid if $\widehat{f }|_{\mathbb{R }^n \setminus 0 } \in L^1 (\mathbb{R }^n , |\xi |^{\kappa - q } \wedge |\xi |^{\kappa } ) $, on replacing $|| \widehat{f } ||^{\circ } _{1, \kappa } $ by $|| \widehat{f } ||^{\circ } _{1, \kappa - q ,  \kappa } $: compare the proof of theorem \ref{thm:conv_RBF_scheme_1}.   
\medskip   
   
\medskip   
   
\begin{versionA}   
\noindent \textcolor{blue}{What if ${\rm Re } \, a (\xi )$ is hypo-elliptic in the sense that 
$$   
{\rm Re } \, a(\xi ) \geq c |\xi |^{q - \varepsilon }    
$$   
for some $\varepsilon \in [0 , q ) $, $c > 0 $ and for all sufficiently big $|\xi | ? $}   
\medskip   
   
\end{versionA}   
\medskip   
   
   
The argument can be refined to give a rough approximate approximation estimate: the proof above shows that if 
$r \leq \min (\rho , \pi ) $ and if $s > \kappa- q \geq 0 $,            
\begin{equation} \label{eq:rough_approx_approx_1}   
|| u_h (\cdot , t ) - u (\cdot , t ) ||_A \leq C(r) \, t \cdot h^{\kappa - q } || \widehat{f } ||^{\circ } _{\kappa } + \frac{h ^s }{r^s } || \widehat{f } ||^{\circ } _s + || \widehat{s_h [f ] } - \widehat{f  } ||_1 , 
\end{equation}   
where   
$$   
C(r) := \sup _{ h |\xi | \leq r } \frac{| a(\xi ) - G_a ^* (\xi ; h ) | }{ h^{\kappa - q } |\xi |^{\kappa } } .   
$$   
On closer examination, the proof of proposition \ref{prop_G_a} shows that   
$$   
C(r) \leq (4 \pi )^q 
|| a ||_{\infty , -q } \left( \sup _{|\eta | \leq r } \frac{1 - \widehat{L }_1 (\eta ) }{|\eta |^{\kappa } } + \sum _{k \neq 0 }  |k |^q \cdot \sup _{|\eta | \leq r } \frac{\widehat{L }_1 ( \eta + 2 \pi k ) }{|\eta |^{\kappa } } \right) ,   
$$   
where we recall that $|| a ||_{\infty , -q } = \sup _{\xi } (1 + |\xi |)^{-q } |a (\xi ) | $ (the constant in front is not optimal). From this it easily follows that 
\begin{equation} \label{eq:rough_approx_approx_2}   
\limsup _{r \to 0 } C(r) \leq (4 \pi )^q 
|| a || _{\infty , -q } \cdot \overline{l }_{\kappa , q } ,   
\end{equation}   
where we've put   
\begin{equation} \label{eq:rough_approx_approx_3}   
\overline{l }_{\kappa , q } 
:= \overline{l }_{\kappa , q } (\varphi ) := \underline{A }^{-1 }\sum _{k \neq 0 } (1 + |k |^q ) \widehat{\varphi } (2 \pi k ) ,   
\end{equation}   
with $\underline{A } $ defined by (\ref{eq:def_A(varphi)}); note that 
the series converges since $N > n + q $, by assumption.   
\medskip   
   
\begin{versionA}   
\noindent \textcolor{blue}{Some details: if $h |\xi | \leq r $, then   
\begin{eqnarray*}   
&&\frac{| a(\xi ) - G_a ^* (\xi ; h ) | }{ h^{\kappa - q } |\xi |^{\kappa } } \\   
&\leq & \frac{|a (\xi ) | (1 - \widehat{L }_1 (h \xi ) ) }{h^{\kappa - q } |\xi |^{\kappa } } + \sum _{k \neq 0 } \frac{ |a (\xi + 2 \pi h^{-1 } k ) | \widehat{L }_1 (h \xi + 2 \pi k ) }{h^{\kappa - q } |\xi |^{\kappa } } \\   
&\leq & 
|| a ||_{\infty , -q } \left( (1 + h^{-1 } r )^q \frac{(1 - \widehat{L }_1 (h \xi ) ) }{h^{\kappa - q } |\xi |^{\kappa } } + \sum _{k \neq 0 } (1 + (2 \pi |k |  + r ) h^{-1 } )^q \frac{ \widehat{L }_1 (h \xi + 2 \pi k ) }{h^{\kappa - q } |\xi |^{\kappa } } \right) \\   
&\leq &   
(4 \pi )^q 
|| a ||_{\infty , -q } \left( \frac{(1 - \widehat{L }_1 (h \xi ) }{h^{\kappa } |\xi |^{\kappa } } + \sum _{k \neq 0 }  |k |^q \frac{ \widehat{L }_1 (h \xi + 2 \pi k ) }{h^{\kappa }  |\xi |^{\kappa } } \right) ,   
\end{eqnarray*}   
where we used that if (wlog) $h < \pi $, then $1 < \pi / h $, and, since $r < \pi $,       
$$   
1 + \frac{r }{h } \leq \frac{\pi + r }{h } < \frac{2 \pi }{h } ,   
$$   
while   
$$   
1 + \frac{2 \pi |k | + r }{h } < \frac{2 \pi (|k |  + 1 ) }{h } \leq \frac{4 \pi |k | }{h } ,   
$$   
and the estimate follows.   
}   
\end{versionA}   
   
Note that 
$2 \overline{l }_{\kappa , 0 } (\varphi ) = \overline{l }_{\kappa } (\varphi ) = \overline{l }_{\kappa } . $ Equations (\ref{eq:rough_approx_approx_1}) , (\ref{eq:rough_approx_approx_3}) and theorems \ref{thm:approx_approx_RBF} (when $s > \kappa $) and \ref{thm:convergence_RBF_bis} (when $\kappa - q < s \leq \kappa $) now imply the following rough approximate approximation estimate:
   
\begin{theorem} \label{thm:rough_approx_approx} Suppose that $a \in S^q _0 (\mathbb{R }^n ) $, $0 \leq q \leq \kappa $, with ${\rm Re } \, a(\xi ) \geq 0 . $ Then there exists, for each $\varepsilon > 0 $ 
a constant $C_{\varepsilon } $ such that for $s > \kappa - q $ and $\widehat{f } |_{\mathbb{R }^n \setminus 0 } \in 
\mathring{L }^1 _{\kappa } (\mathbb{R }^n ) \cap \mathring{L }^1 _s (\mathbb{R }^n ) $, the error $|| u_h (\cdot , t ) - u (\cdot , t ) ||_A $ can be bounded by   
\begin{equation} \label{eq:rough_approx_approx_4}
(1 + \varepsilon )\left( C_a \, \overline{l }_{q, \kappa } t \cdot h^{\kappa - q } + \overline{l }_{\kappa } h^{\kappa } \right) || \widehat{f } ||^{\circ } _{\kappa } + C_{\varepsilon } h^s || \widehat{f } ||^{\circ } _s   
\end{equation}   
with $C_a = (4 \pi )^q  || a ||_{S^q _0 } $, where the second term can obviously be absorbed in the final one if $s \leq \kappa . $
\end{theorem}   
   
We call this a rough approximate approximation theorem since the bound does not correctly reflect the $a $-dependence. In particular,
there is no corresponding lower bound for the approximation error. It does however imply that if $\overline{l }_{q, \kappa } (\varphi ) $ is small, then the error will appear to be smaller that $O(h^{\kappa - q } ) $ for not-too-small $h > 0 . $  It can in particular be applied with basis functions $\varphi = \phi _c $ depending on a large shape parameter $c $, since one easily shows that   
$\overline{l }_{\kappa , q } (\phi _c ) \to 0 $ as $c \to \infty $ if $\phi \in \mathfrak{B }_{\kappa , N } (\mathbb{R }^n ) $ with $N > \max (n + q , \kappa ) $; 
cf. corollary \ref{corr:formulas_l_kappa-a_kappa} and its proof. 
One also shows that the constant $C_{\varepsilon } $ in (\ref{eq:rough_approx_approx_4}) (for any fixed $\varepsilon > 0 $) 
again behaves like $c^s $ for the $\phi $'s we considered in subsection \ref{subsection:c-dependence}.   
\medskip   
   
To obtain a more precise result, and in particular an asymptotic lower bound for the approximation error, we put an additional hypothesis on $a . $   

\begin{definition} \rm{We will say that 
$a \in S^q _0 ( \mathbb{R }^n ) $ is {\it asymptotically homogeneous at infinity} if 
$\partial _{\xi } a \in S^{q - \varepsilon } _0 (\mathbb{R }^n ) $ for some $\varepsilon > 0 $ and if   \begin{equation}   
\lim _{\lambda \to \infty } \frac{a(\lambda \xi ) }{\lambda ^q } =: a_{\infty } (\xi )      
\end{equation} 
exists for all $\xi \in \mathbb{R }^n \setminus 0 . $   
}   
\end{definition}   
   
\noindent We will also assume for simplicity that $\lim _{\eta \to 0 } |\eta |^{\kappa } \widehat{\varphi } (\eta ) | =: A(\varphi ) =: A $ exists: if not, one has to replace $A $ by the corresponding liminf or limsup, and replace equalities by inequalities at the appropriate places in theorem \ref{thm:cov_RBF_gen_symbols_bis} below. Recall that $A > 0 $ by the definition of the Buhmann class.      
   
\begin{lemma} \label{lemma:g_a_kappa} If $a \in S_0 ^q (\mathbb{R }^n ) $, $q > 0 $, is asymptotically homogeneous at infinity, then   
\begin{equation} \label{eq:lemma_g_a_kappa}   
\lim _{h \to 0 } \frac{ G_a ^* (\xi ; h ) - a(\xi ) }{h^{\kappa - q } | \xi |^{\kappa } } =: g_{a, \kappa }   
\end{equation}   
exists and is equal to $A^{-1 } \sum _{k \neq 0 } a_{\infty } (2 \pi k ) \widehat{\varphi } (2 \pi k ) . $     
\end{lemma}   
   
\noindent {\it Proof}. The hypotheses on $a $ easily imply that   
$$   
\lim _{h \to 0 } h^q a(\xi + 2 \pi k / h ) = a_{\infty } (2 \pi k ) ,   
$$   
 for all $k \neq 0 $ and all $\xi . $ It follows that for each $\xi \in \mathbb{R }^n $,    
$$   
\lim _{h \to 0 } \sum _k h^q (a(\xi + 2 \pi k / h ) - a(\xi ) ) \widehat{\varphi } (h \xi + 2 \pi k ) = \sum _{ k \neq 0 } a_{\infty } (2 \pi k ) \widehat{\varphi }(2 \pi k ) ,   
$$   
where the interchange of summation and limit can be justified by Lebesgue's dominated convergence theorem, using that $N > n + q . $  Since   
$$   
\frac{ G_a ^* (\xi ; h ) - a(\xi ) }{h^{\kappa - q } | \xi |^{\kappa } } = \frac{ h^q \sum _k (a(\xi + 2 \pi k / h ) - a(\xi ) ) \widehat{\varphi } (h \xi + 2 \pi k ) }{ h^{\kappa } |\xi |^{\kappa } \sum _k \widehat{\varphi } (h \xi + 2 \pi k ) }   
$$   
the lemma follows. \hfill $\Box $   
   
\begin{theorem} \label{thm:cov_RBF_gen_symbols_bis} 
Suppose that $a \in S^q _0 (\mathbb{R }^n ) $, $q > 0 $, is asymptotically homogeneous and let the initial value $f $ 
be as in theorem \ref{thm:cov_RBF_gen_symbols}, in particular, $\widehat{f }|_{\mathbb{R }^n \setminus 0 } \in \mathring{L }^1 _{\kappa } (\mathbb{R }^n ) . $ 
Then if $\kappa > q $, 
\begin{equation} \label{eq:thm:cov_RBF_gen_symbols_bis}   
\lim _{h \to 0 } h^{q - \kappa } || u_h (\cdot , t ) - u(\cdot , t ) ||_A = |g_{a , \kappa } | \int _{\mathbb{R }^n } t |\xi |^{\kappa } e^{- t {\rm Re } \, a(\xi ) } |\widehat{f }(\xi ) | \, d\xi ,   
\end{equation}   
while if $\kappa = q $ this limit equals 
\begin{equation}   
\int _{\mathbb{R }^n } \left | 1 - e^{ - t g_{a , \kappa } |\xi |^{\kappa } } \right |  e^{- t {\rm Re } \, a(\xi ) } \, |\widehat{f } (\xi ) | \, d\xi ,      
\end{equation} \label{eq:cov_RBF_gen_symbols_bis1}   
where we note that ${\rm Re } \, g_{a, \kappa } \geq 0 . $   
\end{theorem}   
   
\noindent {\it Proof}. First of all, since ${\rm Re } \, G_a ^* \geq 0 $ and $|| \widehat{s_h [f ] } - \widehat{f } ||_1 = O(h^{\kappa } ) $, the limit we have to compute is equal to   
$$  
\lim _{h \to 0 } \int _{\mathbb{R }^n } h^{- (\kappa - q ) } \left | e^{ - t (G_a ^* - a ) } - 1 \right | e^{ - t {\rm Re } \, a } |\widehat{f } | d\xi .   
$$   
If $\kappa > q $ then for each $\xi $, $G_a ^* (\xi , h ) - a(\xi ) \to 0 $ as $h \to 0 $, and consequently, using the Taylor expansion of the exponential function,   
$$   
\lim _{h \to 0 } h^{- (\kappa - q ) } \left | e^{ - t (G_a ^* - a ) } - 1 \right | = | g_{a , \kappa } | \, |\xi |^{\kappa } ,   
$$   
while if $\kappa = q $, this limit equals $| e^{ - t g_{a , \kappa } |\xi |^{\kappa } } - 1 | . $   
The dominated convergence theorem, using the estimate (\ref{eq:cov_RBF_gen_symbols1}), then shows that      
$$   
\lim _{h \to 0 } h^{- (\kappa - q ) } \int _{\mathbb{R }^n } \left | e^{ - t (G_a ^* - a ) } - 1 \right | e^{ - t {\rm Re } \, a } \, | \widehat{f } | \cdot \mathbf{1 }_{ \{ |\xi | \leq r / h \} } d\xi   
$$   
exists and is equal to the right hand side of (\ref{eq:thm:cov_RBF_gen_symbols_bis}) if $r > 0 $ is sufficiently small. The integral over $h |\xi | \geq r $ can as before be bounded by $C h^{\kappa } || \widehat{f } ||_{1 , \kappa } $ which goes to 0 when multiplied by $h^{ - (\kappa - q ) } $ if $q > 0 $, which proves the theorem. \hfill $\Box $   
   
\medskip   
   
The theorem shows that the approximation order $O(h^{\kappa - q } ) $ is exact if $g_{a , \kappa } \neq 0 $, but the latter can be 0, 
for example if $a_{\infty } (\eta ) $ is odd and $\widehat{\varphi } (\eta ) $ is even, in which case $a_{\infty } $ would have to be purely imaginary 
to comply with (\ref{cond_symbol_2})).   
A concrete example is given by the constant-coefficient transport equation $\partial _t u +v \cdot \nabla  u = 0 $, $v \in \mathbb{R }^n . $   
\medskip   
   
\begin{versionA}   
   
\noindent \textcolor{blue}{ {\bf Quelques observations sur le cas de $q = 0 $} \\ \\   
\noindent 1. Ce cas n'est pas denu\'e d'int\'er\^et: par exemple, les densit\'es de transition d'un processus de Poisson compos\'e (compuns Poisson process) satsifont une \'equation pseudo-diff parabolique du type consider\'ee ici avec un $a $ qui est born\'ee. \\ \\
\noindent 2. Si $q =  $, \c ca ne semble plus tellement essentiel de demander que ${\rm Re } \, a \geq 0 $, mais on va quand m\^eme le supposer dans un premier temps (?) \\ \\   
\noindent 3. Si $q = 0 $, on ne peut plus s\'eparer les contributions de $e^{- t G_a ^* } (\widehat{s_h [f ] } - \widehat{f } ) $ et de $ (e^{- t G_a ^* } - e^{- t a } ) f \mathbf{1 }_{h |\xi | \geq r } $ comme on a fait ci-haut, parce que celles-ci sont de l'ordre de $O(h^{\kappa } ) = O(h^{\kappa - q } ) . $ \\ \\   
\noindent 4. On peut toutefois donner par la m\'ethode utilis\'ee ci-haut une borne sup\'erieure du limsup, si on impose une plus grande d\'ecroissance de $\widehat{f }(\xi ) $ \`a l'infinie, en s'appuyant sur le 
corollaire \ref{corollary:approx_RBF}: si $\widehat{f } \in L^1 ( \max (|\xi |^{\kappa } , |\xi |^s ) $ for some $s >  \kappa $, then   
\begin{equation}   
\limsup _{h \to 0 } h^{- \kappa } || u_h ( \cdot , t ) - u ( \cdot , t ) ||_A \leq l_{\kappa } (\varphi ) || \widehat{f } ||_{\kappa } ^{\circ } + 
\int _{\mathbb{R }^n } \left( g_{a , \kappa } - \frac{1 }{2 } l_{\kappa } (\varphi ) \, a(\xi ) \right) |\xi |^{\kappa } |\widehat{f } (\xi ) | e^{- t a(\xi ) } d\xi .   
\end{equation}   
\\ \\   
\noindent 4. On r\'egarde d'un peu plus pr\`es la limite en supposant (pour simplifier) que $\widehat{f } $ est de support compact. On a alors   
\begin{eqnarray*}   
&&h^{- \kappa } || u_h ( \cdot , t ) - u ( \cdot , t ) ||_A = \int _{\mathbb{R }^n } h^{- \kappa } \left | \, e^{ - t G_a ^* } (\widehat{s_h [f ] } - \widehat{f } ) + \left( e^{- t G_a ^* } - e^{- t a } \right) \widehat{f } \,  \right | d\xi \\   
&=& \int _{\mathbb{R }^n } \left | \left( \left( \frac{\widehat{L }_1 (h \xi ) - 1 }{h ^{\kappa } } \right) e^{- t G_a ^* } + \left( \frac{e^{- t G_a ^* } - e^{-t a } }{h^{\kappa } } \right) \right) \widehat{f } (\xi ) + h^{- \kappa } \sum _{k \neq 0 } \widehat{f }(\xi + 2 \pi k h^{-1 } ) \widehat{L }_1  (h \xi ) e^{- t G_a ^* } \right | \, d\xi \\   
&=&  \int _{\mathbb{R }^n } \left | \left( \frac{\widehat{L }_1 (h \xi ) - 1 }{h ^{\kappa } } \right) e^{- t G_a ^* } + \left( \frac{e^{- t G_a ^* } - e^{-t a } }{h^{\kappa } } \right) \right | |\widehat{f } | \, d\xi + 
\sum _{k \neq 0 } \int _{\mathbb{R }^n }  \left | \widehat{f }(\xi + 2 \pi k h^{-1 } ) \widehat{L }_1  (h \xi ) e^{- t G_a ^* } \right | \, d\xi :   
\end{eqnarray*}   
si ${\rm supp } (\widehat{f } ) \subset B(0, R ) $, alors ${\rm supp } (\widehat{f } (\cdot + 2 \pi k h^{-1 } ) ) \subset B (2 \pi k h^{-1 } , R ) $ et ces supports sont disjoint pour $h $ suffisament petit. Pour le premier int\'egrale on observe que   
$$   
\left( \frac{\widehat{L }_1 (h \xi ) - 1 }{h ^{\kappa } } \right) \to - \frac{1 }{2 } l_{\kappa } (\varphi ) |\xi |^{\kappa } ,   
$$   
tandis que cette expression est uniform\'ement born\'ee (pour $h \leq h_0 $) par $C |\xi |^{\kappa } . $ (En fait, m\^eme si la partie r\'eelle de $a $ n'est pas n\'ecessairement positif, dans le cas de $q = 0 $: on a alors 
$|G_a ^* (\xi ; h ) | \leq || a ||_{\infty } $, puisque la somme des translat\'ees de $\widehat{L }_1 $ est 1. Aussi, l'in\'egalit\'e $\widehat{L }_1 (\eta ) - 1 \leq |\eta |^{\kappa } $ est valable globallement sur $\mathbb{R }^n $, puisque $\widehat{L }_1 (\eta ) $ tend vers 0 \`a l'infinie.) \\   
Ensuite, on a maintenant que   
$$
\frac{G_a ^* (\xi , h )  - a(\xi ) }{h^{\kappa } } = \sum _{k \neq 0 } a(\xi + 2 \pi k h^{-1 } ) \frac{\widehat{L }_1 (h \xi + 2 \pi k ) }{h^{\kappa } } + a(\xi ) \frac{\widehat{L }_1 (h \xi ) - 1 }{h ^{\kappa } }   
\to \left( g_{a, \kappa } - \frac{1 }{2 } l_{\kappa } (\varphi ) a(\xi ) \right) |\xi |^{\kappa }   
$$   
et donc 
$$   
\frac{e^{- t G_a ^* } - e^{- t a } }{h^{\kappa } } \to (g_{a , \kappa } - \frac{1 }{2 } l_{\kappa } (\varphi ) a(\xi ) ) \, |\xi |^{\kappa } e^{- t a (\xi ) } ,   
$$   
en se souvenant que $\kappa > 0 $ (on n'a pas pouss\'e la vice jusqu'\`a regarder le cas $\kappa = q = 0 $). Par la proposition \ref{prop_G_a} (ou plut\^ot son cons\'equence, l'in\'egalit\'e (\ref{eq:cov_RBF_gen_symbols1})),    
$$   
\left | \frac{e^{- t G_a ^* } - e^{- t a } }{h^{\kappa } } \right| |\widehat{f } (\xi ) | \leq C | \xi |^{\kappa } |\widehat{f } (\xi ) | ,   
$$   
qui est int\'egrable puisque $\widehat{f } $ est \`a support compact: le facteur exponentiel $e^{- t {\rm Re } \, a } $ est sans cons\'equence m\^eme si ${\rm Re } \, a $ n'est pas positive puisque $a $ est born\'e si $q = 0 . $ Finalement on observe que   
$$   
G_a ^* (\xi ; h ) = \sum _k a (\xi + 2 \pi k h^{-1 } ) \widehat{L }_1 (h \xi + 2 \pi k ) \to a(\xi ) , \ \ h \to 0 ,   
$$   
car $\widehat{L }_1 (h \xi + 2 \pi k ) \to 0 $ si $k \neq 0 $ ($\kappa > 0 $) et est born\'e en valeur absolue par $C |k |^{-N } $, pour chaque $\xi $ fixe et $k $ suffisament grand ($|k | > |\xi | $, disons, et alors $|k | > h \xi | $ si $h < 1 $); finalement, $a $ est born\'e. 
On peut donc appliquer la convergence domin\'ee pour concl\^ure que   
$$   
\int _{\mathbb{R }^n } \left | \left( \frac{\widehat{L }_1 (h \xi ) - 1 }{h ^{\kappa } } \right) e^{- t G_a ^* } + \left( \frac{e^{- t G_a ^* } - e^{-t a } }{h^{\kappa } } \right) \right | |\widehat{f } | \, d\xi   
\to \int _{\mathbb{R }^n } \left | g_{a, \kappa } - \frac{1 }{2 } l_{\kappa } (\varphi ) (1 + a(\xi ) ) \right| \, |\widehat{f } (\xi ) | \, |\xi |^{\kappa } e^{- t {\rm Re } a } \, d\xi   
$$   
En ce qui concerne la somme des int\'egrales ci-haut, en faisant comme avant des changement de variables $\xi \to \xi - 2 \pi k / h $ et en observant que $G_a ^* $ est $2 \pi / h $-p\'eriodique, on trouve (en utilisant le the\'eor\`eme de la convergence monotone pour justifier l'interchange de l'int\'egrale et de la somme, et en se souvenant que $\widehat{L }_1 \geq 0 $)   
\begin{eqnarray*}   
&&h^{- \kappa } \sum _{k \neq 0 } \int _{\mathbb{R }^n }  \left | \widehat{f }(\xi + 2 \pi k h^{-1 } ) \widehat{L }_1  (h \xi ) e^{- t G_a ^* } \right | \, d\xi = \int _{\mathbb{R }^n } h^{- \kappa } \left( \sum _{k \neq 0 } \widehat{L }_1  (h \xi - 2 \pi k ) \right) \left | \widehat{f }(\xi ) e^{- t G_a ^* } \right | \, d\xi \\   
&=& \int _{\mathbb{R }^n } h^{- \kappa } (1 - \widehat{L }_1 (h \xi ) ) \left | \widehat{f }(\xi ) e^{- t G_a ^* } \right | \, d\xi \\   
&\to & \int _{\mathbb{R }^n } \frac{1 }{2 } l_{\kappa } (\varphi ) |\widehat{f } | \, |\xi |^{\kappa } e^{-t {\rm Re } a } \, d\xi .   
\end{eqnarray*} \\ \\   
{\bf Conclusion} (\`a v\'erifier): Si $q = 0 $ mais $\kappa > 0 $ et si $\widehat{f } $ est \`a support compact (mais peut toujours avoir une singularit\'e en 0), alors   
\begin{equation}   
\lim _{h \to 0 } h^{- \kappa } || u_h (\cdot , t ) - u( \cdot , t ) ||_A = \int _{\mathbb{R }^n } \left( \frac{1 }{2 } l_{\kappa } (\varphi ) + \left | g_{a, \kappa } - \frac{1 }{2 } l_{\kappa } (\varphi ) (1 + a(\xi ) ) \right | \right) |\xi |^{\kappa } \, |\widehat{f } | e^{- t {\rm Re } a } \, d\xi   
\end{equation}   
une formule qui reste valable (je pense) quand on n'impose pas la condition ${\rm Re } a \geq 0 $ (preuve rapide: sinon on remplace $a $ par $a + C $, $C > || a ||_{\infty } ? $) Par rapport \`a l'estmation plus grossi\`ere ci-haut (valable uniq\'ement si ${\rm Re } a \geq 0 $, par ailleurs) on a le facteur $e^{- t {\rm Re } a } $ dans l'int\'egrale, plus le fait que   
$$   
\frac{1 }{2 } l_{\kappa } (\varphi ) + \left | g_{a, \kappa } + \frac{1 }{2 } l_{\kappa } (\varphi ) (1 - a(\xi ) \right | \leq l_{\kappa } (\varphi ) + | g_{a , \kappa } - a(\xi ) | ;   
$$   
au moins ces deux r\'esulats sont pas en contradiction, quoique le gain avec le r\'esultat plus pr\'ecis semble mineure, et on ne devrait pas en faire tout un plat dans le papier   
}
   
\end{versionA}   

\begin{remark} \rm{We include some observations on the case of $q = 0 $, which has some relevance for applications: for example, the transition probability densities of a compound Poisson process satisfy a pseudo-differential equation of the type considered here, with a symbol $a $ which is bounded together with its derivatives under suitable hypotheses on the probabilty distribution of the jumps (e.g. when the latter has moments of all, or sufficiently high, order).   
\begin{versionA}   
\textcolor{blue}{and one might want to compute these densities by solving (\ref{eq:PS_1}), as an altaernative to, say, inverting the, known, characteristic function of the process - attention: ceci implique imposer une condition initiale du type $\delta $, et on ne sait pas d\'efinir $s_h [\delta ] ... $   
}   
\end{versionA}   
   
If $a \in S^q _0 $ with $q = 0 $, we have to substract 
$\frac{1 }{2 } l_{\kappa } (\varphi ) a(\xi ) $ from the right hand side of (\ref{eq:lemma_g_a_kappa}) (where we recall that $l_{\kappa } (\varphi ) = 2 A^{-1 } \sum _{k \neq 0 } \widehat{\varphi } (2 \pi k ) $), making the limit depend on $\xi . $ Also, 
the contributions of the $L^1 $-norms of $e^{- t G_a ^* } (\widehat{s_h [f ] } - \widehat{f } ) $ et de $ (e^{- t G_a ^* } - e^{- t a } ) f \mathbf{1 }_{h |\xi | \leq r } $ are now of the same order $O(h^{\kappa } ) $, as is our estimate for the norm of $ (e^{- t G_a ^* } - e^{- t a } ) f \mathbf{1 }_{h |\xi | \geq r } $, and cannot be separated anymore as before. We can nevertheless impose a stronger decay at infinity to control the latter, and use corollary \ref{corollary:approx_RBF} to show that if $\widehat{f } \in L^1 ( \max (|\xi |^{\kappa } , |\xi |^s ) $ for some $s >  \kappa $, then   
\begin{equation} \label{eq:S^0}   
\limsup _{h \to 0 } h^{- \kappa } || u_h ( \cdot , t ) - u ( \cdot , t ) ||_A \leq l_{\kappa } (\varphi ) || \widehat{f } ||_{\kappa } ^{\circ } + 
\int _{\mathbb{R }^n } \left( g_{a , \kappa } - \frac{1 }{2 } l_{\kappa } (\varphi ) \, a(\xi ) \right) |\xi |^{\kappa } |\widehat{f } (\xi ) | e^{- t {\rm Re } \, a(\xi ) } d\xi .   
\end{equation}   
A further analysis shows that if, for example, $\widehat{f } $ has compact support, then we have the more precise result
\begin{equation}   
\lim _{h \to 0 } h^{- \kappa } || u_h (\cdot , t ) - u( \cdot , t ) ||_A = \int _{\mathbb{R }^n } \left( \, \frac{1 }{2 } l_{\kappa } (\varphi ) + \left | g_{a, \kappa } - \frac{1 }{2 } l_{\kappa } (\varphi ) (1 + a(\xi ) ) \right | \, \right) |\xi |^{\kappa } \, |\widehat{f } | e^{- t {\rm Re } \, a } \, d\xi ,   
\end{equation}  
a result which remains true if the real part of $a \in S^0 _0 $ is not necessarily positive, and which implies (\ref{eq:S^0}) if it is. Since the gain with respect to the latter seems relatively modest, we skip the details.   
}   
\end{remark}   
   
   
We cannot deduce from lemma \ref{lemma:g_a_kappa} a precise approximate approximation  estimate 
similar to theorem \ref{thm:approx_approx_kappa > 2} since the limit isn't "uniform in $h |\xi | $". In fact, as we will show now, lower order parts of the symbol will also occur in such an estimate. We 
consider the case of a constant coefficient partial differential operator, $p = p(\xi ) $ a polynomial of degree $q \in \mathbb{N } $ with non-negative real part: ${\rm Re } \, p(\xi ) \geq 0 . $ 
We write   
\begin{equation} \label{eq:poly_homogeneous}   
p = \sum _{j = 0 } ^q p_{q - j } ,   
\end{equation}   
with $p_{q - j } (\xi ) $ homogeneous of degree $q - j . $ Then $G_{p_{\nu } } ^* (\xi ; h ) = h^{-\nu } G_{p_{\nu } } (h\xi ) $, where   
$$   
G_{p_{\nu } } (\eta ) := \sum _k p_{\nu } (\eta + 2 \pi k ) \widehat{L }_1 (\eta + 2 \pi k ) .   
$$   
and, as is easily verified,   
$$   
\lim _{\eta \to 0 } \frac{G_{p_{\nu } } (\eta ) - p_{\nu } (\eta ) }{|\eta |^{\kappa } } = g _{p_{\nu } , \kappa } ,   
$$   
for $\nu > 0 $, while $G_{p_0 } = p_0 $ since $p_0 $ is a constant. 
Since   
$$   
G_p ^* (\xi ; h ) - p(\xi ) = \sum _{j = 0 } ^{q - 1 }  h^{- (q - j ) } \left( G_{p_{q - j } } (h\xi ) - p_{q - j } (h \xi ) \right) ,   
$$   
it follows that for each $\varepsilon > 0 $ there exists a $\rho = \rho (\varepsilon) > 0 $ such that   
$$   
\sup _{| h\xi | \leq \rho } \frac{ | G_p ^* (\xi ; h ) - p(\xi ) | }{h^{\kappa } |\xi |^{\kappa } } \leq \sum _{j = 0 } ^{q - 1 } h^j \left( | g_{p_{q - j } , \kappa }  | + \varepsilon \right) .      
$$   
Since ${\rm Re } \left( p (\xi ) - G_p ^* (\xi ; h ) \right) \leq |h \xi |^{\kappa } {\rm Re } \, p (\xi ) < \varepsilon {\rm Re } \, p(\xi ) $ if $h |\xi | $ is sufficiently small, the arguments above can be used to prove the following estimate:   
   
\begin{theorem} \label{thm:approx_approx_polyhomogeneous} For all $\varepsilon > 0 $ there exists a $C_{\varepsilon } > 0 $ such that if $\widehat{f } |_{\mathbb{R }^n \setminus 0 } \in \mathring{L }^1 _{\kappa } (\mathbb{R }^n ) \cap \mathring{L }^1 _s (\mathbb{R }^n ) $,   
\begin{eqnarray*}     
|| u_h (\cdot , t ) - u (\cdot , t ) ||_A &\leq &\sum _{j = 0 } ^{q - 1 } h^{\kappa - q + j } \left( |  g_{p_{q - j } , \kappa } | + \varepsilon \right) \int _{\mathbb{R }^n } |\xi |^{\kappa } \, |\widehat{f } (\xi ) | e^{- (1 - \varepsilon ) {\rm Re } p(\xi ) } \, d\xi  \\   
&& + \overline{l }_{\kappa } h^{\kappa }  || \widehat{f } ||^{\circ } _{\kappa } + C_{\varepsilon } h^s || \widehat{f } ||^{\circ } _s ,      
\end{eqnarray*}   
where $ |  g_{p_{q - j } , \kappa } | + \varepsilon $ can be replaced by $(1 + \varepsilon ) | g_{p_{q - j } , \kappa } | $ if $g_{p_{q - j } , \kappa } \neq 0 . $   
\end{theorem}   

\begin{examples} \label{ex:pseudo-diff_evolution_eq} \rm{We finally mention 
some examples of pseudo-differential evolution equations 
which are of interest for applications, and to which our results apply.      
\medskip   
   
\noindent (i) The Kolmogorov - Fokker - Planck equation associated to a L\' evy process $(X_t )_{t \geq 0 } $ on $\mathbb{R }^n . $ Recall that  according to the L\'evy - Khintchine theorem such a process is completely characterized by its characteristic function, 
$\mathbb{E } \left( e^{ i (\xi , X_t ) } \right) = e^{t \psi (\xi ) } $ with   
$$   
\psi (\xi ) = i (\mu , \xi ) - \frac{1 }{2 } (\Sigma \, \xi , \xi ) + \int _{\mathbb{R }^n \setminus 0 } \left( e^{i (x , \xi ) } - 1 - i (x , \xi ) \chi (x ) \right) d\nu (x ) ,   
$$   
where $\Sigma $ is a positive semi-definite linear operator and where $\nu  $ is a positive Borel measure on $\mathbb{R }^n \setminus 0 $ such that   
$$   
\int _{\mathbb{R }^n \setminus 0 } (|x |^2 \wedge 1 ) d\nu (x)  < \infty ,   
$$   
called the L\'evy measure; here $\chi $ is a compactly supported function which is equal to 1 on a neighborhood of 0, and which can be taken smooth, if necessary.   
   
If, for a given $f $, we let $u(x, t ) = \mathbb{E } (f (x + X_t ) ) $, then $u $ satisfies (\ref{eq:PS_1}) with $a(\xi ) := - \psi (\xi ) $ and initial value $f . $ Note that $a(\xi ) $ satisfies (\ref{cond_symbol_2}) since   
$$   
{\rm Re } \, \psi (\xi ) = - \frac{1 }{2 } \left( \Sigma \, \xi , \xi \right) + \int _{\mathbb{R }^n \setminus 0 } ( \cos (x \xi ) -1 ) d \nu (x ) \leq 0 .   
$$   
Under appropriate hypothese on the L\'evy-measure $\nu $ one can derive symbol-type estimates for $a(\xi ) . $ For example, when 
$d\nu (x ) = |x |^{- q } h(x) dx $ with $q < n + 2 $, and $h(x) $ a rapidly decreasing continuous function, 
then $a \in S^2 _0 $ if $\Sigma \neq 0 $, and in $S^{q - n } _0 $ if $V = 0 $: cf. remark \ref{remark:symbol_estimates} in Appendix A below. Examples of such processes are the jump-diffusion processes and the CGMY-processes of mathematical finance, which were treated numerically in \cite{RB_RC}, \cite{RC} and \cite{RC_SH} with different choices of basis functions (respectively the multi-quadric, inverse multi-quadric and the cubic spline). If $\Sigma = 0 $ and if $\nu $ is a finite measure on $\mathbb{R }^n $ of total mass $\lambda $, 
we have the special case of a (multi-dimensional) compound Poisson process with intensity $\lambda $ and drift $\mu + \int _{\mathbb{R }^n } x  \chi d\nu (x) $, $\lambda ^{-1 } \nu $ then being the probability law of the jumps. The symbol will then be in $C_b ^{\lceil \kappa \rceil + n + 1 } (\mathbb{R }^n ) $ if $\nu $ has moments of this order.   
\medskip   
   
\noindent (ii) The special case of {\it symmetric stable} L\' evy-processes leads to the fractional heat equation:   
$$   
\partial _t u + (- \Delta )^s u = 0 ,   
$$   
with $s \in (0, 1 ) . $ Here the symbol, $a(\xi ) = |\xi |^{2s } $ is not smooth in 0, so our theory does not apply directly, but it applies to regularized versions with $|\xi |^{2s } $ replaced by $\chi (\xi ) |\xi |^{2s } $, where $\chi \in C^{\infty } (\mathbb{R }^n ) $ is equal to 1 outside of a small ball. The theory of this section can nevertheless be extended to certain families of symbols having singularities in $\xi = 0 . $ The presence of such singularities cause a number of technical problems, notably with regards to the initial values $f $ which are allowed in theorem \ref{thm:cov_RBF_gen_symbols} and the other theorems of this section: in the case of the fractional heat equation, $e^{ - t |\xi |^{2s } } \widehat{f } $ needs to be defined as a tempered distribution, so it becomes natural to require that $\widehat{f } |_{\mathbb{R }^n \setminus 0 } \in L^1 ( \mathbb{R }^n , (|\xi |^{2s } \wedge 1 ) d\xi ) . $ On the other hand one can use lemma \ref{lemma:Stein} to show that the inverse Fourier transform of $e^{ - t |\xi |^{2s } } $ decays as $(1 + |x | )^{n + 2s } $ at infinity, so $f $ will have to be taken integrable with respect to this weight. The details of this case will be treated elsewhere. Note that once such a theory is in place, we can consider generalizations of theorem \ref{thm:approx_approx_polyhomogeneous} to general polyhomogeneous symbols.   
\medskip   
   
\noindent (iii) Certain non-parabolic operators, such as the transport equation $\partial _t u = v \cdot \nabla u $, 
the so-called "half wave equation" $\partial _t u = i \sqrt{ - \Delta } u $, modified so as to make its symbol smooth in a neighborhood of 0, or the free Schr\"odinger equation, $i \partial _t u = - \Delta u $: for these, the symbol $a(\xi ) $ is purely imaginary.
}   
   
\end{examples}   
   
\begin{versionA}   
   
\textcolor{blue}{Discussion of the numerical results of \cite{RB_RC} in the light of the results of this paper - not yet complete} In \cite{RB_RC} used the RBF-scheme to compute put option prices in the so-called CGMY model (called temperate L\'evy models or some such thing in Cont and Tankov \cite{CT}). The RBF used was the multi-quadric in dimension 1, which is in $\mathfrak{B }_{2, N } $ for all $N > 1 . $ Our results indicates that the over the range of $h $' used, the scheme converges at a rate of $h^2 $, which is odd in the light of theorem \ref{?}: the operator is of order 2 if $\Sigma \neq 0 $ and of order $Y $ otherwise, and the rate if convergence according to the theorem therefore is either 0 or $2 - Y . $ We examine whether the numerical results of \cite{RB_RC} can be explained as an approximate approximation phenomenon.   
   
The initial vale $f(x) = \max (K - e^x , 0 ) $, corresponding to a put in log-price coordinates. Its Fourier transform can readily be computed as   
$$   
\widehat{f }(\xi ) = - \frac{i K^{1 - i \xi } }{(\xi + i )(\xi + i 0+ ) } ,   
$$   
from which it is apparent that $\widehat{f } $ is integrable.   
   
\end{versionA}

\appendix   
   
\section{\bf Proof of theorem \ref{thm_L_1} 
}   
\label{Appendix:Lagrange_function}   
     
   
We prove the existence and main properties of the cardinal function associated to a basis function $\varphi \in \mathfrak{B }_{\kappa , N } (\mathbb{R }^n ) $ as stated in theorem \ref{thm_L_1}. As already mentioned, this was done by Buhmann \cite{Bu1}, \cite{Bu_book} for a more restricted class of radial basis functions.   
%
%
The main difference in our treatment and that of \cite{Bu1} 
is the use of 
lemma \ref{lemma:Stein} below, relating the decay at infinity of the Fourier transform of a function with its behavior in 0, which allows one to go beyond the  class of basis functions considered by Buhmann. 

\medskip   
  
Before embarking upon the proof, it may be interesting to observe that the estimates (\ref{eq:wB2}) are similar to the symbol conditions of pseudodifferential calculus, except that the latter concern the behavior at infinity\footnote{Indeed, if (\ref{eq:wB3}) were required for all orders $\alpha $ (with constants which may then depend on $\alpha $), then $\chi (\xi ) \widehat{\varphi } (\xi /|\xi |^2 ) \in S^{\kappa } _{1, 0 } (\mathbb{R }^n ) $, where $\chi $ is a $C^{\infty } $-function such that $1 - \chi $ is compactly supported, and $S^p _{1, 0 } (\mathbb{R }^n ) $ is the standard symbol class of order $p $ (cf. \cite{Stein}). } instead of at 0. From this point of view, (\ref{eq:wB3}) corresponds to having an elliptic symbol, whence our terminology. \medskip  
   
We note that as a consequence of conditions (ii) and (iii) of definition \ref{def:wBuhmann class},      
\begin{equation} \label{eq:wB4}   
|\partial _{\eta } ^{\alpha } (\widehat{\varphi }^{ \, -1 } ) | \leq C |\eta |^{\kappa - |\alpha | } , \ \ |\eta | \leq 1 , |\alpha | \leq n + \lfloor \kappa \rfloor + 1 .      
\end{equation}   
\medskip   
   

Turning to the proof of theorem \ref{thm_L_1}, following Buhmann we start by {\it defining} $L_1 $ as the inverse Fourier transform of the right hand side of (\ref{Fourier_L_1}), observing that since the latter is an integrable function, by definition \ref{def:wBuhmann class}(iv),  $L_1 $ is a well-defined continuous function. We first show that $L_1 (x) $ has the proper decay at infinity.   
   
\begin{theorem} \label{thm:L_wBuhmann} Let $\varphi \in \mathfrak{B }_{\kappa , N } (\mathbb{R }^n ) $, and let   
\begin{equation}   
L_1 := \mathcal{F }^{-1 } \left( \frac{\widehat{\varphi } (\cdot ) }{ \sum _{k \in \mathbb{Z }^n } \widehat{\varphi } (\cdot + k ) } \right) .   
\end{equation}   
Then there exists a positive constant $C $ such that        
\begin{equation} \label{eq:L_wBuhmann}   
|L_1 (y) | \leq C (1 + |y| )^{-\kappa - n } , \ \ y \in \mathbb{R }^n .   
\end{equation}   
\end{theorem}   
   
The proof will use the following lemma, which basically is a special case of a classical estimate for kernels of convolution operators: see Stein \cite{Stein}, proposition 2 of Chapter VI, section 4.4. 
   
\begin{lemma} \label{lemma:Stein} Let $p > - n  $ and let $a \in C^{\lfloor p \rfloor + n + 1 } (\mathbb{R }^n \setminus 0 ) $ be supported in some ball $B(0, R ) $ such that\footnote{Note that $\lfloor p \rfloor + n + 1 \geq 1 $ since $p > - n . $ }         
\begin{equation} \label{eq:Stein1}   
|\partial _{\xi } ^{\alpha } a(\xi ) | \leq C |\xi |^{p - |\alpha | } , \ \ , 
|\xi | \leq R , \ |\alpha | \leq \lfloor p \rfloor + n + 1 .         
\end{equation}   
Then the inverse Fourier transform $k = \mathcal{F }^{-1 } (a) $ satisfies     
\begin{equation} \label{eq:Stein_at_inf}   
|k(x) | \leq C_1 (1 + |x|)^{ - p - n } , \ \ x \neq 0 ,     
\end{equation}   
with a constant $C_1 \leq c_n C $, where $c_n $ only depends on $n . $   
\end{lemma}   
   
   
Stein in fact shows 
that if (\ref{eq:Stein1}) is satisfied at all orders, without $a $ necessarily being compactly supported, 
then $k $ can be identified with a $C^{\infty } $-function away from 0, satisfying $|\partial _x ^{\alpha }  k (x) | \leq C_{\alpha } |x |^{ - p - n - |\alpha | } $ for all $\alpha $ and all $x . $ This result is stated and proven there for $p = 0 $, but the proof generalizes to any $p > -n . $ We only need this estimate for $k(x) $ itself, in which case we only need (\ref{eq:Stein1}) for the limited number of derivatives of $a $ indicated, and we furthermore only need it for large $|x | $ (note that if $a $ has compact support, $k $ is continuous, even $C^{\infty } $, and Stein`s estimate for $k(x) $ at 0 becomes trivial). The proof in \cite{Stein} uses the Paley-Littlewood decomposition. An elementary prove of lemma \ref{lemma:Stein} can be given by writing   
$$   
k(x) = (2 \pi )^{-n } \int _{\mathbb{R }^n } \, \chi (|x| \xi ) a(\xi ) e^{i (x , \xi ) } \, d\xi + (2 \pi )^{-n } \int _{\mathbb{R }^n } \, (1 - \chi (|x| \xi ) ) a(\xi ) e^{i(x , \xi ) } \, d\xi .   
$$   
where $\chi \in C^{\infty }(\mathbb{R }^n ) $ with bounded derivatives such that $\chi (\xi ) = 0 $ for $|\xi | \leq 1 $, $\chi (\xi ) = 1 $ for $|\xi | \geq 2 $, and integrating the first integral by parts $\lfloor p \rfloor + n + 1 . $   
   
\begin{versionA}   
\textcolor{blue}{(detailed verification in Appendix C)}   
\end{versionA}   
\medskip   
   
\noindent {\it Proof of theorem \ref{thm:L_wBuhmann}}. Let $\chi _0 \in C^{\infty } _c (\mathbb{R }^n ) $ such that $\chi _0 (\eta ) = 1 $  in a neighbourhood of $0 $ and ${\rm supp }(\chi _0 ) \subset (- \pi, \pi ) ^k . $ For $k  \in \mathbb{Z }^n $, define $\chi _k $ by $\chi _k (\eta ) := \chi (\eta + 2 \pi k ) $ and note that the supports of the $\chi _k $ are disjoint. Finally, let $\chi _c := 1 - \sum _k \chi _k $ ("c" for "complement"), so that $\chi _c $ together with the $\chi _k $'s form a partition of unit. Then   
\begin{equation}   
L_1 (x) = \ell _c (x) + \sum _{k \in \mathbb{Z }^n } \ell _k (x) ,   
\end{equation}   
where   
\begin{equation}  
\ell _k = \mathcal{F }^{-1 } \left( \chi _k (\eta ) \frac{\widehat{\varphi } (\eta ) }{\sum _{\nu } \widehat{\varphi }(\eta + 2 \pi \nu ) } \right), \ \ k \in \mathbb{Z }^n \ \mbox{or } k = c .     
\end{equation}   
We examine the decay in $x $ of each term separately. 
\medskip   
   
\noindent {\it Decay of $\ell _c $}. The function $\chi _c (\eta ) / \sum _k \widehat{\varphi} (\eta + 2 \pi k ) $ is in $C_b ^{\lfloor \kappa \rfloor + n + 1 } (\mathbb{R }^n ) $, 
since the denominator is a strictly positive periodic function which is $C_b ^{\lfloor \kappa \rfloor + n + 1 } $ on the complement of $(2 \pi \mathbb{Z } )^n $ and therefore on the support of $\chi _c . $ 
Multiplying with $\widehat{\varphi } $, we find that $\widehat{\ell }_c (\eta ) $ 
is $C^{\lfloor \kappa \rfloor + n + 1  } $ with integrable derivatives of all orders, which implies by the usual integration by parts argument that $|\ell _c (x) | \leq C (1 + |x| )^{-(\lfloor \kappa \rfloor + n + 1 ) } \leq C (1 + |x| )^{- \kappa - n } . $   
   
\begin{versionA}   
\textcolor{blue}{   
\begin{remark} \textcolor{blue}{regarding theorem \ref{thm:L_wBuhmann} (perhaps put this remark at the end of the proof?} \noir \rm{If $\kappa \in \mathbb{N } $, we need $\kappa + n $ derivatives to obtain the upper bound of $C (1 + |x| )^{- \kappa - n } $ for $\ell _c (x) . $ Hence we only need the estimate at infinity, (\ref{eq:wB1}), for derivatives up till order $\lceil \kappa \rceil + n $, which is equal to $\lfloor \kappa \rfloor + n + 1 $ if $\kappa \notin \mathbb{N } $, but one less if $\kappa \in \mathbb{N } . $ However, to obtain the correct degree of decay of the $\ell _k (x) $ we do need $\lfloor \kappa \rfloor + n + 1 $ derivatives in (\ref{eq:wB2}), see below, mainly because of lemma \ref{lemma:Stein}. This somewhat contrasts with \cite{Bu1}, who only needs $\lceil \kappa \rceil + n $ derivatives there also, but who has stronger assumptions on $\widehat{\varphi } (\eta ) $ near 0:  a homogeneous top order term plus an error which is (asymptotically) sufficiently small: $\widehat{\varphi }_1 (\eta ) = c |\eta |^{\kappa } \left( 1 + h(\eta ) \right) $ with $|h^{(\alpha ) } | \leq C_{\alpha } |\eta |^{\varepsilon - |\alpha | } $ as $\varepsilon \to 0 $, {\it plus } the condition $\varepsilon | \lceil \kappa \rceil - \kappa $ (which is void if $\kappa \in \mathbb{N } $). The extra derivative is thus a price to pay for a greater flexibility. (In practice it will of course matter little: for all examples of RBFs we are aware of $\widehat{\varphi } $ is $C^{\infty } $ away from the origen and satisfies (\ref{eq:wB1}) and (\ref{eq:wB2}) to all orders.)   
}   
\end{remark}   
}   
\end{versionA}   
   
Note that $\widehat{L }_1 (\eta ) $ is at best $C^{\lfloor \kappa \rfloor} $ in the points of $2 \pi \mathbb{Z }^n $, so integration by parts will not give the required decay for of $\ell _k $, $k \in \mathbb{Z }^n . $ We use lemma \ref{lemma:Stein} instead.  
\medskip   
   
\noindent {\it Decay of $\ell _0 $}. Since      
$$   
\widehat{\ell }_0 (\eta ) - \chi _0 (\eta ) = \chi_0 (\eta ) \left( \widehat{L }_1 (\eta ) - 1 \right)   
=  - \chi _0 (\eta ) \left( \frac{ \widehat{\varphi } (\eta ) ^{-1 } \sum _{k \neq 0 } \widehat{\varphi } (\eta + 2 \pi k ) }{1 + \widehat{\varphi }(\eta ) ^{-1 } \sum _{k \neq 0 } \widehat{\varphi } (\eta + 2 \pi k ) } \right)     
$$   
and since $\sum _{k \neq 0 } \widehat{\varphi } (\eta + 2 \pi k ) $ is $C^{\lfloor \kappa \rfloor + n + 1 } $ on the suport of $\chi _0 $, the estimates 
(\ref{eq:wB4}) easily imply that $\widehat{\ell }_0  (\eta ) - \chi _0 (\eta ) $ satisfies condition (\ref{eq:Stein1}) of lemma \ref{lemma:Stein} with $p = \kappa $ (for $\psi := \widehat{\varphi } ^{-1 } \sum _{k \neq 0 } \widehat{\varphi } (\cdot + 2 \pi k ) $ does,  and then also $\psi / (1 + \psi ) $). It follows that $|\ell _0 (x) | \leq C (1 + |x |)^{- \kappa - n } $, since $\mathcal{F }^{-1 } (\chi _0 ) $ is rapidly decreasing.   
\medskip   
   
\noindent {\it Decay of $\ell _k $, $k \neq 0, c $}. This is similar, except that we have to pay attention to the size of the constant in front of the $(1 + |x| )^{ - \kappa - n } . $ The Fourier transform $\widehat{\ell }_k (\eta ) $ will now be supported near $\eta = - 2 \pi k . $ Shifting by $2 \pi k $, we see that   
$$   
\widehat{\ell }_k (\eta - 2 \pi k ) = \chi _0 (\eta ) \varphi (\eta - 2 \pi k ) \frac{\widehat{\varphi } (\eta )^{-1 } }{1 + \sum _{\nu \neq 0 } \widehat{\varphi } (\eta )^{-1 } \widehat{\varphi } (\eta + 2 \pi \nu ) }   
$$   
is supported in a small neighbourhood of 0, with derivatives of order $|\alpha | \leq \lfloor \kappa \rfloor + n + 1  $ bounded by   
$C (1 + |k | )^{-N } |\eta |^{\kappa - |\alpha | } $, with $C $ independent of $k . $ Lemma \ref{lemma:Stein} then implies that   
$$   
|\ell _k (x) | = \left |\ell _k (x) e^{2 \pi i (k, x ) } \right | \leq C (1 + |k| )^{-N } (1 + |x |)^{- \kappa - n } .   
$$   
Since $N > n $ by assumption, summation over $k \in \mathbb{Z }^n $ completes the proof. 
   
\hfill $\Box $          
\medskip   
   
The same arguments prove theorem \ref{lemma:a(D)L_1}: $\widehat{a(D)(L_1)} = a \widehat{\ell}_c + \sum _k a \widehat{\ell _k } $, and $a \widehat{\ell }_c $ has $\lceil \kappa \rceil + n + 1 $ derivatives which are integrable since $N + q < n $, as has $a \chi _0 . $ Finally, $a \chi _0 (\widehat{L } - 1 ) $ and $a \chi _k \widehat{L }_1 $, $k \neq 0 $ satisfy the hypothesis (\ref{eq:Stein1}) of lemma \ref{lemma:Stein}, the latter with constants bounded by $C |k |^{-(N - q ) } $ which are summable since $N > n + q . $   
\medskip

We 
continue with the proof of theorem \ref{thm_L_1}.   
The expressions for $\widehat{\ell }_k $ above also show that $\widehat{L }_1 $ satisfies the Strang-Fix conditions (\ref{eq:SFC}). Once we have defined $L_1 $ through its Fourier transform, it is immediate to check that $L_1 (k) = \delta _{0k } $ for $k \in \mathbb{Z }^n $: indeed, by the $2 \pi $-periodicity of the denominator, writing the integral over $\mathbb{R }^n $ as a sum of integrals over translates of $(-\pi , \pi )^n $,   
\begin{eqnarray*}   
L_1 (k) &=& \int _{\mathbb{R }^n } \frac{\widehat{\varphi }(\eta ) }{\sum _{\nu } \widehat{\varphi }(\eta + 2 \pi \nu ) } e^{i k \eta } \frac{d \eta }{(2 \pi )^n } \\   
&=& \int _{(- \pi , \pi )^n } \frac{\sum _{\nu ' } \widehat{\varphi } (\eta + 2 \pi \nu ' ) }{\sum _{\nu } \widehat{\varphi }(\eta + 2 \pi \nu ) } e^{i k \eta } \frac{d \eta }{(2 \pi )^n } \\   
&=& \int _{(- \pi , \pi )^n } e^{i k \eta } \frac{d \eta }{(2 \pi )^n } \\   
&=& \delta _{0k } .   
\end{eqnarray*}   
It remains to recognise $L_1 $ as a sum of translates of $\varphi . $ 
To show this we first write the denominator of (\ref{Fourier_L_1}) 
as a Fourier series:   
\begin{equation} \label{eq:Fourier-series}   
\left( \sum _k \widehat{\varphi }(\eta + 2 \pi k ) \right) ^{-1 } = \sum _k c_k e^{i k \eta } .   
\end{equation}   
One verifies by the similar arguments as the ones of the proof of theorem \ref{thm:L_wBuhmann} that   
\begin{equation} \label{eq:est_c_k}   
|c_k | \leq C (1 + |k |)^{- \kappa - n } ,            
\end{equation}   
write   
$$   
c_k = (2 \pi )^{-n } \int _{(- \pi , \pi )^n } \frac{\chi _0 (\eta ) }{\sum _{\nu } \widehat{\varphi } (\eta + 2 \pi \nu ) } e^{i (\eta , k ) } d\eta + (2 \pi )^{-n } \int _{(- \pi , \pi )^n } \frac{1 - \chi _0 (\eta ) }{\sum _{\nu } \widehat{\varphi } (\eta + 2 \pi \nu ) } e^{i (\eta , k ) } d\eta ,   
$$   
and estimate the first integral using lemma \ref{lemma:Stein} and the second by integrating by parts.   
   
\begin{versionA}   
   
\noindent \textcolor{blue}{Indeed, the, periodic, function   
\begin{eqnarray*}   
a(\eta ) &:=& \left( \sum _{\nu } \widehat{\varphi } (\eta + 2 \pi \nu ) \right)^{-1 } \\   
&=& \varphi (\eta )^{-1 } \left( 1 + \sum _{\nu \neq 0 } \widehat{\varphi } (\eta )^{-1 } \widehat{\varphi } (\eta + 2 \pi \nu ) \right)^{-1 }     
\end{eqnarray*}   
satisfies the estimates (\ref{eq:Stein1}) in a neighborhood of 0 and is $C^{\lfloor \kappa \rfloor + n + 1 } $ on the complement of $(2 \pi \mathbb{Z } )^n . $ Since       
$$   
(2 \pi )^n c_k = \int _{(- \pi , \pi )^n } \chi_0 (\eta ) a(\eta ) e^{ - i (x, k ) } d\eta + \int _{(- \pi , \pi )^n } (1 - \chi_0 (\eta ) ) a(\eta ) e^{- i (x, k ) } d\eta ,   
$$   
estimating the first integral using lemma \ref{lemma:Stein} and the second by integrating by parts yields (\ref{eq:est_c_k})   
}   
\medskip   
   
\end{versionA}   
It follows that (\ref{eq:Fourier-series}) converges absolutely. We then claim that   
\begin{equation} \label{eq:formula_L_1}   
L_1 (x) = \sum _k c_k \varphi (x - k ) ,   
\end{equation}   
where the series converges absolutely and uniformly on compacta, by (\ref{eq:est_c_k}), since $\varphi (x) $ grows at most as $(1 + |x \ )^{\kappa - \varepsilon } $, by assumption. Formally, (\ref{eq:formula_L_1}) follows by writing   
\begin{eqnarray*}   
L_1 (x) = \int _{\mathbb{R }^n } \left( \sum _k c_k e^{i k \eta } \right) \widehat{\varphi } (\eta ) e^{i \eta x } \frac{d \eta }{(2 \pi )^n } = \sum _k c_k \varphi (x + k ) ,   
\end{eqnarray*}   
except that the final step does not make sense for an arbitrary $\varphi \in \mathfrak{B }_{\kappa , N } (\mathbb{R }^n ) $ since $\widehat{\varphi } (\eta ) $ will not be integrable in 0 if $\kappa > n $ and even if it is, when $\kappa < n $, 
it might differ from integration against the function $\widehat{\varphi } (\eta ) $ by a distribution supported in 0.   
   
\begin{versionA}   
\medskip   
   
\noindent \textcolor{blue}{{\bf Comment}: If I understood correctly, Buhman in his book seems to claim that, as tempered distributions, the {\it function }   
$$   
\frac{\widehat{\varphi } (\cdot ) }{ \sum _{k \in \mathbb{Z }^n } \widehat{\varphi } (\cdot + k ) }   
$$   
equals the {\it distribution}   
$$   
\left( \sum _k c_k e^{i (\cdot , k ) } \right) \widehat{\varphi } ,   
$$   
but this does not seem at all obvious to me: for one thing, we are not multiplying the distribution with a $C^{\infty } $-function, so one has to show the product makes sense, which requires thinking about the order of $\widehat{\varphi } $ as a distribution. For another, even if true on $\mathbb{R }^n \setminus 0 $ (for that one has to discuss convergence of the Fourier series?) that doesn't make it true as distributions on $\mathbb{R }^n $, which is what is needed for the inverse Fourier transform. I give a proof of (\ref{eq:formula_L_1}) below, which is however a bit long and technical. Can one perhaps find a shorter proof?}   
\medskip   
   
\end{versionA}   
We have to carefully distinguish between the tempered distribution $\widehat{\varphi } $ and the locally integrable function $\eta \to \widehat{\varphi } (\eta ) $ with which it can be identified on $\mathbb{R }^n \setminus 0 . $ The relation between the two is given by the following identity: 
there exist constants $c_{\alpha } $, $|\alpha | \leq \lceil \kappa \rceil - 1 $ such that for all $\psi \in \mathcal{S }(\mathbb{R }^n ) $,   
\begin{eqnarray} \label{eq:widehat_varphi}   
\langle \widehat{\varphi } , \psi \rangle &=& \int _{|\eta | \leq 1 } \, \widehat{\varphi } (\eta ) \left( \psi (\eta ) - \sum _{|\alpha | \leq \lfloor \kappa \rfloor - n } \frac{\psi ^{(\alpha ) }(0) }{\alpha ! } \eta ^{\alpha } \right) \, d\eta  \\   
&& + \int _{|\eta | \geq 1 } \widehat{\varphi } (\eta ) \psi (\eta ) \, d\eta \label{eq:FT_of_varphi} \nonumber \\  
&& + \sum _{|\alpha | \leq \lceil \kappa \rceil - 1 } (-1 )^{|\alpha | } c_{\alpha } \psi ^{(\alpha ) }(0) , \nonumber   
\end{eqnarray}   
where the sum is interpreted as empty if $\kappa < n . $ Indeed, the first integral converges since $|\eta |^{\lfloor \kappa \rfloor - n + 1 } \widehat{\varphi } (\eta ) $ has an integrable singularity at 0. The sum of the two integrals on the right defines a tempered distribution. If we denote this distribution by $\Lambda _{\widehat{\varphi } } $ then the restriction of $\Lambda _{\widehat{\varphi } } $ 
to $\mathbb{R }^n \setminus 0 $ can be identified with the function $\widehat{\varphi } (\eta ) . $ The difference $\widehat{\varphi } - \Lambda _{\widehat{\varphi } } $ is then supported in 0, and therefore a linear combination $\sum _{|\alpha | \leq p } c_{\alpha } \delta _0 ^{(\alpha ) } $ of derivatives of the delta distribution in 0. To bound $p $, we use the following lemma, whose proof we postpone till the end of his section:    
   
\begin{lemma} \label{lemma:growth_Lambda_varphi} If $\kappa \geq n $, then the inverse Fourier transform $\mathcal{F }^{-1 } \left( \Lambda _{\widehat{\varphi } } \right) $ is a continuous function which is bounded by $C (|x |^{\kappa - n } +1 ) $ for non-integer $\kappa $ 
and by $C (|x |^{\kappa - n } \log |x | + 1 ) $ if $\kappa $ is a positive integer. 
\end{lemma}   
   
\noindent Since the inverse Fourier transform of $\sum _{|\alpha | \leq p } c_{\alpha } \delta _0 ^{(\alpha ) } $ is a polynomial of order $p $, and since, by assumption, $\varphi (x) $, has polynomial growth of order strictly less than $\kappa $, it follows that $p < \kappa $, which is equivalent to $p \leq \lceil \kappa \rceil - 1 . $
\smallskip   
   
We now use (\ref{eq:widehat_varphi}) 
to prove that (\ref{eq:formula_L_1}) holds as tempered distributions, that is, if $\psi \in \mathcal{S }(\mathbb{R }^n ) $, then   
\begin{equation} \label{eq:proof-FTL}   
\langle L_1 , \widehat{\psi } \rangle = \left \langle \sum _k c_k \varphi ( \cdot + k ) , \widehat{\psi } \right \rangle .   
\end{equation}   
If we let   
$$   
\Psi (\eta ) := \frac{\psi (\eta ) }{\sum _k \widehat{\varphi } (\eta + 2 \pi k ) } . 
$$   
then $\Psi $ is $C^{\lfloor \kappa \rfloor } $ if $\kappa \notin \mathbb{N } $, and $C^{\kappa - 1 , 1 } $ if $\kappa \in \mathbb{N }^* $, with all its derivatives rapidly decreasing. To obtain a function in the  Schwartz class $\mathcal{S }(\mathbb{R }^n ) $ we convolve with $\chi _{\varepsilon } (x) := \varepsilon ^{-n } \chi (x/\varepsilon ) $, where $\chi \in C^{\infty }_c (\mathbb{R }^n ) $ with integral 1. Let $\Psi _{\varepsilon } := \chi _{\varepsilon } * \Psi . $ Then we first claim that   
\begin{equation} \label{eq:proof_FTL_1}
\langle \widehat{\varphi } , \Psi _{\varepsilon } \rangle \to 
\int _{\mathbb{R }^n } \, \widehat{L }_1 (\eta ) \psi (\eta ) \, d\eta ,   
\end{equation}   
To show this it suffices to consider the case that $\kappa \geq n $, since $\widehat{\varphi } $ is integrable if $\kappa < n $ and $\Psi _{\varepsilon } $ converges uniformly.   
If we assume for example that ${\rm supp }(\chi ) \subset B(0, 1 ) $ then by the Taylor expansion with remainder 
there exists a constant $C > 0 $ such that for all $\varepsilon \leq 1 $,   
\begin{eqnarray*}   
\left | \Psi _{\varepsilon } (\eta ) - \sum _{|\alpha | \leq \lfloor \kappa \rfloor - n } \frac{ \Psi _{\varepsilon } ^{(\alpha ) } (0) }{\alpha ! } \eta ^{\alpha } \right | &\leq & C \max _{|\beta | = \lfloor \kappa \rfloor - n + 1 } \sup _{B(0, 1 ) } | \Psi _{\varepsilon } ^{(\beta ) } | \cdot |\eta |^{\lfloor \kappa \rfloor - n + 1 }   
\\   
&\leq & C \max _{|\beta | = \lfloor \kappa \rfloor - n + 1 } \sup _{B(0, 2 ) } |\Psi ^{(\beta ) } | \cdot |\eta |^{\lfloor \kappa \rfloor - n + 1 } ,   
\end{eqnarray*}   
where we note that if $n \geq 2 $ or if $\kappa \notin \mathbb{N }^* $, then derivatives of $\Psi $ of order $\lfloor \kappa \rfloor - n + 1 $ exist, while if $n = 1 $ and $\kappa \in \mathbb{N }^* $, these derivatives exist a.e. but are uniformly bounded, and the estimate remains true.  Next, $\Psi ^{(\alpha ) } _{\varepsilon } (x) \to \Psi ^{(\alpha ) } (x) $ for $|\alpha | \leq \lceil \kappa \rceil - 1 $ since $\Psi $ is $ C^{\lceil \kappa \rceil - 1 } . $ Furthermore, since $(\sum _k \widehat{\varphi } (\eta + 2 \pi k ) )^{-1 } $ vanishes of order $\kappa $ in 0,  it follows that $\Psi ^{(\alpha ) } (0) = 0 $ for $|\alpha | \leq \lceil \kappa \rceil - 1 $ and (\ref{eq:proof_FTL_1}) follows by dominated convergence, 
\medskip   
   
Since $\Psi _{\varepsilon } $ is Schwartz-class, we have $\langle \widehat{\varphi } , \Psi _{\varepsilon } \rangle = \langle \varphi , \widehat{\Psi }_{\varepsilon } \rangle . $ By (\ref{eq:Fourier-series}) $\Psi = \left( \sum _k c_k e^{-i (k, \eta ) } \right) \psi (\eta ) $, which can be interpreted as the product of a tempered distribution and a test function, whose Fourier transform equals      
$$   
\widehat{\Psi } (x) = 
\sum _k c_k \widehat{\psi } (x - k ) .   
$$  
One easily verifies using (\ref{eq:est_c_k}) that 
$|\widehat{\Psi } (x) | \leq C (1 + |x| )^{- \kappa - n } . $   

\begin{versionA}   
\medskip   
   
\textcolor{blue}{For example, since $(1 + |x| ) \leq (1 + |k |) (1 + |x - k | ) $,   
$$   
(1  + |x| )^{\kappa + n } \sum _k | c_k \widehat{\psi } (x - k ) | \leq C_p \sum _k \frac{(1 + |x | )^{\kappa + n } }{(1 + |k | )^{\kappa + n } (1 + |x - k | )^p } \leq 
C_p \sum _k \frac{1 }{(1 + |x - k | )^{p - \kappa - n } } ,   
$$   
which converges if we take $p $ sufficiently big, and then is uniformly bounded in $x . $ }      
\medskip   

\end{versionA}   
Since $\widehat{\Psi }_{\varepsilon } (x) = \widehat{\chi } (\varepsilon x ) \widehat{\Psi } (x) $, and since $|\varphi (x) | \ \leq C (1 + |x | )^{\kappa - \rho } $ for some $\rho > 0 $, Lebesgue's dominated convergence theorem then shows that   
\begin{eqnarray*}   
\langle \varphi , \widehat{\Psi }_{\varepsilon } \rangle &=& \int _{\mathbb{R }^n } \, \varphi (x) \widehat{\Psi } (x) \widehat{\chi }(\varepsilon x ) \, dx \\   
&\to & \int _{\mathbb{R }^n } \varphi (x) \left( \sum _k c_k \widehat{\psi } (x - k ) \right) \, dx . 
\end{eqnarray*}   
Finally, one checks that the functions $(x, k ) \to c_k \varphi (x) \widehat{\psi } (x - k ) $ and $(x, k ) \to c_k \varphi (x + k ) \widehat{\psi } (x) $ are integrable on $\mathbb{R }^n \times \mathbb{Z }^n $ with respect to the product of the Lebesgue measure and the counting measure.   
\begin{versionA}   
\textcolor{blue}{For example, for any $p \in \mathbb{N } $,    
\begin{eqnarray*}   
| c_k \varphi (x) \widehat{\psi } (x - k ) | &\leq & C_p \frac{(1 + |x| )^{\kappa - \rho } }{(1 + |k | )^{n + \kappa } (1 + |x - k |)^p } \\   
&\leq & C_p \frac{(1 + |k| )^{\kappa - \rho } (1 + |x - k| )^{\kappa - \rho } }{(1 + |k | )^{n + \kappa } (1 + |x - k |)^p } \\   
&=& \frac{C_p }{(1 + |k | )^{n + \rho }  (1 + |x - k |)^{p - \kappa + \rho } } .   
\end{eqnarray*}   
}   
\end{versionA}   
A double application of Fubini's theorem then shows that the right hand side equals   
$$   
\int _{\mathbb{R }^n } \left( \sum _k c_k \varphi (x + k ) \right) \widehat{\psi }(x) \, dx ,  
$$   
which proves (\ref{eq:proof-FTL}). The pointwise identity (\ref{eq:formula_L_1}) follows since both sides are continuous.   
\smallskip   
   
\begin{versionA}   
   
\noindent \textcolor{blue}{{\bf Parenthetical remark}. Strictly speaking, we have shown that $L_1 $, defined as the inverse Fourier transform of the function $\widehat{\varphi } / \sum _k \widehat{\varphi }(\cdot + 2 \pi k ) $, is a.e. equal to the right hand side of (\ref{eq:formula_L_1}). Since both sides are continuous ($L_1 $ being the inverse Fourier transform of an integrable function, and the series $\sum _k c_k \varphi (x + k ) $ converging uniformly on compact subsets of $\mathbb{R }^n $) they are equal everywhere.}   
   
\end{versionA}   
\medskip   
   
\noindent {\it Proof of lemma \ref{lemma:growth_Lambda_varphi}}. The lemma presumably is classical, but since we could not locate a suitable reference (apart from the well-known case of homogeneous $\widehat{\varphi }  $), we sketch a proof for convenience of the reader. If $\kappa < n $, the inverse Fourier transform is a bounded function, so suppose that $\kappa \geq n . $ Since $\mathbf{1 }_{\{ |\eta | \geq 1 \} } \widehat{\varphi } (\eta ) $ is integrable, its inverse Fourier transform is a bounded continuous function, and it therefore suffices to examine the inverse Fourier transform of the tempered distribution defined by the first integral on the right hand side of (\ref{eq:widehat_varphi}). This distribution being of compact support, its inverse Fourier transform is the function $k(x) $ obtained by taking $\psi (\eta ) = (2\pi )^{-n } e^{i(x, \eta ) } $:   
\begin{equation} \label{eq:lemma_growth_Lambda_varphi}   
k(x) := (2 \pi )^{-n } \int _{|\eta | \leq 1 } \, \widehat{\varphi } (\eta ) \left( e^{i(x, \eta ) } - \sum _{j \leq \nu } \frac{i^j (x, \eta ) ^j }{j! } \right) \, d\eta ,   
\end{equation}   
where we put $\nu := \lfloor \kappa \rfloor - n . $ This can be bounded by 
\begin{eqnarray*}   
|k(x) | &\leq & C \int _{|\eta | \leq 1  } \, |\eta |^{- \kappa } \left | e^{i(x, \eta ) } - \sum _{j \leq \nu } \frac{i^j (x, \eta ) ^j }{j! } \right | \, d\eta \\   
&=& C |x| ^{\kappa - n } \int _{|\eta | \leq |x| } |\eta |^{- \kappa } \left | e^{i(\tfrac{x }{|x | } , \eta ) } - \sum _{j \leq \nu } \frac{i^j (\tfrac{x }{|x | } , \eta ) ^j }{j! } \right | \, d\eta .   
\end{eqnarray*}   
Split the integral into an integral over $|\eta | \leq c $ and one over the complement, where $c > 0 $ is some fixed number and where we assume wlog that $|x | > c . $ The first integral converges absolutely, since    
$$   
\left | e^{i(\tfrac{x }{|x | } , \eta ) } - \sum _{j \leq \nu }\frac{i^j (\tfrac{x }{|x | } , \eta ) ^j }{j! } \right | \leq \frac{1 }{\nu ! } | (\tfrac{x }{|x | } , \eta ) |^{\nu + 1 } \leq \frac{|\eta |^{\nu + 1 } }{\nu ! } ,   
$$   
and we can bound its contribution to $k(x) $ by $C |x |^{\kappa - n } . $ As for the 
integral over $| :eta \ > c $, it can be bounded by a constant times   
$$   
|x |^{\kappa - n } \sum _{j = 0 } ^{\nu } \int _c ^{|x | } r^{- \kappa + j + n - 1 } dr = |x |^{\kappa - n } \sum _{j = 0 } ^{\nu } \frac{1 }{j - \kappa + n } \left( |x|^{j - \kappa + n } - c^{j - \kappa + n } \right) ,   
$$   
assuming that $\kappa \notin \mathbb{N } . $ Since $j - \kappa + n \leq \nu - \kappa + n  \leq 0  $ by the definition of $\nu $, this will be bounded by $C |x|^{\kappa - n } . $ Finally, if  $\kappa \in \mathbb{N } $, $\kappa \geq n $, then $\nu = \kappa - n $ and 
\begin{eqnarray*}   
&&|x| ^{\kappa - n } \sum _{j = 0 } ^{\kappa - n } \int _c ^{|x | } r^{j - (\kappa - n ) - 1 } dr \\   
&= &   
|x|^{\kappa - n } \sum _{j = 0 } ^{\kappa - n - 1 } \frac{1 }{j - \kappa + n } \left( |x| ^{j - (\kappa - n ) } - c^{j - (\kappa - n ) } \right) + \log (|x | / c ) \\   
&\leq & C |x|^{\kappa - n } (\log |x | + 1 ) .   
\end{eqnarray*}   
\hfill $\Box $   
   
\begin{remark} \label{remark:symbol_estimates} \rm{The only hypotheses on $\widehat{\varphi }(\eta ) $ we needed for this lemma is that it be integrable on $\{ |\eta | \geq 1 \} $ and that $\widehat{\varphi } (\eta ) = O(|\eta |^{- \kappa } ) $ near 0. If we strengthen the first assumption to   
\begin{equation} \label{eq:remark_symbols_estimates}   
|\eta |^r |\widehat{\varphi } (\eta ) | \mathbf{1 }_{ \{ |\eta | \geq 1 \} } \in L^1 (\mathbb{R }^n ) ,   
\end{equation}   
where $r \in \mathbb{N } $, then $k $ will be $r $-times differentiable, and we will have that 
$$   
|\partial _x ^{\alpha } k(x) | \leq \left \{ \begin{array}{ll} C (|x |^{ \max ( \kappa - n - |\alpha | , 0 ) } + 1 ) &\kappa \notin \mathbb{N } \\   
C (|x |^{ \max ( \kappa - n - |\alpha | , 0 ) } \log |x | + 1 ) &\kappa \in \mathbb{N } ,   
\end{array}   
\right.   
$$   
for $|\alpha | \leq r $: 
it suffices to observe that if $k(x) $ is given by (\ref{eq:lemma_growth_Lambda_varphi}) then its derivative of order $\alpha $ is given by the same formula with $\widehat{\varphi } (\eta ) $ replaced by $(i \eta )^{\alpha } \widehat{\varphi }(\eta ) . $   
   
These estimates can be used to obtain symbol estimates for the generator of pure-jump L\'evy processes 
with L\' evy measure 
$$   
d\nu (\eta ) = \frac{h(\eta ) }{|\eta |^q } d\eta ,   
$$   
with $q < n + 2 $, and $h(\eta ) $ a rapidly decreasing continuous function satisfying (\ref{eq:remark_symbols_estimates}) for all $r . $ The inverse Fourier transform of 
$\Lambda _{ |\eta |^{-q } h } $ then is, modulo a function in $C^{\infty } _b $, equal to the symbol of the generator of the L\'evy process, and the estimates show this symbol to be in $S^{\max (q - n , 0 ) } _0 $ if $q \notin \mathbb{N } $, and in $S_0 ^{\max (q - n , 0 ) + \varepsilon } $ for any $\varepsilon > 0 $ otherwise (even a bit better, since the first few derivatives will decay relative to the symbol itself). Examples are given by the CGMY-processes 
which are used in financial modeling.   
\medskip   
   
\begin{versionA}   
   
\noindent \textcolor{blue}{{\bf Question}: under what conditions on $h $, that is, on the L\'evy measure, will the symbol of the generator be elliptic? For example, under what conditions on $\widehat{\varphi }(\eta ) $ will $k(x) $ above satisfy $|k(x) | \geq c |x |^{\kappa - n } $ when $\kappa \notin \mathbb{N } . $ Concretely, as an example, is the generator of the CGMY-process elliptic if $Y \in (0, 2 ) \setminus \{ 1 \} ? $   
As a further question, motivated by the CGMY with $Y = 0, 1 $, consider symbol spaces with arbitrary weight $w(\xi ) $, instead of $(1 + |\xi |)^p ? $   
}   
\end{versionA}   
}   
\end{remark}   
   
\section{\bf Some technical proofs}  
\label{Appendix_proof_thm_conv_RBF_bis}   
   
\subsection{Proof of lemma \ref{lemma:convergence_RBF_bis}} Let $F  \in L^1 \left( \mathbb{R }^n \setminus 0 , (|\xi |^{\kappa } \wedge 1 ) d\xi \right) $, where $a \wedge b := \min (a, b ) $ and $\kappa \geq 0 . $ Then $F $ gives rise to a tempered distribution $\Lambda _F \in \mathcal{S }' (\mathbb{R }^n ) $ defined as follows: if $g \in C^{\infty } _c (\mathbb{R }^n ) $ be equal to 1 on a neighbourhood of 0, we put       
\begin{eqnarray} \label{eq:u_f_hat}   
\langle \Lambda _F, \psi \rangle &:= &\int _{\mathbb{R }^n } \left( \psi (\xi ) - \sum _{|\alpha | \leq \lceil \kappa \rceil - 1 } \psi ^{(\alpha ) } (0) \frac{\xi ^{\alpha } }{\alpha ! } \right) g(\xi ) F (\xi ) \, d\xi \\   
&& + \int _{\mathbb{R }^n } (1 - g(\xi ) ) F(\xi ) \psi (\xi ) \, d\xi , \ \ \psi \in \mathcal{S }(\mathbb{R }^n ) . \nonumber   
\end{eqnarray}
The integral converges since $\psi - \sum _{|\alpha | \leq \lceil \kappa \rceil - 1 } \psi ^{(\alpha ) } (0) \xi  ^{\alpha } / \alpha ! = O(|\xi |^{\lceil \kappa \rceil } ) = O(|\xi |^{\kappa } ) $ in a neighbourhood of 0 and defines a distribution of order $\lceil \kappa \rceil - 1 .  $ Note that $\Lambda _F $ coincides on $\mathbb{R }^n \setminus 0 $ with the locally integrable function $F . $ 

      
We next observe that $\Lambda _F $ extends to a continuous linear functional on the H\"older space $C_b ^{\lceil \kappa \rceil - 1 , \lambda }  := C_b ^{\lceil \kappa \rceil - 1 , \lambda }  (\mathbb{R }^n ) $ with $\lambda = \kappa - (\lceil \kappa \rceil - 1 ) . $ Indeed, if $\psi \in C^{K , \lambda } (\mathbb{R }^n ) $, 
then the Taylor expansion formula with integral remainder term easily implies that   
\begin{equation} \label{eq:Taylor_remainder}   
\left | \psi (\xi ) - \sum _{ |\alpha | \leq K } \psi ^{(\alpha ) } (0) \xi ^{\alpha } / \alpha ! 
\right | \leq C \left( \sum _{ |\beta | = K } || \psi ^{(\beta ) } ||_{0 , \lambda } \right) |\xi |^{K + \lambda } ,   
\end{equation}   
which shows, with $K = \lceil \kappa \rceil - 1 $ and $\lambda = \kappa - (\lceil \kappa \rceil - 1 ) $, that $\langle \Lambda _F , \psi \rangle $ is well-defined and continuous.   
\medskip   
   
\begin{versionA}   
   
\noindent \textcolor{blue}{We check (\ref{eq:Taylor_remainder}) (not to be retained in the final version): if $n = 1 $, this follows from the Taylor formula with integral remainder applied at order $K - 1 $: since    
\begin{eqnarray*}   
\psi (\xi ) - \sum _{\alpha \leq K - 1 } \frac{\psi ^{(\alpha ) } (0) }{\alpha ! } \xi ^{\alpha } &=& \int _0 ^{\xi } \psi ^{(K ) } (\eta ) \frac{\ \ (\xi - \eta )^{K - 1 } }{(K - 1 )! } d\eta \\   
&=& \frac{\psi ^{(K ) }(0) }{K! } \xi ^K + \int _0 ^{\xi } (\psi ^{(K ) }(\eta ) - \psi ^{(K ) }(0) ) \frac{(\xi - \eta )^{K - 1 } }{(K - 1 )! } d\eta ,   
\end{eqnarray*}    
the remainder can be bounded by   
$$   
|| \psi ^{(K ) } ||_{0, \lambda } \int _0 ^{\xi } |\eta |^{\lambda } \frac{(\xi - \eta )^{K - 1 } }{(K - 1 )! } d\eta \leq | \psi ^{(K ) } ||_{ 0, \lambda } \frac{|\xi |^{K + \lambda } }{K! } .   
$$   
On $\mathbb{R }^n $ the estimate follows by applying the one-dimensional estimate to the function restricted to rays through the origin.   
}  
\medskip   
   
\end{versionA}   
   
We can, in particular, let $\Lambda _F $ act on the imaginary exponentials $\xi \to e^{i (x, \xi ) } . $ The function   
$$   
\check{F } : x \to (2 \pi )^{-n } \left \langle \Lambda _F , e^{i (x, \xi ) } \right \rangle .   
$$   
is then found to be bounded by $ C (1 + |x | )^{\kappa } $, since $|| e^{i(x, \xi ) } ||_{K , \lambda } \leq C (1 + |x |^{K + \lambda } ) $, 
\begin{versionA}   
\textcolor{blue}{as follows from   
$$   
\sup _{\xi } \frac{| e^{i (x, \xi ) } - 1 | }{ |\xi |^{\lambda } } = |x |^{\lambda } \sup _{\xi } \frac{| e^{i (x, \xi ) } - 1 | }{ (|x | |\xi |)^{\lambda } } \leq |x |^{\lambda } \sup _{\xi } \frac{\min ( |x | |\xi | , 2 ) }{ (|x | |\xi |)^{\lambda } } ,   
$$   
}   
\end{versionA}   
and one easily verifies that the inverse Fourier transform of $\Lambda _F $ coincides with $\check{F } . $ 
If $\kappa \in \mathbb{N } $ one has the stronger estimate   
\begin{equation}\label{eq:decay_F_check}   
|\check{F } (x) | = o(|x |^{\kappa } ) , \ \ |x | \to \infty ,   
\end{equation}    
which can be seen as follows: write $F= \chi F  + (1 - \chi ) F $ with $\chi $ the characteristic function of a small ball around 0. Since 
$(1 - \chi ) F $ is integrable, its inverse Fourier transform tends to 0 at infinity, by the Riemann-Lebesgue lemma. We can 
therefore wlog assume that $F $ is supported in $\{ g = 1 \} . $ If we apply (\ref{eq:u_f_hat}) with $\psi (\xi ) = e^{i (x, \xi ) } $ then\footnote{e.g. by using the Taylor formula with integral remainder in the form   
$$   
\psi (\xi ) - \sum _{|\alpha | \leq \kappa - 1 } \psi ^{(\alpha ) } (0) \frac{\xi ^{\alpha } }{\alpha ! } = \int _0 ^1 \frac{((1 - s )^{\kappa - 1 } }{(\kappa - 1 ) ! } \frac{d ^{\kappa } }{d s^{\kappa } } \psi _{\xi } (s) ds ,   
$$   
where $\psi _{\xi } (s) := \psi (s \xi ) $ }   
\begin{eqnarray*}   
\check{F } (x) &=& \sum _{|\alpha |= \kappa } \int _{\mathbb{R }^n } F(\xi ) \frac{(ix )^{\alpha } \xi ^{\alpha } }{\alpha ! } \left( \int _0 ^1 \frac{ \ \ (1 - s )^{\kappa - 1 } }{(\kappa - 1 ) ! } e^{i s (x, \xi ) } ds \right) \frac{d\xi }{(2 \pi )^n } \\   
&=: & \sum _{|\alpha | = \kappa } (i x )^{\alpha } \int _0 ^1 \check{F }_{\alpha } (sx ) \frac{ \ \ (1 - s )^{\kappa - 1 } }{(\kappa - 1 ) ! } ds ,   
\end{eqnarray*}   
where $\check{F }_{\alpha } (x) $ is the inverse Fourier transform of the $L^1 $-function $\xi \to \xi ^{\alpha } F(\xi ) . $ By the Riemann-Lebesgue lemma, $\check{F }_{\alpha } (sx) \to 0 $ as $x \to \infty $, for all $s \in (0, 1 ] $, and the same is true for the integral over $s \in [0, 1 ] $, 
by the dominated convergence theorem (the $\check{F }_{\alpha } $ are bounded). Hence $\check{F }(x) / |x |^{\kappa } \to 0 $ for $x \to \infty $, as claimed.   
\medskip   
   
Now let $f $ be a measurable function on $\mathbb{R }^n $ of polynomial growth of order strictly less than $\kappa $, such that its Fourier transform $\widehat{f } $ (in the sense of tempered distributions) satisfies   
$$   
\widehat{f } \big{| }_{\mathbb{R }^n \setminus 0 } \in L^1 \left( \mathbb{R }^n , (|\xi |^{\kappa } \wedge 1 ) d\xi \right).   
$$   
We write $\Lambda _{\widehat{f } } $ for $\Lambda _{\widehat{f } |_{\mathbb{R }^n \setminus 0 } } . $ Then $\widehat{f } - \Lambda _{\widehat{f } } $ is a distribution which is supported in 0, and therefore of the form $\sum _{|\alpha | \leq N } c_{\alpha } \delta _0 ^{(\alpha ) } $ for certain $N \in \mathbb{N } $ and $c_{\alpha } \in \mathbb{C } $ with $\sum _{|\alpha | = N } |c_{\alpha } | \neq 0 . $    Since the inverse Fourier transform of $\widehat{f } - \Lambda _{\widehat{f } } $ is a polynomial of degree $N $, it follows that $N \leq \lceil \kappa \rceil - 1 $, the largest integer which is strictly smaller than $\kappa $, since otherwise $|f(x) | $ would grow at a rate of at least $|x |^{\lceil \kappa \rceil } $ 
in certain directions. If $\kappa \notin \mathbb{N } $  this would contradict    
the bound $\check{F }(x) = 0(|x |^{\kappa } ) $, and if $\kappa \in \mathbb{N } $ this would contradict (\ref{eq:decay_F_check}).   

In follows that $\widehat{f } =  \Lambda _{\widehat{f } } + \sum _{|\alpha | \leq N } c_{\alpha } \delta _0 ^{(\alpha ) } $ also extends to a continuous linear functional on $C^{\lceil \kappa \rceil - 1 , \kappa - (\lceil \kappa \rceil - 1 ) } . $ We exploit this to define $\Sigma _h (\widehat{f } ) $ by duality.   
\medskip   
   
If $\psi \in \mathcal{S }(\mathbb{R }^n ) $, we let      
\begin{equation}   
\Sigma _h ' (\psi ) := \sum _k \psi (\xi + 2 \pi h^{-1 } k ) \, \widehat{L }_1 (h \xi + 2 \pi k ) .   
\end{equation}   
Note that $\Sigma _h ' $ is the formal (real) adjoint of $\Sigma _h . $ By lemma \ref{lemma:der_L_hat}, $\Sigma _h ' (\psi ) $ 
is $C_b ^{\lceil \kappa \rceil - 1 , \lambda }  $ with $\lambda = \kappa - (\lceil \kappa \rceil - 1 ) $ and uniformly bounded together with all its derivatives, since $2 \pi h^{-1 } $-periodic. In fact, this is true even if $\psi \in C_b ^{\lceil \kappa \rceil - 1 , \lambda }  $ with the same $\lambda $, on account of the decay at infinity of $\widehat{L }_1 . $ We can then define $\Sigma _h (\widehat{f } ) $, as a tempered distribution and, more generally, as a bounded linear functional on $C_b ^{\lceil \kappa \rceil - 1 , \lambda }  (\mathbb{R }^n ) $ by      
\begin{equation}   
\left \langle \Sigma _h (\widehat{f } ) , \psi \right \rangle := \left \langle \widehat{f } , \Sigma _h ' (\psi ) \right \rangle .   
\end{equation}   
We next check that $\Sigma _h (\widehat{f } ) $ is the Fourier transform, in distribution sense, of $s _h [f ] . $ This is done by a standard approximation argument, with some care with the spaces in which the approximating sequence converges. We first note that we can assume without loss of generality that $\widehat{f } $ is compactly supported: indeed, we can write $f = f_1 + f_2 $ with $\widehat{f }_1 $ compactly supported and $\widehat{f }_2 \in L^1 (\mathbb{R }^n ) $, and we know already that $\widehat{s _h [f_2 ] } = \Sigma _h (\widehat{f }_2 ) . $   
   
So let $\widehat{f } $ be compactly supported, and let $\chi \in C^{\infty } _c (\mathbb{R }^n ) $ be a non-negative symmetric function with $\int _{\mathbb{R }^n } \chi d\eta = 1 . $ Let $\chi _{\varepsilon } (\eta ) := \varepsilon ^{-n } \chi (\eta / \varepsilon ) . $ Then $\widehat{f }* \chi _{\varepsilon } \in C^{\infty } _c (\mathbb{R }^n ) . $   
   
\begin{lemma} $\widehat{f } * \chi _{\varepsilon } \to \widehat{f } $ in the dual of $C^{K , \lambda } $, with $K = \lceil \kappa \rceil - 1 $ and $\lambda = \kappa - K . $   
\end{lemma}   
   
\noindent {\it Proof}. On account of the symmetry of $\chi $,      
$$   
\langle \widehat{f } * \chi _{\varepsilon } , \psi \rangle = \langle \widehat{f } , \psi * \chi _{\varepsilon } \rangle ,   
$$   
which is valid both for Schwarz-class functions $\psi \in \mathcal{S } $ and for $\psi \in C^{K , \lambda } . $ Write $\psi _{\varepsilon } := \psi * \chi _{\varepsilon } . $ If $\psi \in C^{K , \lambda } $, then $\psi _{\varepsilon } ^{(\alpha ) } (x) \to \psi ^{(\alpha ) } (x) $ pointwise on $\mathbb{R }^n $ for all $|\alpha | \leq K $, while a trivial estimate shows that $|| \psi _{\varepsilon } ^{(\alpha ) } ||_{0, \lambda } \leq || \psi ^{(\alpha ) } ||_{0 , \lambda } $, uniformly in $\varepsilon > 0 $, for $|\alpha | = K . $ This, together with the remainder estimate (\ref{eq:Taylor_remainder}), the integrability of $\widehat{f } (\xi ) (|\xi |^{\kappa } \wedge 1 ) $ and Lebesgue's dominated convergence theorem, implies that   
$ \langle \Lambda _{\widehat{f } } , \psi _{\varepsilon } \rangle \to \langle \Lambda _{\widehat{f } } , \psi \rangle . $ Since also $\langle \delta _0 ^{(\alpha ) } , \psi _{\varepsilon } \rangle \to \langle \delta ^{(\alpha ) } _0 , \psi \rangle $ for all $|\alpha | \leq K $, the lemma follows.   
   
\hfill $\Box $   
\medskip   
   
\noindent The lemma immediately implies that if $\psi \in \mathcal{S } (\mathbb{R }^n ) $, then $\langle \widehat{f } * \chi _{\varepsilon } , \Sigma _h ' (\psi ) \rangle \to \langle \widehat{f } , \Sigma _h ' (\psi ) \rangle $, so $\Sigma _h ( \widehat{f } * \chi _{\varepsilon } ) \to \Sigma _h (\widehat{f } ) $ in $\mathcal{S }' (\mathbb{R }^n ) $ and even in $(C^{K , \lambda } ) ' $ with $K $ and $\lambda $ as above.   
   
On the other hand, if we let $f_{\varepsilon } $ be the inverse Fourier transform of $\widehat{f } * \chi _{\varepsilon } $, then $f_{\varepsilon } \in \mathcal{S } (\mathbb{R }^n ) $ since $\widehat{f } * \chi _{\varepsilon } $ is, and $\widehat{s_h [ f_{\varepsilon } ] } = \Sigma _h ( \widehat{f } * \chi _{\varepsilon } ) . $ We have that $f_{\varepsilon } (x) = (2 \pi )^{n } f(x) \check{\chi } (\varepsilon x ) $, with $\check{\chi } $ the inverse Fourier transform of $\chi $, so $\check{\chi } \in \mathcal{S }(\mathbb{R }^n ) $ and $(2 \pi )^n \check{\chi } (0) = 1 . $  By hypotheses, $f \in L^{\infty } _{- p } $ for some $p < \kappa . $ If $a > 0 $ such that $p + a < \kappa $, then by (\ref{eq:s_h_assertion_1}), writing $\widetilde{\chi } := (2 \pi )^n \check{\chi } $,   
\begin{eqnarray*}   
|| s_h [f_{\varepsilon } ] - s_h [f ] ||_{\infty , -(p + a ) } &\leq & C || f \left( \widetilde{\chi } (\varepsilon \cdot ) - 1 \right) ||_{\infty , -(p + a ) } \\   
&\leq & C || f ||_{\infty , -p } \, \sup _{x \in \mathbb{R }^n } \frac{ |\widetilde{\chi } (\varepsilon x ) - 1 | }{(1 + |x | )^a } \to 0 ,   
\end{eqnarray*}   
as $\varepsilon \to 0 $, using for example the first order Taylor expansion for $\widetilde{\chi } $ for $|x | \leq \varepsilon ^{- 1 / 2 } $ plus a trivial estimate for $|x | > \varepsilon ^{-1 / 2 } . $ This certainly implies that $s_h [ f_{\varepsilon } ] \to s_h [f ] $ in $\mathcal{S }' (\mathbb{R }^n ) $, so we conclude that $\widehat{s_h [f_{\varepsilon } ] } =   
\Sigma _h ( \widehat{f } * \chi _{\varepsilon } ) \to  \widehat{s_h [f ] } $ and therefore $\widehat{s_h [f ] } = \Sigma _h (\widehat{f } ) $ as distributions.   
\medskip   
   
We finally show that $\Sigma _h (\widehat{f } ) = \widehat{f } + F $, where $F $ is the (distribution obtained by integrating against the) function   
\begin{equation}   
F (\xi ) = \widehat{f }(\xi ) (\widehat{L }_1 (h \xi ) - 1 ) + \sum _{k \neq 0 } \widehat{f }(\xi + 2 \pi h^{-1 } k ) \widehat{L }_1 (h \xi ) .   
\end{equation}   
We first check that $F $ is well-defined and in $L^1 $: first of all, each of the terms on the right hand side is in $L^1 $, on account of the Fix-Strang conditions satisfied by for $\widehat{L }_1 $ and the integrability of $(|\xi |^{\kappa } \wedge 1 ) \widehat{f }(\xi ) . $ Next, the function $(\xi , k ) \to \widehat{f }(\xi + 2 \pi h^{-1 } k ) (\widehat{L }_1 (h \xi ) - \delta _{0k } ) $ is absolutely integrable on $\mathbb{R }^n \times \mathbb{Z }^n $ with respect to the product of Lebesgue measure and the counting measure, since   
\begin{eqnarray*}   
&& \sum _k \int _{\mathbb{R }^n } | \widehat{f }(\xi + 2 \pi h^{-1 } k ) | \, |(\widehat{L }_1 (h \xi ) - \delta _{0k } ) | \, d\xi \\   
&=& \int _{\mathbb{R }^n } (1 - \widehat{L }_1 (h \xi ) ) |\widehat{f }(\xi ) | d\xi + \sum _{k \neq 0 } \widehat{L }_1 (h \xi + 2 \pi k ) |\widehat{f }(\xi ) | d\xi \\   
&=& 2 \int _{\mathbb{R }^n } (1 - \widehat{L }_1 (h \xi ) ) |\widehat{f } (\xi ) | d\xi .   
\end{eqnarray*}   
Fubini's theorem then implies that $F(\xi ) $ is well-defined for almost all $\xi \in \mathbb{R }^n $ and that $F \in L^1 (\mathbb{R }^n ) . $  If $\psi \in \mathcal{S } (\mathbb{R }^n ) $, then a double application of Fubini 
will show that   
\begin{eqnarray*}   
&& \int _{\mathbb{R }^n } F(\xi ) \psi (\xi ) d\xi \\   
&& = \int _{\mathbb{R }^n } \left( \psi (\xi ) (\widehat{L }_1 (h \xi ) - 1 ) + \sum _{k \neq 0 } \psi (\xi + 2 \pi h^{-1 } k ) \widehat{L }_1 (h \xi + 2 \pi k ) \right) \widehat{f }(\xi ) \, d\xi \\   
&=& \int _{\mathbb{R }^n } (\Sigma _h ' (\psi ) - \psi ) \widehat{f }(\xi ) d\xi .   
\end{eqnarray*}   
Since, by the Fix-Strang conditions (\ref{eq:SFC}),  
all derivatives of order $\leq \lceil \kappa \rceil - 1 $ of $\Sigma _h (\psi ) - \psi $ in 0 are 0, the last integral is equal to $\langle \widehat{f } , \Sigma _h ' (\psi ) - \psi \rangle = \langle \Sigma _h (\widehat{f } )- \widehat{f } , \psi \rangle $,  and therefore $\Sigma _h (\widehat{f } ) - \widehat{f } = F $, which finishes the proof of lemma \ref{lemma:convergence_RBF_bis}.   
   
\begin{remark} \rm{The lemma and its proof generalizes to $f $'s such that $\widehat{f } |_{\mathbb{R }^n \setminus 0 } $ is a finite Borel measure with respect to which the function $|\xi |^{\kappa } \wedge 1 $ is integrable, provided that $\kappa \notin \mathbb{N } $ (the reason being that we then no longer have (\ref{eq:decay_F_check})).   
\begin{versionA}   
\medskip   
   
\textcolor{blue}{We briefly verify the different steps of the proof (not to be retained for the final version, but for own reassurance only):   
\begin{enumerate}     
\item [1. ] If $\nu $ is a locally finite Borel measure on $\mathbb{R }^n \setminus 0 $ such that $\int ( |\xi |^{\kappa } , 1 ) d|\nu | (\xi ) < \infty $, then (\ref{eq:u_f_hat}), with $F(\xi ) d\xi $ replaced by $d\nu (\xi ) $ and integrals taken over $\mathbb{R } ^n \setminus 0 $ instead of $\mathbb{R }^n $, still extends $\nu $ to a tempered distribution $\Lambda _{\nu } $ on $\mathbb{R }^n $, and to an element of the dual of the appropriate H\"older space.   
\item [2. ] The inverse Fourier transform of $\Lambda _{\nu } $ is then a function which is $O(|x |^{\kappa } ) $ at infinity.      
\item [3. ] It follows that if $|f(x) | = O(|x |^{\kappa } ) $ at infinity, and if $\widehat{f } |_{\mathbb{R }^n \setminus 0 } $ is a measure $\nu _f $ as 1 above, then if $\kappa \notin \mathbb{N } $, $\widehat{f } = \Lambda _{\nu _f } + \sum _{|\alpha | \leq \lfloor \kappa \rfloor } c_{\alpha } \delta _0 ^{(\alpha ) } $, and $\widehat{f } $ extends to a continuous linear functional on $C_b ^{\lceil \kappa \rceil - 1 , \lambda }  $ with $\lambda = \kappa - (\lceil \kappa \rceil - 1 ) . $ This allows to define $\Sigma _h (\widehat{f } ) $ by duality.   
\item [4. ] If $|f(x) | = O(|x |^p ) $ at infinity, then $s_h [f ] $, as before, is well-defined and of polynomial growth, and $\widehat{s_h [f ] } = \Sigma _h [\widehat{f } ] . $   
\item [5. ] If $\psi \in C_b ^{\lceil \kappa \rceil - 1 , \lambda }  $, then since $\widehat{f } = \Lambda _{\nu _f } + \sum _{|\alpha | \leq \lfloor \kappa \rfloor } c_{\alpha } \delta _0 ^{(\alpha ) } $ and since all derivatives of $\Sigma _h ' (\psi ) - \psi $ in 0 are 0 (by the Fix-Strang conditions satisfied by $\widehat{L }_1 $),   
\begin{eqnarray*}   
\langle \Sigma _h [\widehat{f } ] - \widehat{f } , \psi \rangle &=& \langle \widehat{f } , \Sigma _h ' (\psi ) - \psi \rangle \\   
&=& \int _{\mathbb{R }^n \setminus 0 } \left( \Sigma _h ' (\psi ) - \psi \right) d\nu _f (\xi ) .   
\end{eqnarray*}   
\item [6. ] Finally, with $\nu = \nu _f $,   
\begin{eqnarray*}   
&&\left | \int _{\mathbb{R }^n \setminus 0 } \left( \Sigma _h ' (\psi ) - \psi \right) d\nu _f (\xi ) \right | \\   
&&\leq \left | \int _{\mathbb{R }^n \setminus 0 } (\widehat{L }_1 (h \xi )  - 1 ) \psi (\xi ) d \nu (\xi ) \right | + \sum _{k \neq 0 } \left | \int _{\mathbb{R }^n \setminus 0 } \widehat{L }_1 (h\xi + 2 \pi k ) \psi (\xi + 2 \pi h^{-1 } k ) \, d \nu (\xi ) \right | 
\\   
&\leq & || \psi ||_{\infty } \left( \int _{\mathbb{R }^n \setminus 0 } (1 - \widehat{L }_1 (h \xi ) ) d |\nu | (\xi ) + \sum _{k \neq 0 } \widehat{L }_1 (h\xi + 2 \pi k ) \, d |\nu | (\xi ) \right) \\   
&=& 2 || \psi ||_{\infty } \int _{\mathbb{R }^n \setminus 0 } (1 - \widehat{L }_1 (h \xi ) ) \, d |\nu |(\xi ) ,   
\end{eqnarray*}   
which implies that $\Sigma _h (\nu ) - \nu $ is a finite measure whose variation norm (= norm as an element of the dual of $C_b (\mathbb{R }^n $, in which $C_b ^{\lceil \kappa \rceil - 1 , \lambda } $  is dense) is bounded 
$2 \int _{\mathbb{R }^n } (1 - \widehat{L }_1 (h \xi ) ) \, d |\nu |(\xi ) . $ If $|\xi |^{\kappa } \in L^1 (\mathbb{R }^n \setminus 0 , d|\nu | ) $ this, together with the Fix-Strang condition in 0, then implies that 
$$   
|| s_h [f ] - f ||_{\infty } \leq C h^{\kappa } \int _{\mathbb{R } ^n \setminus 0 } |\xi |^{\kappa } d|\nu | (\xi ) .   
$$   
\end{enumerate}   
}   
\end{versionA}   
}   
\end{remark}   
   
\begin{versionA}   
\subsection{An extension of lemma \ref{lemma:convergence_RBF_bis}} \textcolor{blue}{This subsection contains an alternative, soft-analytic, argument for the lemma - written just for own reassurance and not to be retained for final version (RB)}   
   
Part of the argument above is soft analytic, and can be used to show that the relation $\widehat{\sigma _h [f ] } = \Sigma _h (\widehat{f } ) $ generalises to $f \in L^{\infty } _{-p } (\mathbb{R }^n ) $ for some $p < \kappa $ such that   
\begin{enumerate}   
   
\item [1. ] its Fourier transform $\widehat{f } $ extends to a continuous linear functional on $C_b ^{\lceil \kappa \rceil - 1 , \lambda } (\mathbb{R }^n ) $, where $\lambda = \kappa + 1 - \lceil \kappa \rceil $, such that   
   
\item [2. ] for any uniformly bounded sequence $g_{\nu } \in C_b ^{\lceil \kappa \rceil - 1 , \lambda }  (\mathbb{R }^n ) $ such that $\partial ^{\alpha } _{\xi } g_{\nu } $ converges pointwise to $\partial ^{\alpha } _{\xi } g  $ for some $g \in C_b ^{\lceil \kappa \rceil - 1 , \lambda } (\mathbb{R }^n ) $, we have that $\langle \widehat{f } , g _{\nu } \rangle \to \langle \widehat{f } , g \rangle $;   
   
\end{enumerate}   
\noindent The second property expresses continuity of $\widehat{f } $ with respect to a weaker topology than the norm topology. A continuous linear functional $u $ on $C^{K , \lambda } $ satisfying property 2 will be completely determined by its restriction to $\mathcal{S }(\mathbb{R }^n ) $ or $C^{\infty } _c (\mathbb{R }^n ) $, which is not the case for a general element of the dual (which include Banach-limit type linear functionals).      
 
\medskip   
   
To prove that $\widehat{\sigma _h [f ] } = \Sigma _h (\widehat{f } ) $ we approximate such $\widehat{f } $ in (the weak-* topology on) $(C_b ^{\lceil \kappa \rceil - 1 , \lambda } ) ' $ by  a sequence of Schwarz-class functions $\widehat{f }_{\nu } \in \mathcal{S }(\mathbb{R }^n ) . $ This can be done in two steps:   
\begin{itemize}   
\item If $\chi _1 \in C^{\infty } _c (\mathbb{R }^n ) $ such that $\chi _1 (\xi ) = 1 $ in a neighbourhood of 0, then one checks, using property 2 above, that   
$\chi _1 (\cdot / \nu ) \widehat{f } \to \widehat{f } $ in $\left ( C_b ^{\lceil \kappa \rceil - 1 , \lambda }  \right) ' $: in fact, if $g \in C_b ^{\lceil \kappa \rceil - 1 , \lambda }  $, then $g_{\nu } (\xi ) := \chi _1 (\xi / \nu ) g(\xi ) $ is uniformly bounded (it may be easiest to observe that $C_b ^{\lceil \kappa \rceil - 1 , \lambda }  $ is a    normed algebra with its natural norm, and that $\chi (\cdot / \nu ) $ is uniformly bounded in norm; recall that by definition all derivatives of order $\leq \lceil \kappa \rceil - 1 $ are uniformly bounded on $\mathbb{R }^n . $)   
   
\item We can therefore suppose that $\widehat{f } $ has compact support. If $\chi _2 \in C^{\infty } _c (\mathbb{R }^n ) $ with $\chi _2 $ symmetric and $\int _{\mathbb{R }^n } \chi _2 d\xi = 1 $, we let $\chi _{2, \nu } (\xi ) := \nu ^n \chi _2 (\nu \xi ) $ and $F_{\nu } := \widehat{f } * \chi _{2 , \nu } \in    \mathcal{S } (\mathbb{R }^n ) . $ We claim that $F_{\nu } \to \widehat{f } $ in $\left( C_b ^{\lceil \kappa \rceil - 1 , \lambda }  \right) ' $: if $g \in C_b ^{\lceil \kappa \rceil - 1 , \lambda }   $, then $\langle F_{\nu } , g \rangle = \langle \widehat{f } , \chi _{2, \nu } * g \rangle $, and   
$$   
g_{\nu } (\xi ) := \chi _{2, \nu } * g (\xi ) = \int _{\mathbb{R }^n } \chi _2 (\eta ) g(\xi - \eta / \nu ) d\eta ,   
$$   
so $\partial ^{\alpha } _{\xi } g_{\nu } (\xi ) \to \partial ^{\alpha } _{\xi } g $ for $|\alpha | \leq \lceil \kappa \rceil - 1 $, while one easily checks that $g_{\nu } $ is uniformly norm-bounded in $C_b ^{\lceil \kappa \rceil - 1 , \lambda }  $: see the lemma below. Hence, by the hypothesis on $\widehat{f } $, $\langle \widehat{f } , g_{\nu } \rangle \to \langle \widehat{f } , g \rangle $, and therefore $F_{\nu } \to \widehat{f } . $ Now let $f_{\nu } $ be the inverse Fourier transform of $F_{\nu } . $   
   
\end{itemize}   
   
\noindent N.B. We used the following trivial observation:   
   
\begin{lemma} If $h $ is (uniformly) H\"older continuous of exponent $\lambda $, then $ || h * \chi _{\nu } || _{0, \lambda } $ is uniformly bounded in $\nu $, for any integrable function $\chi . $   
\end{lemma}      
   
\noindent {\it Proof}. If we let $\varepsilon = \nu ^{-1 } $, then   
\begin{eqnarray*}   
| h * \chi _{\nu } (x) - h * \chi _{\nu } (y) | &\leq & \int _{\mathbb{R }^n } \left | h (x - \varepsilon z ) - h (y - \varepsilon z ) \right | \, |g(z) | \, dz \\   
&\leq & ||g ||_1 \, ||h ||_{0, \lambda } .   
\end{eqnarray*}   
\hfill $\Box $   
\medskip   
   
Returning to the situation of lemma \ref{lemma:convergence_RBF_bis}, if we take $g = \Sigma _h ' (\psi ) \in C_b ^{\lceil \kappa \rceil - 1 , \lambda }  $ with $\psi \in \mathcal{S } (\mathbb{R }^n ) $, we find as before that 
$$   
\langle \Sigma _h (\widehat{f }_{\nu } ) , \psi \rangle = \langle \widehat{f }_{\nu } , \Sigma _h ' (\psi ) \rangle \to \langle \widehat{f } , \Sigma _h ' (\psi ) \rangle = \langle \Sigma _h (\widehat{f } ) , \psi \rangle ,   
$$   
so that   
$$   
\Sigma _h (\widehat{f }_{\nu } ) \to \Sigma _h (\widehat{f } ) \ \mbox{in } \mathcal{S }' (\mathbb{R }^n ) .   
$$   
To complete the proof that $\widehat{s _h [f ] } = \Sigma _h (\widehat{f } ) $, we must check that for the approximating sequence of Schwarz-class functions $\widehat{f }_{\nu } $ we constructed above, we also have that $s_h [f_{\nu } ] \to s_h [f ] $ in $\mathcal{S }' (\mathbb{R }^n ) . $ We will in fact verify that for any $0 < a < \kappa - p $,       
$$   
|| s_[ [f_{\nu } ] - s_h [f ] ||_{\infty , -p - a } \to 0 .   
$$   
The operations of convolution with $\chi _{2, \nu } $ and multiplication with $\chi _1 (\cdot / \nu ) $ on the Fourier-side correspond to, respectively, multiplying respectively convolving with a function of the form $\chi (x / \nu ) $ respectively $\chi _{\nu } (x) := \nu ^n \chi (\nu x ) $ before Fourier transform, with $\chi (0) = 1 $ respectively $\int \chi dx = 1 $ \rouge (check that this is indeed OK!) \noir We already checked that if $f \in L^{\infty } _{-p } (\mathbb{R }^n ) $, then   
$$   
s_h \left[ f \chi (\cdot / \nu ) \right ] \to s_h [f ] \ \mbox{in } L^{\infty } _{- p - a } (\mathbb{R }^n ) .   
$$   
The same is true when we convolve with an approximation of the identity $\chi _{\nu } $: we first observe that $|| f * \chi _{\nu } ||_{\infty , -p } $ is uniformly bounded in $\nu $: this follows from   
\begin{eqnarray*}   
| f * \chi _{\nu } (x) | &\leq & || f ||_{\infty , -p } \int _{\mathbb{R }^n } \, |\chi (y) | (1 + |x - \nu ^{-1 } y | )^p \, dy \\   
&\leq & || f ||_{\infty , -p } (1 + |x| )^p \int _{\mathbb{R }^n } \, (1 + |x| )^p |\chi (y) | \, dy .   
\end{eqnarray*}   
Since $f * \chi _{\nu } (x ) \to f(x) $ uniformly on compact  \rouge (N.B. Doesn't $f $ to be continuous here, or does this hold a.e. for locally bounded $f $'s)\noir, it then easily follows that $|| f * \chi _{\nu } - f ||_{\infty , - p - a } \to 0 $ when $a > 0 $ which, using (\ref{eq:s_h_assertion_1}), then implies that $|| \, s_h [f_{\nu } ] - s_h [f ] \, ||_{\infty , -p - a } \to 0 $ if $p + a < \kappa . $  This completes the proof.   
\medskip   
   
\end{versionA}   
   
\subsection{Proof of lemma \ref{lemma:u_h_hat_bis}} It again suffices to consider the case of compactly supported $\widehat{f } $'s. We use the notations of the proof of lemma \ref{lemma:convergence_RBF_bis} above: in particular, let $f_{\varepsilon } $ be the inverse Fourier transform of $\widehat{f } * \chi _{\varepsilon } $ where  $\chi _{\varepsilon } = \varepsilon ^{-n } \chi ( \cdot / \varepsilon ) $ ia an approximation of the identity.   
We have seen that $\Sigma _h (\widehat{f }_{\varepsilon } ) \to \Sigma _h (\widehat{f } ) $ in $\left( C^{K , \lambda } \right) ' $, where $K = \lceil \kappa \rceil - 1 $ and $\lambda = \kappa - K . $ Since $e^{- h^{-2 } t G (h \cdot ) } \in C^{K , \lambda } $ this implies that   
$$   
e^{- h^{-2 } t G (h \cdot ) } \Sigma _h \left( \widehat{f }_{\varepsilon } \right) \to e^{- h^{-2 } t G (h \cdot ) } \Sigma _h (\widehat{f } )   
$$   
in $\left( C^{K , \lambda } \right) '  $ and hence in $\mathcal{S }' (\mathbb{R }^n ) . $   
      
On the other hand, we have seen in the proof of lemma \ref{lemma:convergence_RBF_bis} that $f_{\varepsilon } \to f $ in $L^{\infty } _{- p - a } $ if $a > 0 . $ 
Hence by lemma \ref{lemma:ODE-system}, if $a < \kappa - p $ then $u_h [f_{\varepsilon } ] \to u_h [f ] $ in $L^{\infty } _{- p - a } $ and therefore as tempered distributions. This implies that       
$$   
e^{- h^{-2 } t G (h \cdot ) } \Sigma _h \left( \widehat{f }_{\varepsilon } \right) = \widehat{u_h [f_{\varepsilon } ] } \to \widehat{u_h [f ] } ,   
$$   
where we used lemma \ref{lemma:u_h_hat}. Hence $\widehat{u_h [f ] } = e^{- h^{-2 } t G (h \cdot ) } \Sigma _h (\widehat{f } ) =  e^{- h^{-2 } t G (h \cdot ) } \widehat{s_h [f ] } $ as tempered distributions, as claimed. 
We finally prove (\ref{eq:proof_thm:conv_RBF_scheme}): 
if we let   
$$   
g(\xi , t ; h ) := e^{- t (h^{-2 } G(h \xi ) - |\xi |^2 ) } - 1 ,   
$$   
then $g $ is a $C^{K , \lambda } $ 
-function of $\xi $ and   
$$   
\widehat{u_h [f ] } (\cdot , t ) - \widehat{u } (\cdot , t ) = e^{- t h^{-2 } G(h \cdot ) } \left( \, \widehat{s_h [f ] } - \widehat{f } \, \right) + (g(\cdot , t ; h ) - 1 ) \, e^{ - t |\cdot |^2 } \widehat{f } .   
$$   
Since $g (\xi , t ; h ) $ vanishes to order $|\xi |^{\kappa } $ in $\xi = 0 $, by proposition \ref{prop_G}(ii), the representation $\widehat{f } = \Lambda _{\widehat{f } }+ \sum _{|\alpha | \leq \lceil \kappa \rceil - 1 } c_{\alpha } \delta ^{(\alpha ) } $ from the proof of lemma \ref{lemma:convergence_RBF_bis} shows that the distribution $g (\cdot , t , h ) \widehat{f } $ can be identified with the locally integrable function $\xi \to g(\xi , t , h ) \widehat{f }(\xi ) . $ 
\hfill $\Box $   
   
\begin{versionA}   
   
\section{\bf Proof of lemma \ref{lemma:Stein}}   
   
\textcolor{blue}{This appendix contains alternative proofs of the Fourier-transform estimates from Stein's book we used: since these results are anyhow known, this is not to be retained for the final version, but for own information and re-assurance only}   
\medskip   
   
\noindent {\it Proof}. Stein \cite{Stein} uses the Paley-Littlewood decomposition, and actually proves a stronger result: if $|\partial _{\xi } ^{\alpha } a (\xi ) | \leq C_{\alpha } |\xi |^{p - |\alpha | } $ on $\mathbb{R }^n \setminus 0 $ for $|\alpha | \leq \lfloor p \rfloor + n + 1 $, then the inverse Fourier transform $k := \mathcal{F }^{-1 } (a) $ (in distribution sense) can be identified with a continuous function $k(x) $ on $\mathbb{R }^n \setminus 0 $ which satisfies the estimate $|k(x) | \leq C |x |^{- p - n } $ there. (Stein states and proves this for $p = 0 $, which is the case of greatest interest for classical harmonic analysis\footnote{when $a $ is the symbol of a singular integral operator}, but the proof trivially generalises to any $p > -n . $) For $a $'s with compact support $k $ will automatically be a bounded function and only the estimate at infinity matters.   

Alternatively, one can give a demonstration of the lemma using integration by parts, as follows. Let $a $ be as in the lemma, and let $\chi \in C^{\infty }(\mathbb{R }^n ) $ with all derivatives of $\chi $ bounded such that $\chi (\xi ) = 0 $ for $|\xi | \leq 1 $, $\chi (\xi ) = 1 $ for $|\xi | \geq 2 . $ Let $x \in \mathbb{R }^n $, $x \neq 0 . $ Then   
$$   
k(x) = (2 \pi )^{-n } \int _{\mathbb{R }^n } \, \chi (|x| \xi ) a(\xi ) e^{i (x , \xi ) } \, d\xi + (2 \pi )^{-n } \int _{\mathbb{R }^n } \, (1 - \chi (|x| \xi ) ) a(\xi ) e^{i(x , \xi ) } \, d\xi .   
$$   
Since $1 - \chi (|x| \xi ) $ is supported in $|\xi | \leq 2 |x|^{-1 } $, the second integral can be bounded by a constant times $\int _0 ^{2|x|^{-1 } } r^{p + n - 1 } = C |x|^{- p - n } $, using (\ref{eq:Stein1}) with $\alpha = 0 . $ (Recall that $p > - n . $). To estimate the first integral, let $\alpha $ be a multi-index of length $|\alpha | =  \lfloor p \rfloor + n + 1 $ such that $x^{\alpha } \neq  0 . $ An integration by parts starting off from $i^{- |\alpha | } x^{- \alpha } \partial _{\xi } ^{\alpha } e^{i(x , \xi ) } = e^{i (x , \xi ) } $ then shows the integral to be equal to   
\begin{equation} \label{eq:Stein_proof_1}   
i^{|\alpha | } x^{- \alpha } \int _{\mathbb{R }^n } \, \partial _{\xi } ^{\alpha } \left( \chi (|x| \xi ) a(\xi ) \right) e^{i(x , \xi ) } \, d\xi .   
\end{equation}   
Now   
\begin{eqnarray*}   
\left | \partial _{\xi } ^{\alpha } \left( \chi (|x| \xi ) a(\xi ) \right) \right | &=& \left | \sum _{\beta \leq \alpha } \left( \begin{array}{cc} \alpha \\ \beta \end{array} \right) |x|^{|\beta | } \partial _{\xi }^{\beta } \chi (|x| \xi ) \partial ^{\alpha - \beta } _{\xi } a (\xi ) \right | \\   
&\leq & C \sum _{\beta \leq \alpha } |x|^{|\beta | } (\partial _{\xi } ^{\beta } \chi )(|x| \xi ) |\xi | ^{p - |\alpha | + |\beta | } .   
\end{eqnarray*}   
The terms with $\beta \neq 0 $ are supported in $|\xi | \sim |x|^{-1 } $ and will, after integration with respect to $\xi $ and multiplication by $|x^{- \alpha } | $, make a contribution to (\ref{eq:Stein_proof_1}) of order at most $|x |^{-|\alpha | + |\beta | } \cdot |x|^{- p + |\alpha | - |\beta | } \cdot |x|^{-n } = |x|^{- p - n } $ to $k(x) . $ As for the term with $\beta = 0 $, we note that   
\begin{eqnarray*}   
|x^{- \alpha } | \int _{\mathbb{R }^n } \chi (|x| \xi ) |\xi |^{p - |\alpha | } d\xi &\leq & C |x^{- \alpha } | \int _{|x|^{-1 } } ^{\infty } r^{p - |\alpha | + n - 1 } dr \\   
&\simeq & |x|^{- p - n } ,   
\end{eqnarray*}   
where we used that $p - |\alpha | + n < 0 . $ This estimate will in fact hold for all $x' $ such that $\frac{1 }{2 } |x' | ^{\alpha } \leq |x^{\alpha } | $, 
say, and a standard partition of unity on the unit sphere then shows that $k(x) | \leq C |x|^{- p - n } . $ Since $k(x) $ is trivially bounded, we only keep this estimate for $|x| $, and the lemma follows.     
\medskip   
   
\noindent \textcolor{blue}{{\bf Note for self}: if $p \geq 0 $, then need to take (at least) $|\alpha | = \lfloor p \rfloor + n + 1 $; if $p < 0 $, $p = - q $, say, then need   
$$   
|\alpha | > p + n = n - q .   
$$   
Note that $q < n $ since we assumed that $p > - n . $ So need   
$$   
|\alpha | = \lfloor n - q \rfloor + 1 = n - \lceil q \rceil + 1 = n + \lfloor p \rfloor + 1 .   
$$   
(E.g. $\lfloor - \frac{1 }{2 } \rfloor =  - 1 = - \lceil \frac{1 }{2 } \rceil . $) }   
\medskip   
      
\begin{remarks} \rm{(i) More generally, if $a $ is $C^{\lfloor p \rfloor  + n + K + 1 } $ on $\mathbb{R }^n \setminus 0 $ and (\ref{eq:Stein1}) holds for all $|\alpha | \leq \lfloor p \rfloor + n + K + 1 $, then $|\partial _x ^{\beta } k (x) | \leq (1 + |x|)^{- p - n - |\beta | } $ for all $|\beta | \leq K . $   
   
\noindent (N. B. $k $ is trivially $C^{\infty } $, even analytic,  since $a $ has compact support; what matters here is its decay, and that of some of its derivatives, at infinity.)   
\medskip   
   
\noindent (ii) \rouge (probably not to be included in final version - just included for own satisfaction) \noir The proof of the lemma can easily be modified so as to include non-compactly supported $a $'s satisfying (\ref{eq:Stein1}). If $f \in \mathcal{S } (\mathbb{R }^n ) $ is a rapidly decreasing function supported away from 0, then   
\begin{eqnarray*}   
\langle k, f \rangle &=& \langle a , \mathcal{F }^{-1 } (f) \rangle \\   
&=& (2 \pi )^{-n } \int _{\mathbb{R }^n } \int _{\mathbb{R }^n } a(\xi ) f(x) e^{i (x, \xi ) } \, dx \, d\xi \\   
&=& \lim _{\varepsilon \to 0 } (2 \pi )^{-n } \int _{\mathbb{R }^n } \int _{\mathbb{R }^n } a(\xi ) f(x) e^{- \varepsilon |\xi |^2 } e^{i (x, \xi ) } \, dx \, d\xi ,   
\end{eqnarray*}   
where we introduce the $e^{- \varepsilon |\xi^2 } $ so as to have an absolutely convergent "double" integral. Let $\chi := \chi _1 $ be as in the proof of the lemma, and put $\chi _2 := 1 - \chi _1 . $ Then we can write   
\begin{eqnarray*}   
(2 \pi )^n \langle k , f \rangle &=& \lim _{\varepsilon \to 0 } \int _{\mathbb{R }^n } \int _{\mathbb{R }^n } a(\xi ) \left( \chi _1 (|x| \xi ) + \chi _2 (|x| \xi ) \right) f(x) e^{- \varepsilon |\xi |^2 } e^{i (x, \xi ) } \, dx \, d\xi \\   
&=& \lim _{\varepsilon \to 0 } \int _{\mathbb{R }^n } \int _{\mathbb{R }^n } (i x )^{- \alpha } \partial _{\xi } ^{\alpha } \left( a(\xi ) \chi _1 (|x| \xi ) e^{- \varepsilon |\xi |^2 } \right) f(x) e^{i (x, \xi ) } \, dx \, d\xi \\   
&\ & +  \int _{\mathbb{R }^n } \int _{\mathbb{R }^n } a(\xi ) \chi _2 (|x| \xi ) f(x) e^{i (x, \xi ) } \, dx \, d\xi \\   
&=& \int _{\mathbb{R }^n } \, (k_1 (x) + k_2 (x) ) f(x) \, dx ,   
\end{eqnarray*}   
where   
$$   
k_1 (x) := (i x )^{- \alpha } \int _{\mathbb{R }^n } \partial _{\xi } ^{\alpha } \left( a(\xi ) \chi _1 (|x| \xi ) \right) e^{i (x, \xi ) } \, d\xi , \ \ k_2 (x) := \int _{\mathbb{R }^n } a(\xi ) \chi _2 (|x| \xi ) e^{i (x, \xi ) } \, d\xi ,    
$$   
and where we reduce the support of $f $ to lie in a cone on which $x^{\alpha } \neq 0 $ for a given $\alpha $ of length $\lfloor p \rfloor + n + 1 . $ Taking into account the supports of $\chi _1 $ and $\chi _2 $, straightforward estimates show as before that $k_1 (x) $ and $k_2 (x) $ are both $O(|x|^{- p - n } ) $, at least on the cone where $|x^{\alpha } | \geq \frac{1 }{2 } |x |^{|\alpha | } $ for our given $\alpha . $ This shows, modulo a partition of unity argument, that $k $ can be identified with a continuous function on $\mathbb{R }^n \setminus 0 $ satisfying the required estimate.   
\bigskip   
   
The restriction to compactly supported $a(\xi ) $ is not essential: if $a \in C^{\lfloor p \rfloor + n + 1 } (\mathbb{R }^n \setminus 0 ) $ satisfies (\ref{eq:Stein1}), then $k $, as a tempered distribution, can be identified with a continuous function $k(x) $ on $\mathbb{R }^n \setminus 0 $ satisfying      
\begin{equation} \label{eq:Stein2}   
|k(x) | \leq C |x|^{ - p - n } , \ \ x \neq 0 ,     
\end{equation}   
and analogously for derivatives $\partial _x ^{\beta } k (x) $ (with $|x |^{-p - n } $ replaced by $|x |^{- p - n - |\beta | } $), assuming (\ref{eq:Stein1}) holds to order $\lfloor p \rfloor + n + |\beta | + 1 . $ The estimate (\ref{eq:Stein2}) tells us something about the behaviour of $k(x) $ at both zero and infinity. For the proof of theorem \ref{thm:L_wBuhmann} we will only need the information at infinity. 
}   
\end{remarks}

%
%
%
   

\noindent \textcolor{blue}{{\bf Proof of (\ref{eq:S^p}) - not to be included in final version.} The proof in Stein's book uses a Paley-Wiener decomposition. One can give the following simpler proof of the result we stated: we first note that if $N = \lfloor p \rfloor + n + 1 $, then $\widehat{f } $ can be indentified on $\mathbb{R } \setminus 0 $ with the continuous function (also denoted by $\widehat{f } $, by abuse of notation):   
\begin{equation} \label{eq:example}   
\widehat{f } (\xi ) := |\xi |^{-2 N } \int _{\mathbb{R }^n } \, e^{- i (x, \xi ) } \, (\xi , \partial _x )^N \left [f (x) \right ] \, dx ,  
\end{equation}   
the integral being absolutely convergent since $N > p + n $ ($N $ is the smallest integer having this property). Indeed, if $\chi \in C^{\infty } _c (\mathbb{R }^n ) $ (or even $\mathcal{S }(\mathbb{R }^n ) $) such that $\chi (0 ) = 0 $, then $\chi (\varepsilon x ) f(x) \to f(x) $ as tempered distributions, and it follows that   
$$   
\langle \widehat{f } , \psi \rangle = \lim _{\varepsilon \to 0 } \int _{\mathbb{R }^n } \int _{\mathbb{R }^n } \chi (\varepsilon x ) f(x) \psi (\xi ) e^{- i (x, \xi ) } \, dx d\xi ,   
$$   
each of the integrals now being absolutely convergent. Since $(\xi , \partial _x ) e^{- i (x, \xi ) } = - i |\xi |^2 e^{- i (x, \xi ) } $, an $N $-fold integration by parts shows that if ${\rm supp }(\psi ) \subset \mathbb{R }^n \setminus 0 $, then   
$$   
\langle \widehat{f } , \psi \rangle = \lim _{\varepsilon \to 0 } \int _{\mathbb{R }^n } \int _{\mathbb{R }^n } \, (\xi , \partial _x )^N \left [ \chi (\varepsilon x ) f(x) \right ] \frac{\psi (\xi ) }{|\xi |^{2N } } e^{- i (x, \xi ) } \, dx d\xi .   
$$   
When the differential operator hits $\chi (\varepsilon x ) $ we will obtain terms such as $\varepsilon ^{|\alpha | } \chi ^{(\alpha ) } (\varepsilon x ) f^{(\beta ) } (x) $ with $|\alpha | + |\beta | = N $,  and $\chi ^{(\alpha ) } (\varepsilon x ) f^{(\beta ) } (x) \to \chi ^{(\alpha ) } (0) f^{(\beta ) } (x) $ as tempered distributions, and to 0 after multiplication by $\varepsilon ^{|\alpha | } $, and (\ref{eq:example}) follows. It then follows immediately that $|\widehat{f } (\xi ) | \ \leq C |\xi |^{- N } $ on $\mathbb{R }^N . $ To obtain the improved estimate at $\xi = 0 $, suppose that $\chi = 1 $ on $B(0, 1 ) = \{ \xi : |\xi \ = 1 \} $ and 0 outside of $B(0, 2 ) $, and split the integral as     
\begin{eqnarray*}   
\widehat{f } (\xi ) &=& |\xi |^{-2N } \int _{\mathbb{R }^n } \, e^{- i (x, \xi ) } \, \chi (|x | |\xi | ) \, (\xi , \partial _x )^N f (x) \, dx + \\   
&& |\xi |^{-2N } \int _{\mathbb{R }^n }  \, e^{- i (x, \xi ) } (1 - \chi (|x | |\xi | ) )  \, (\xi , \partial _x )^N f (x)  \, dx .   
\end{eqnarray*}   
The second integral extends over $|x | \geq 1 / |\xi $ and can be bounded by a constant times   
$$   
|\xi |^{-N } \int _{|\xi |^{-1 } } ^{\infty } r^{p - N + n - 1 } = \frac{|\xi |^{p - n } }{p + n - N } .   
$$   
As for the 
first integral, after integrating $N $ times by parts in the other direction we end up with an expression which can be bounded by a constant times a sum of terms of the form   
\begin{equation} \nonumber   
\int _{\mathbb{R }^n }  | \chi ^{(\alpha ) } (|x | |\xi | ) | \, | f(x) | \, dx , \ \ |\alpha | + |\beta | = N .   
\end{equation}   
which for small $|\xi | $ can be bounded by a constant times   
$$   
\int _{|x | \leq 2 |\xi |^{-1 } } (1 + |x |)^p \, dx \sim |\xi |^{-n - p } .   
$$   
This completes the proof of $|\widehat{f } (\xi ) | \leq C |\xi |^{- n - p } . $ 
}   
   
\end{document}